\documentclass{article}
\pdfoutput=1
\makeatletter
\def\ps@pprintTitle{%
	\let\@oddhead\@empty
	\let\@evenhead\@empty
	\def\@oddfoot{\reset@font\hfil\thepage\hfil}
	\let\@evenfoot\@oddfoot
}
\makeatother
\usepackage{graphicx}
\usepackage{bm}
\usepackage{pstricks}
\usepackage{caption}
\usepackage{subcaption}
\usepackage{rotating}
\usepackage{inputenc}
\usepackage{wrapfig}
\usepackage{framed}
\usepackage{url}
\usepackage{hyperref}
\usepackage{amssymb,latexsym,amsmath,amsthm,verbatim}
\usepackage{graphicx,epsfig,epstopdf,amssymb,color}
\usepackage{lineno}

\newcommand{\w}{\omega}

\newcommand{\tB}{\tilde{B}}
\newcommand{\tW}{\tilde{W}}
\newcommand{\tX}{\tilde{X}}

\newcommand{\ta}{\tilde{a}}
\newcommand{\A}{\mathcal{A}}
\newcommand{\B}{\mathcal{B}}
\newcommand{\C}{\mathcal{C}}
\newcommand{\D}{\mathcal{D}}
\newcommand{\E}{\mathcal{E}}
\newcommand{\I}{\mathcal{I}}
\newcommand{\J}{\mathcal{J}}
\newcommand{\HH}{\mathcal{H}}
\newcommand{\M}{\mathcal{M}}
\newcommand{\OO}{\mathcal{O}}
\newcommand{\R}{\mathcal{R}}
\newcommand{\SSS}{\mathcal{S}}
\newcommand{\T}{\mathcal{T}}
\newcommand{\W}{\mathcal{W}}
\newcommand{\eps}{\varepsilon}
\newcommand{\beq}{\begin{eqnarray}}
\newcommand{\feq}{\end{eqnarray}}

\usepackage[affil-it]{authblk}
\newtheorem{theorem}{Theorem}[section]
\newtheorem{lemma}[theorem]{Lemma}

\newtheorem{corollary}[theorem]{Corollary}
\newtheorem{remark}[theorem]{Remark}

\numberwithin{equation}{section}

\title{The existence of localized vegetation patterns in a systematically reduced model for
dryland vegetation}

\author[1,*]{Olfa Jaibi}
\author[1]{Arjen Doelman}
\author[1]{Martina Chirilus-Bruckner}
\author[2,3]{Ehud Meron}
\affil[1]{Mathematisch Instituut - Universiteit Leiden, P.O. Box 9512, 2300 RA, Leiden, The Netherlands}
\affil[2]{The Blaustein Institutes for Desert Research, Ben-Gurion University, Sede Boqer Campus 8499000, Israel}
\affil[3]{Department of Physics, Ben-Gurion University, Beer-Sheva, Israel}
\affil[*]{Corresponding author: \textsf{\textnormal{o.jaibi@math.leidenuniv.nl}}}

\begin{document}
\maketitle

\begin{abstract}
In this paper we consider the 2-component reaction-diffusion model that was recently obtained by a systematic reduction of the 3-component Gilad et al. model for dryland ecosystem dynamics \cite{GiladEtAl2004}. The nonlinear structure of this model is more involved than other more conceptual models, such as the extended Klausmeier model, and the analysis a priori is more complicated. However, the present model has a strong advantage over these more conceptual models in that it can be more directly linked to ecological mechanisms and observations. Moreover, we find that the model exhibits a richness of analytically tractable patterns that exceeds that of Klausmeier-type models. Our study focuses on the 4-dimensional dynamical system associated with the reaction-diffusion model by considering traveling waves in 1 spatial dimension. We use the methods of geometric singular perturbation theory to establish the existence of a multitude of heteroclinic/homoclinic/periodic orbits that `jump' between (normally hyperbolic) slow manifolds, representing various kinds of localized vegetation patterns. The basic 1-front invasion patterns and 2-front spot/gap patterns that form the starting point of our analysis have a direct ecological interpretation and appear naturally in simulations of the model. By exploiting the rich nonlinear structure of the model, we construct many multi-front patterns that are novel, both from the ecological and the mathematical point of view. In fact, we argue that these orbits/patterns are not specific for the model considered here, but will also occur in a much more general (singularly perturbed reaction-diffusion) setting. We conclude with a discussion of the ecological and mathematical implications of our findings.\\
\begin{tabbing}
\>singular perturbations
\end{tabbing}

\end{abstract}

\section{Introduction}
\label{sec:Intro}

Ecosystems consist of organisms that interact among themselves and with their environment. These interactions involve various kinds of feedback processes that may combine to form positive feedback loops and instabilities when environmental conditions change \cite{Meron2018ann_rev,Meron2019}. In many ecosystems -- drylands, peatlands, savannas, mussel beds, coral reefs, and ribbon forests –- the leading feedback processes have different spatial scales: a short-range facilitation by local modification of the environment versus a long-range competition for resources \cite{RietkerkKoppel2008}. Like the well-established activator-inhibitor principle in bio-chemical systems \cite{Murray2013}, the combination of these scale-dependent feedback mechanisms can induce instabilities that result in large-scale spatial patterns, which are similar to a wide variety of vegetation patterns observed in drylands, peatlands, savannas and undersea \cite{Deblauwe2008geb,ShefferEtAl2013,BastiaansenEtAl2018,GowdaEtAl2018,Ruiz2017sci_adv,RietkerkEtAl2004,Getzin2016}. Varying climatic conditions and human disturbances may continue to propel ecosystem dynamics. Ecosystem response to decreasing rainfall, for example, may take the form of abrupt collapse to a nonproductive `desert state' \cite{Scheffer2001nature,vHardenbergEtAl2001,RietkerkEtAl2004}, or involve gradual desertification, consisting of a cascade of state transitions to sparser vegetation \cite{Siteur2014eco_comp,BastiaansenEtAl2020}, or gradual vegetation retreat by front propagation \cite{BelEtAl12,Zelnik2018eco_ind}. Understanding the dynamics of spatially extended ecosystems has become an active field of research in the last two decades -- within communities of ecologists, environmental scientists, mathematicians and physicists. Apart from its obvious environmental and societal relevance, the phenomena exhibited pose fundamental challenges to the research field of pattern formation.
\\ \\
Several models of increasing complexity have been proposed in the past two decades. Of these, the models that have received most attention are the one-component model by Lefever and Lejeune \cite{Lefever1997bmb}, the two-component models by Klausmeier \cite{Klausmeier1999} and von Hardeberg et al. \cite{vHardenbergEtAl2001}, and the three-component models by Rietkerk et al. \cite{RietkerkEtAl2002} and Gilad et al. \cite{GiladEtAl2004}. A basic difference between these models is the manner by which they describe water dynamics. The Lefever-Lejeune model does not describe water dynamics at all, the Klausmeier model does describe water dynamics but does not make a clear distinction between soil water and surface water \cite{vdSteltetal13}, while the von Hardenberg et al. model only takes soil water into account. The Rietkerk model and the Gilad et al. model describe both soil water and surface water dynamics and, therefore, capture more aspects of real dryland ecosystems. A major difference between these two models is the inclusion of water conduction by laterally spread roots, as an additional water-transport mechanism, in the Gilad et al. model.
\\ \\
Despite these differences, all models appear to share a similar bifurcation structure, as analytical and numerical-continuation studies reveal \cite{Lejeune2004ijqc,GowdaEtAl2014,Dijkstra2011IJBC,Zelnik2013philtransA}, except the Klausmeier model. This structure includes, in particular, a stationary uniform instability (i.e. involving the monotonic growth of spatially uniform perturbations) of the bare soil (zero biomass) state as the precipitation rate exceeds a threshold value. The Klasumeier model fails to capture that instability, leaving the bare soil state stable at all precipitation values. This behavior limits the applicability of the Klausmeier model to ecological contexts where the bare soil state is stabilized at relatively high precipitation rates, e.g. by high evaporation rates. Nevertheless, of all models, the Klausmeier model and its extension to include water diffusion have been studied to a greater extent \cite{CarterD18, Sherratt2010, Sherratt2013,vdSteltetal13}, partly because the extended form coincides with the much studied Gray-Scott model for autocatalytic chemical reactions -- see \cite{BastiaansenD19,KolokolnikovEtAl2015,SewaltD17} and the references therein.
\\ \\
All models have been analyzed mathematically to various extents. Two main analytical approaches can be distinguished in these studies (see however Goto et al. \cite{GotoEtAl2011}); linear stability and weakly nonlinear analysis near instability points \cite{Lejeune2004ijqc,DawesWilliams2016,GowdaEtAl2014,GowdaEtAl2016,Fernandez-Oto2019prl,vdSteltetal13}, and singular perturbation analysis, based on the disparate length scales associated with biomass (short) and water (long) \cite{CarterD18,BastiaansenD19,KolokolnikovEtAl2015,SewaltD17}. Studies of the first category are strictly valid only near instability points, although they do capture essential parts of the bifurcation structure even far from these points and are quite insightful in this respect. By contrast, studies of the second category apply to the strongly nonlinear `far-from-equilibrium' regime, where desertification transitions take place, and are, potentially, of higher ecological interest. So far, however, these studies have been limited to the simpler and less realistic Klausmeier model.
\\ \\
In this paper we apply a geometric singular perturbation analysis to a reduced version of the Gilad et al. model in order to study the existence of various forms of localized patterns. Singular perturbation theory has already been applied to three-component models -- see for instance \cite{DVeerman2015,vHeijsterSandstede2014} -- and could be applied, in principle, to the non-local three-component Gilad et al. model. Here we choose to consider ecological contexts that allow to reduce that model to a local two-component model for the vegetation biomass and the soil water content. Specifically, we assume soil types characterized by high infiltration rates of surface water into the soil, such as sandy soil, and plant species with laterally confined root zones (see \ref{Appendix:derivation} for more details). These conditions are met, for example, by Namibian grasslands showing localized and extended gap patterns (`fairy circles') \cite{Zelnik2015}. We further simplify the problem by assuming one space dimension. The reduced model reads:
\begin{eqnarray}
\label{eq:meron}
\left\{	
\begin{aligned}
\frac{\partial \tB}{\partial T} &= \Lambda \tW \tB(1-\tB/K)(1+E\tB) - M\tB + D_B \frac{\partial^2 \tB}{\partial \tX^2}, \\
\frac{\partial \tW}{\partial T} &= P - N(1-R\tB/K)\tW - \Gamma \tW\tB(1+E\tB) + D_W \frac{\partial^2 \tW}{\partial \tX^2},
\end{aligned}
\right.
\end{eqnarray}
where $\tB(\tX,T) \geq 0$ and $\tW(\tX,T) \geq 0$ represent areal densities of biomass and soil water, respectively, and $\tX \in \mathbb{R}$,  $T \in \mathbb{R}^+$ are the space and time coordinates.
In the biomass ($\tB$) equation, $\Lambda$ represents the biomass growth rate coefficient, $K$ the maximal standing biomass, $E$ is a measure for the root-to-shoot ratio, $M$ the plant mortality rate and $D_B$ the seed-dispersal or clonal growth rate, while in the water ($\tW$) equation, $P$ represents the precipitation rate, $N$ the evaporation rate, $R$ the reduction of the evaporation rate due to shading, $\Gamma$ the water-uptake rate coefficient and $D_W$ the effective soil water diffusivity. Notice that the power of the factor $(1+E\tB)$ in both equations is unity, whereas in the reduced model in \cite{Zelnik2015} the power is two. This difference stems from the consideration in this study of one space dimension rather than two (see \ref{Appendix:derivation}).
\\ \\
From the ecological point of view, the advantage in studying model (\ref{eq:meron}) over the much analyzed Klausmeier model lies in the fact that it has been systematically derived from a more extended model that better captures relevant ecological processes, such as water uptake by plant roots (controlled by $E$), reduced evaporation by shading (controlled by $R$), and late-growth constraints, such as self-shading (controlled by $K$) -- see  \cite{GiladEtAl2004,GiladEtAl2007,Meron2015,ShefferEtAl2013}. As a consequence, (mathematical) insights in (\ref{eq:meron}) can be linked to ecological observations and mechanisms in a  direct fashion. Naturally, there also is a disadvantage to analyzing a model that incorporates concrete ecological mechanisms: the more involved --  algebraically more complex -- nonlinear structure of (\ref{eq:meron}) a priori makes it less suitable for an analytical study than the Klausmeier model (or other more conceptual models). However, that apparent disadvantage turned around into an advantage: we will find that the reduced model transcends by far the Klausmeier model in terms of richness of analytically tractable pattern solutions.
\\ \\
The model equations (\ref{eq:meron}) represent a singularly perturbed system, because of the low seed-dispersal rate as compared with soil water diffusion, that is, $D_B \ll D_W$ \cite{GiladEtAl2007,vdSteltetal13,Zelnik2015}. To make this explicit and to simplify (\ref{eq:meron}) as much as possible, we introduce the following scalings,
\begin{eqnarray}
B = \frac{\tB}{\alpha}, \;\; W = \frac{\tW}{\beta}, \;\; t = \delta T, \;\; x =\gamma \tX,
\label{eq:biomasswater}
\end{eqnarray}
and set,
\begin{eqnarray}
\alpha = K - \frac{1}{E}, \; \;
\beta = \frac{MK}{\alpha^2 \Lambda E}, \; \;
\gamma = \sqrt{\frac{\alpha^2 \beta \Lambda E}{K D_B}}, \; \;
\delta = \frac{\alpha^2 \beta \Lambda E}{K}.
\label{scalingpars}
\end{eqnarray}
By also introducing our main parameters,
\begin{eqnarray}
a  = \frac{KE}{(KE -1)^2}, \; \; \varepsilon^2 = \frac{D_B}{D_W} \ll 1,
\label{def:aeps}
\end{eqnarray}
we arrive at,
\begin{eqnarray}
\left\{	
\begin{aligned}
\label{eq:scaled}
B_t &= \left(aW- 1\right)B + WB^2 - WB^3 +  B_{xx}, \\
W_t &=\Psi -\left[\Phi + \Omega B + \Theta B^2\right]W + \frac{1}{\varepsilon^2}W_{xx},
\end{aligned}
\right.
\end{eqnarray}
in which,
\begin{eqnarray}
\Psi = \frac{\alpha^2P\Lambda E}{M^2K }, \; \;
\Phi = \frac{N }{M}, \; \;
\Omega= \frac{\alpha}{M}\left(\Gamma - \frac{R}{K}\right), \; \;
\Theta = \frac{\alpha^2 \Gamma E}{M}.
\label{def:pars}
\end{eqnarray}
A more detailed derivation of the scaled equations (\ref{eq:scaled}) from (\ref{eq:meron}) is given in \ref{Appendix:scaling}. Since the signs of the parameters in (\ref{eq:scaled}) will turn out to be crucial in the upcoming analysis, we note explicitly that $a, \Psi, \Phi, \Theta \geq 0$ while $\Omega \in \mathbb{R}$, i.e. $\Omega$ may be negative.
\\
\begin{figure}[t]
    \centering
     \includegraphics[width=\textwidth]{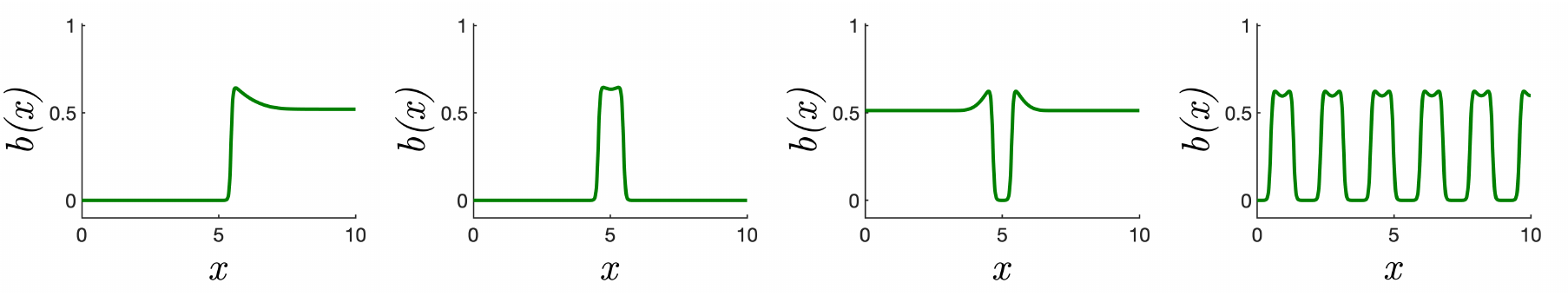}
    \caption{The 4 basic patterns exhibited by numerical simulations of model (\ref{eq:scaled}): a traveling (heteroclinic) invasion front (Theorem \ref{T:primaryfronts}), a stationary, homoclinic, 2-front vegetation spot (Theorem \ref{T:standing2fspots}), a stationary homoclinic, 2-front vegetation gap (Theorem \ref{T:standing2fgaps}), and a stationary, spatially periodic multi-front spot/gap pattern (Theorem \ref{T:periodics}) -- see Remark \ref{rem:parsFig1} for the precise parameter values.}
    \label{fig:basicpatterns}
\end{figure}
\\
We study pattern formation in (\ref{eq:meron}) by analyzing (\ref{eq:scaled}) using the methods of (geometric) singular perturbation theory \cite{Jones1995,Kaper1999} and thus `exploit' the fact that $\eps \ll 1$ (\ref{def:aeps}). In fact -- apart from some observations in section \ref{sec:crithom} and the discussion section \ref{sec:disc} -- we focus completely on the `spatial' 4-dimensional dynamical system that is obtained from (\ref{eq:scaled}) by considering `simple' solutions that are stationary in a co-moving frame traveling with constant speed $c$. More specifically, in this paper we study the existence of traveling (and stationary) solutions -- in particular localized (multi-)front solutions connecting a (uniform) bare soil state to a uniform vegetation state, or a bare soil state to itself (with a `passage' along a vegetated state), etc. -- by taking the classical approach of introducing a (uniformly) traveling coordinate $\xi = x - ct$, with speed $c \in \mathbb{R}$ an a priori free ${\cal O}(1)$ parameter (w.r.t. the asymptotically small parameter $\eps$). By setting $(B(x,t),W(x,t)) = (b(\xi), w(\xi))$ and  introducing $p = b_{\xi}$ and $q = \frac{1}{\varepsilon}w_{\xi}$, PDE (\ref{eq:scaled}) reduces to
\begin{eqnarray}
\left\{	
\begin{aligned}
b_{\xi}\;\; &= \;\;p,\\
p_{\xi}\;\; &= \;\;wb^3 - wb^2 + (1-aw)b - cp,\\
w_{\xi} \;\;&=\;\;\varepsilon q,\\
q_{\xi}\;\;&=\;\;\varepsilon\left(-\Psi +\left[\Phi + \Omega b + \Theta b^2\right]w \right) - \eps^2 cq.
\end{aligned}
\right.
\label{eq:fast}
\end{eqnarray}
Fig. \ref{fig:basicpatterns} shows four basic patterns that naturally appear in simulations of (\ref{eq:scaled}) and have identifiable ecological counterparts: vegetation fronts (ecotones), isolated vegetation spots and gaps, and periodic patterns \cite{Meron2007chaos,Fernandez-Oto2014phil_tran_A,Escaff2015pre,Getzin2016}. These patterns are rigorously constructed by the methods of singular perturbation theory in section \ref{sec:locfront}. From the geometrical point of view, these constructions are natural and thus relatively straightforward: all patterns in Fig. \ref{fig:basicpatterns} `jump' between two well-defined slow manifolds (of (\ref{eq:fast})) -- see Theorems  \ref{T:primaryfronts}, \ref{T:standing2fspots}, \ref{T:standing2fgaps}, and \ref{T:periodics}. Therefore, the main work in establishing these results lies in resolving technical issues (which can be done by the preparations of section \ref{sec:setup}). However, the preparations of section \ref{sec:setup} also form the origin of the construction of a surprisingly rich `space' of  traveling and/or stationary patterns that goes way beyond those exhibited in Fig. \ref{fig:basicpatterns} -- see for instance the sketches of Fig. \ref{fig:sketchessection1}. These are novel patterns, at least from the point of view of explicit rigorous mathematical constructions in multi-component reaction-diffusion equations. However, similar patterns have been analyzed as (perturbations of) heteroclinic networks in a more abstract framework -- see \cite{Radem2005,Radem2010} and the references therein. Moreover, patterns similar to those of Fig. \ref{fig:sketchessection1} have been observed in simulations of the Klausmeier-Gray-Scott model \cite{ZelnikEtAl2018}, although with parameter settings beyond that for which the mathematical singular-perturbation approach can be applied.
\\ \\
Here, our motivation to study these patterns is primarily ecological; however, we claim that patterns like these must also occur generically in the setting of a completely general class of singularly perturbed 2-component reaction-diffusion systems -- as we will motivate in more detail in section \ref{sec:disc}. Thus, our explicit analysis of model (\ref{eq:scaled}) provides novel mathematical insights beyond that of the present ecological setting. The driving mechanism behind these patterns originates from the perturbed integrable flow on the slow manifolds associated with (\ref{eq:fast}) -- see sections \ref{sec:slowred} and \ref{sec:slowfull}. The perturbation terms are generically introduced by the $\OO(\eps)$ differences between the slow manifold and its $\eps \to 0$ limit, and they transform the (Hamiltonian) integrable reduced slow flow to a (planar) `nonlinear oscillator with nonlinear friction' that can be studied by explicit Melnikov methods. Typically, one for instance expects (and finds: Theorem \ref{T:slowflow}) persistent periodic solutions on the slow manifold. Associated with these persisting periodic solutions, one can subsequently construct heteroclinic 1-front connections between a critical point -- representing the uniform bare soil state in the ecological setting -- and such a periodic pattern (Theorems \ref{T:frontstoperiodic} and \ref{T:statfrontstoperiodic} and Fig. \ref{fig:sketchessection1}b) and a countable family of `higher order' heteroclinic 1-fronts between critical points that limits on these orbits (Theorem \ref{T:count1fronts} and Fig. \ref{fig:sketchessection1}a -- where we note that Fig. \ref{fig:basicpatterns}a represents the very first -- primary -- member of this family). In the case of (stationary) localized spot patterns, one can construct a countable family of connections that follow the periodic orbit for arbitrarily many `spatial oscillations' (Theorem \ref{T:higherorder2fronts} and Fig. \ref{fig:sketchessection1}c, \ref{fig:sketchessection1}d). Combining these insights with the ideas of \cite{DKvdPloeg01}, one may even construct many increasingly complex families of spatially periodic and aperiodic multi-spot/gap patterns (Corollary \ref{cor:higherorderperiodics} and section \ref{sec:perpatt}). Moreover, we can explicitly study the associated bifurcation scenarios: in section \ref{sec:frontsbyM+eps} we present a scenario of cascading saddle-node bifurcations of heteroclinic 1-front connections starting from no such orbits to countably many -- all traveling with different speed (Theorem \ref{T:count1fronts} and Figs. \ref{fig:basicpatterns}a, \ref{fig:sketchessection1}a and \ref{fig:sketchessection1}b) -- back to 1 unique 1-front pattern (of the type presented in Fig. \ref{fig:basicpatterns}a) -- see Fig. \ref{fig:cPhi2} in section \ref{sec:frontsbyM+eps}.
\\
\begin{figure}[t]
    \centering
    \includegraphics[width=\textwidth]{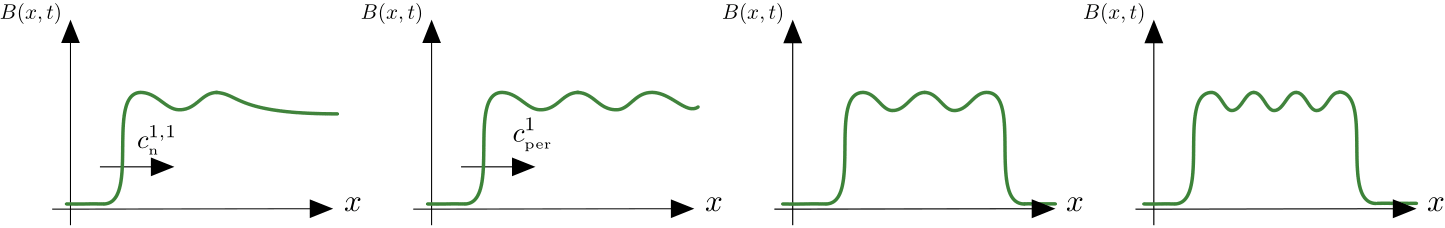}
    \caption{Four sketches of `higher order' localized patterns constructed in this paper. (a) A secondary traveling 1-front, the second one in a countable family of traveling 1-fronts between the bare soil state and a uniform vegetated state -- all traveling with different speeds -- that starts with the primary 1-front of Fig. \ref{fig:basicpatterns}(a) (Theorem \ref{T:count1fronts}). (b) The limiting orbit of the family sketched in (a): a 1-front connection between the bare soil state and a spatially periodic vegetation state (Theorem \ref{T:frontstoperiodic}). (c,d) The first 2 representations of a (countable) `higher order' family of localized (stationary, homoclinic 2-front) spot patterns with an increasing number of `spatial oscillations' (Theorem \ref{T:higherorder2fronts}).}
    \label{fig:sketchessection1}
\end{figure}
\\
Finally, we illustrate our analytical findings by several numerical simulations of PDE model (\ref{eq:meron})/(\ref{eq:fast}) -- see also Fig. \ref{fig:basicpatterns}. We did not systematically investigate the question whether all heteroclinic/homoclinic/periodic (multi-)front orbits of (\ref{eq:fast}) constructed here indeed may be (numerically) observed as stable patterns in (\ref{eq:fast}), either for general parameter combinations in (\ref{eq:scaled}) or for the more restricted class of ecologically relevant parameter combinations. This will be the subject of future work, as will be the analytical question about the spectral stability of the constructed patterns. These issues will be discussed more extensively in section \ref{sec:disc}, where we will also discuss further implications of our findings -- both from the mathematical as well as from the ecological point of view.
\\ \\
The set-up of this paper is as follows. Section \ref{sec:setup} is a preparatory section: in section \ref{sec:fast} and \ref{sec:slowred} we consider the fast and slow reduced problems associated with (\ref{eq:fast}), followed by a brief section -- section \ref{sec:crithom} -- in which we discuss the nature (and stability) of the critical points of (\ref{eq:fast}) as uniform vegetated states in (\ref{eq:scaled}); in section \ref{sec:slowfull} we study the full, perturbed, slow flow on the slow manifolds (leading to Theorem \ref{T:slowflow}). All localized patterns are constructed in section \ref{sec:locfront}, which begins with (another) preparatory section -- section \ref{sec:touchdown} -- in which we set up the geometry of orbits `jumping' between slow manifolds. The primary traveling 1-front patterns of Fig. \ref{fig:basicpatterns}a are constructed in section \ref{sec:primtrav1front}, the associated higher order 1-fronts of Figs. \ref{fig:sketchessection1}a and \ref{fig:sketchessection1}b in section \ref{sec:frontsbyM+eps}. Stationary patterns are considered in \ref{sec:stand1fronts} -- on 1-fronts -- and \ref{sec:stand1fronts} -- on 2-fronts of spot and gap type as shown in Figs. \ref{fig:basicpatterns}b, \ref{fig:basicpatterns}c and Fig. \ref{fig:sketchessection1}(c,d); various families of spatially periodic multi-front patterns -- including the basic ones of Fig. \ref{fig:basicpatterns}d -- are constructed in section \ref{sec:perpatt}. Section \ref{sec:SimDis} begins with section \ref{sec:sims} in which we show various numerically obtained patterns -- some of them beyond the analysis of the present paper -- and ends with discussion section \ref{sec:disc}.

\begin{remark}
While the original model (\ref{eq:meron}) has 8 parameters -- $(\Lambda, \Gamma, R, K, E, M, N, P)$ -- (neglecting $ D_B, D_W $ which are represented by $\varepsilon$), rescaled model (\ref{eq:scaled}) has only 5 parameters -- $ (a, \Psi, \Phi, \Omega, \Theta)$. We will formulate our results by stipulating conditions on $(a, \Psi, \Phi, \Omega, \Theta) $ and refrain from giving a corresponding range for the original parameters. Moreover, we notice that we have implicitly assumed that $\alpha > 0$, i.e. that $EK > 1$ (\ref{scalingpars}). This is a technical assumption (and not unrealistic from ecological point of view), the case $0< EK < 1$ can be treated in a completely analogous way -- see \ref{Appendix:scaling}.
\end{remark}

\section{Set-up of the existence problem}
\label{sec:setup}
We first notice that (\ref{eq:fast}) is the `fast' description of the `spatial ODE' associated with (\ref{eq:scaled}). By introducing $X = \varepsilon \xi \; (= \eps(x - ct))$ we obtain its equivalent slow form,
\begin{eqnarray}
	\left\{	
\begin{aligned}
\varepsilon b_{X}\;\; &= \;\;p,\\
\varepsilon 	p_{X}\;\; &= \;\;wb^3 - wb^2 + (1-aw)b - cp,\\
w_{X} \;\;&=\;\;q,\\
q_{X} \;\;&=-\Psi +\left[\Phi + \Omega b + \Theta b^2\right]w - \eps cq.
\label{eq:SRSepsilonnonzero}
\end{aligned}
\right.
\label{eq:slow}
\end{eqnarray}
Note that these systems possess a $c \to -c$ symmetry that reduces to a reversibility symmetry for $c=0$,
\begin{eqnarray}
\label{eq:symmc}
(c, \xi, p, q) \rightarrow (-c, -\xi,-p,-q) \; \; {\rm or} \; \; (c, X, p, q) \rightarrow (-c, -X,-p,-q).
\end{eqnarray}

\subsection{The fast reduced problem}
\label{sec:fast}

\begin{figure}[t]
    \centering
     \includegraphics[scale=0.4]{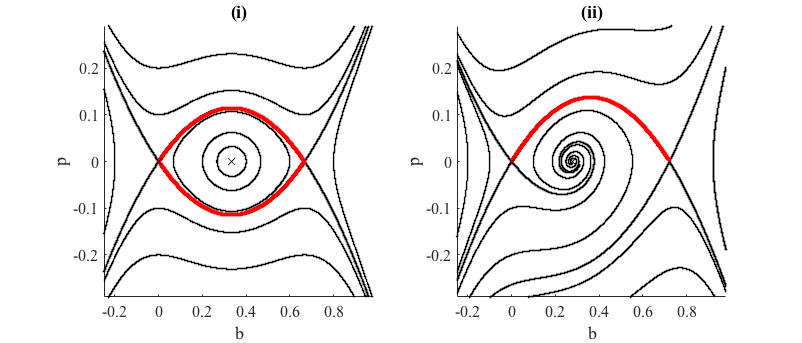}
    \caption{Numerical simulations of dynamics of the fast reduced system \eqref{eq:fast} for $a = \frac14$ and two choices of $w_0 \in \mathcal{U}_a$ (\ref{condw0}), both featuring a heteroclinic orbit between the saddle points $(0,0)$ and $(b_+(w_0),0)$: {\bf (i)} $w_0 = 9/(2+9a), c = c_{\pm}(w_0) = 0 $; {\bf (ii)} $w_0 = 9/(2+9a)+ 0.1, c = c_{+}(w_0) \approx 0.17 $.}
    \label{fig:FRS}
\end{figure}

The fast reduced limit problem associated to (\ref{eq:fast}) is a two-parameter family of planar systems that is obtained by taking the limit $\eps \to 0$ in (\ref{eq:fast}),
\begin{eqnarray}
b_{\xi \xi} = w_0b^3 - w_0b^2 + (1-aw_0)b - cb_{\xi}, \; \; (w,q) \equiv (w_0,q_0) \in \mathbb{R}^2.
\label{eq:fastblue}
\end{eqnarray}
These planar systems can have up to 3 (families of) critical points (parameterized by $(w_0,q_0)$) given by,
\begin{eqnarray}
(b_0,p_{0})  = (0,0),&&
(b_{\pm},p_{\pm}) = (b_{\pm}(w_0),0) = \left(\frac{1}{2}\pm \sqrt{a + \frac14 - \frac{1}{w_0}},0\right).
\label{b0bpm}
\end{eqnarray}
Clearly, $(b_0,w_0)$ represents the (homogeneous) bare soil state $B(x,t) \equiv 0$, the other two solutions correspond to uniform vegetation states and only exist for $w_0 > 4/(1+4a)$. The critical points also determine 3 two-dimensional invariant (slow) manifolds, $\mathcal{M}_0^{0}$ and $\mathcal{M}_0^{\pm}$,
\begin{equation}
\begin{array}{rcl}
\label{def:M00Mpm0}
\mathcal{M}^{0}_{0} &=& \left\{(b,p,w,q) \in \mathbb{R}^4: b=0, p = 0\right\},\\
\mathcal{M}^{\pm}_{0} &=& \left\{ (b,p,w,q) \in \mathbb{R}^4: b=b_{\pm}(w)=\frac{1}{2} \pm \sqrt{a + \frac14 - \frac{1}{w}}, p = 0\right\}.
\end{array}
\end{equation}
A straightforward analysis yields that the critical points $(b_{+},0)$ are saddles for all $c \in \mathbb{R}$ and that the points $(b_0,p_0) = (0,0)$ are saddles for all $c$ as long as $w_0 < 1/a$. Therefore, we consider in this paper $w_0$ such that,
\beq
w_0 \in \mathcal{U}_a = \left\{ w_0 \in \mathbb{R} ~|~ \frac{4}{1+4a} < w_0 < \frac{1}{a} \right\} \, ,
\label{condw0}
\feq
so that (parts of) the manifolds $\mathcal{M}_0^{0}$ and $\mathcal{M}_0^{+}$ are normally hyperbolic for all $w_0$ that satisfy (\ref{condw0}) (and thus persist as $\eps$ becomes nonzero \cite{Jones1995,Kaper1999}); moreover, all stable and unstable manifolds $W^{s,u}(\mathcal{M}_0^{0})$ and $W^{s,u}(\mathcal{M}_0^{+})$ are 3-dimensional. (In this paper, we do not consider the manifold $\mathcal{M}_0^{-}$ for several reasons: {\it (i)} it is not normally hyperbolic in the crucial case of stationary patterns (i.e. for $c = 0$, under the -- natural -- assumption that the water concentration $w_0$ does not become negative), {\it (ii)} critical points for the full system (\ref{eq:fast}) that limit on $\mathcal{M}_0^{-}$ as $\eps \to 0$ cannot correspond to stable homogeneous states of PDE (\ref{eq:scaled}) -- see section \ref{sec:crithom}.)
\\ \\
The manifolds $W^{s,u}(\mathcal{M}_0^{0})$ and $W^{s,u}(\mathcal{M}_0^{+})$ are determined by the stable and unstable manifolds of  $(0,0)$ and $(b_{+},0)$. By the (relatively) simple cubic nature of (\ref{eq:fastblue}) we do have explicit control over these manifolds in the relevant cases that they collide, i.e. that there is a heteroclinic connection between $(0,0)$ and $(b_{+},0)$. Although this is a classical procedure -- see \cite{Murray2013} -- we provide a brief sketch here.
\\ \\
We assume that a heteroclinic solution of (\ref{eq:fastblue}) between $(0,0)$ and $(b_{+},0)$ can also be written as a solution of the first order equation
\beq
\label{eq:firstorder}
b_{\xi} = n b(b_+(w_0)-b),
\feq
where $n$ is a free pre-factor (and we know that this assumption provides all possible heteroclinic connections). Taking the derivative (w.r.t. $\xi$) yields an equation for $b_{\xi\xi}$ that must equal (\ref{eq:fastblue}) -- that we write as $b_{\xi \xi} = w_0b(b-b_-)(b-b_+) - cb_\xi$. Working out the details yield explicit expressions for $n$ and $c$,
\beq
n = n^\pm(w_0) = \pm \sqrt{\frac12 w_0}, \; \; c = c^\pm(w_0) = \pm \sqrt{\frac12 w_0}
\left(3\sqrt{a + \frac14 - \frac{1}{w_0}} - \frac12\right).
\label{explnc}
\feq
Thus, for a given $c$, there is a heteroclinic connection between $\mathcal{M}_0^{0}$ and $\mathcal{M}_0^{+}$ at the `level' $w_0 = w^\pm_h(c)$ if $w_0$ solves (\ref{explnc}). A direct calculation yields that $ c^\pm(w_0) $ are strictly monotonic function of $ w_0 $ with inverse
\beq
\label{wafoc}
w^\pm_h(c) = \frac{4(9 + 2 c^2)^2}{\left(3 \sqrt{2 c^2(1+4a) + 4(2 + 9a)} \mp \sqrt{2} c \right)^2} \, .
\feq
We conclude that for a given $c$, there may be `parabolic' -- by the relation between $b$ and $p$ (\ref{eq:firstorder}) -- two-dimensional intersections $W^u(\mathcal{M}_0^{0}) \cap W^{s}(\mathcal{M}_0^{+})$ and $W^s(\mathcal{M}_0^{0}) \cap W^{u}(\mathcal{M}_0^{+})$ explicitly given by,
\begin{equation}
\label{WsWu0}
\begin{array}{l}
W^u(\mathcal{M}_0^{0}) \cap W^{s}(\mathcal{M}_0^{+})=
\left\{0 < b < b_+(w^+_h), p = n^+(w^+_h)b(b_+(w^+_h)-b), w = w^+_h \right\},\\[1mm]
W^s(\mathcal{M}_0^{0}) \cap W^{u}(\mathcal{M}_0^{+})=
\left\{0 < b < b_+(w^-_h), p = n^-(w^-_h)b(b_+(w^-_h)-b), w = w^-_h \right\}
\end{array}
\end{equation}
(where we recall that $q=q_0 \in \mathbb{R}$ is still a free parameter). See Lemma \ref{lem:Idown} for a further discussion and analysis (for instance on the allowed $c$-intervals for which the heteroclinic connections exist: $w^\pm_h(c)$ must satisfy (\ref{condw0})).
\\ \\
In the case of stationary patterns ($c=0$), fast reduced limit problem (\ref{eq:fastblue}) is integrable, with Hamiltonian  $\mathcal{H}_{f}$ given by,
\begin{eqnarray}
\mathcal{H}_{f}(b,p; w_0) =
\frac{1}{2} p^2 - \frac{1}{2}(1-aw_0)b^2 + \frac{1}{3}w_0b^3 - \frac{1}{4}w_0b^4,
\label{eq:HamiltonianFRS}
\end{eqnarray}
which is gauged such that $\mathcal{H}_{f}(0,0; w_0) = 0$. This system has a heteroclinic connection between $(0,0)$ and $(b_+^0,0)$ for $w_0 = w^\pm_h(0)$ such that $\mathcal{H}_{f}(b_+(w_0),0; w_0) = \mathcal{H}_{f}(0,0; w_0) = 0$.  It follows by (\ref{b0bpm}) and (\ref{eq:HamiltonianFRS}) that $w^+_h(0) = w^-_h(0) = 9/(2+9a)$ (which agrees with (\ref{wafoc})) -- see Fig. \ref{fig:FRS}.

\subsection{The slow reduced limit problems}
\label{sec:slowred}

\begin{figure}[t]
	\centering
	\includegraphics[width=12cm]{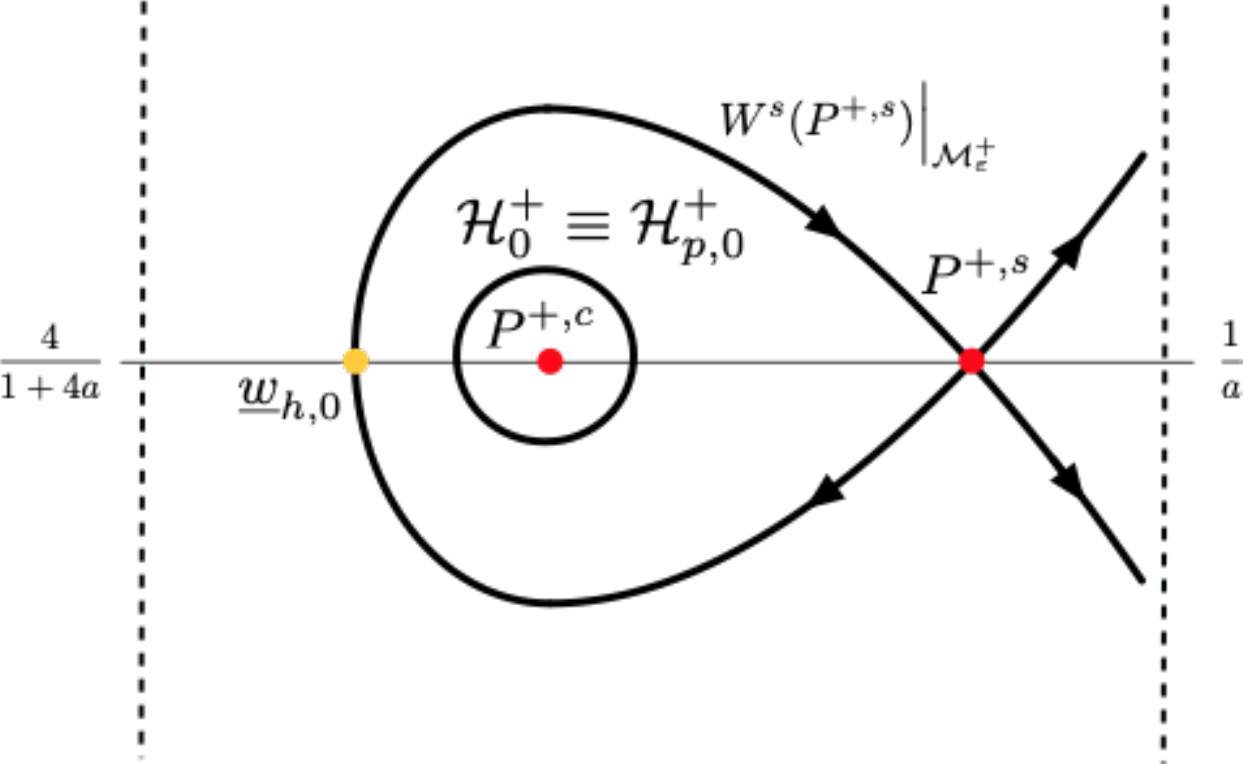}
	\caption{Phase portrait of the unperturbed flow (\ref{eq:SRSplus}) on $\M^+_0$ for parameters $(a, \Psi, \Phi, \Omega, \Theta)$ such that (\ref{asshomM+eps}) holds.}
	\label{fig:slowM+0}
\end{figure}

The slow reduced limit problem is obtained by taking the limit $\eps \to 0$ in (\ref{eq:SRSepsilonnonzero}). It is a planar problem in $(w,q)$,
\begin{eqnarray}
w_{XX} = -\Psi +\left[\Phi + \Omega b + \Theta b^2\right]w.
\label{eq:slowblue}
\end{eqnarray}
restricted to $(p,b)$ such that,
\[
p=0, \; \; wb^3 - wb^2 + (1-aw)b = 0
\]
i.e. (\ref{eq:slowblue}) describes the (slow) flow on the (slow) manifolds $\mathcal{M}_0^{0}$ and $\mathcal{M}_0^{\pm}$ (\ref{def:M00Mpm0}). The flow on $\mathcal{M}_0^{0}$ is linear,
\begin{eqnarray}
w_{XX} & =&- \Psi + \Phi w,
\label{eq:SRSzero}
\end{eqnarray}
with critical point $P^0_0 = (0,0,\Psi/\Phi, 0) \in \mathcal{M}_{0}^0$ of saddle type -- that corresponds to the uniform bare soil state $(B(x,t),W(x,t)) \equiv (0,\Psi/\Phi)$ of (\ref{eq:scaled}) -- that has the stable and unstable manifolds (on $\mathcal{M}_0^{0}$) given by
\begin{equation}
\begin{array}{rcl}
W^s(P^0_0)|_{\mathcal{M}_{0}^0} &:=& \ell^{s}_0 = \left\{(b,p,w,q) \in \M^0_0:
q = -\sqrt{\Phi}\left(w-\frac{\Psi}{\Phi}\right)\right\},\\
W^u(P^0_0)|_{\mathcal{M}_{0}^0} &:=& \ell^{u}_0 = \left\{(b,p,w,q) \in \M^0_0:
q = \sqrt{\Phi}\left(w-\frac{\Psi}{\Phi}\right)\right\}
\label{deflsu}
\end{array}
\end{equation}
Since we focus on orbits -- patterns -- that `jump' between $\mathcal{M}_0^{0}$ and $\mathcal{M}_0^{+}$ (in the limit $\eps \to 0$), we do not consider the flow on $\mathcal{M}_0^{-}$ but focus on (the flow on) $\mathcal{M}_0^{+}$,
\begin{eqnarray}
w_{XX} = - \mathcal{A} + \left(\mathcal{B}+ a\Theta \right) w + \mathcal{C} w \sqrt{a + \frac14 - \frac{1}{w}},
\label{eq:SRSplus}
\end{eqnarray}
where
\beq
\label{defABC}
\A = \Psi + \Theta \geq 0, \; \B = \Phi + \frac12 \Omega + \frac12 \Theta \in \mathbb{R}, \; \; \C = \Omega + \Theta \in \mathbb{R},
\feq
and we notice explicitly that $\B$ and $\C$ may be negative (since $\Omega$ may be negative). For $w$ satisfying (\ref{condw0}), we define,
\beq
\label{defWD}
\W = \sqrt{a + \frac14 - \frac{1}{w}} \geq 0, \; \; \D = \B + a\Theta - \left(a + \frac14\right) \A \in \mathbb{R},
\feq
and conclude that the critical points $P_0^{+,j} = (b_+(w_0^{+,j}),0,w_0^{+,j},0) \in \M^+_0$ are determined as solutions of the quadratic equation,
\beq
\label{critM+0}
\A \W^2 + \C \W + \D = 0.
\feq
Thus, the points $P_0^{+,j}$ exist for parameter combinations such that $\C^2 - 4 \A \D >  0$. There are 2 critical points if additionally $\C < 0$ and $\D > 0$ and only 1 if $\D < 0$.
\\ \\
Clearly, the flow (\ref{eq:SRSplus}) is integrable, with Hamiltonian given by
\begin{equation}
\label{eq:HamiltonianSRS}
\mathcal{H}_{0}^+(w,q) = \frac{1}{2}q^2   + \mathcal{A}w - \frac{1}{2}\left(\mathcal{B}+a\Theta\right)w^2 - \C \mathcal{J}_{0}^+(w),
\end{equation}
with, for $\ta = a + \frac14$,
\begin{equation}
\label{defJ+0}
\mathcal{J}_{0}^+(w) =
\frac{1}{4 \ta}\left(2\ta w - 1\right)\sqrt{\ta w^2 - w} -
\frac{1}{8\ta \sqrt{\ta}}
\ln{\left|\frac12\left(2\ta w - 1\right) + \sqrt{\ta} \sqrt{\ta w^2 - w}\right|}.
\end{equation}
Hence, if non-degenerate, the critical points $P_0^{+,j}$ are either centers -- $P_0^{+,c}$ -- or saddles -- $P_0^{+,s}$. Notice that, except the uniform bare soil state $(0,\Psi/\Phi)$, all critical points correspond to uniform vegetation states $(B(x,t),W(x,t)) \equiv (\bar{B},\bar{W})$ in (\ref{eq:scaled}) -- see section \ref{sec:crithom}. In the case that there is only 1 critical point $P_0^{+} \in \mathcal{M}^+_0$, it can either be of saddle or center type: $P_0^{+}$ is a saddle if,
\beq
\label{saceM+0}
\E = \B + a \Theta + \frac12 \C \left(\W + \frac{a + \frac14}{\W}\right) > 0
\feq
where $\W > 0$ is the solution of (\ref{critM+0}). We notice that the stable and unstable manifolds (restricted to $\M^+_0$) of the saddle point $P_0^{+,s} \in \M^+_0$ are represented by,
\beq
\label{manifoldsP+s}
W^s(P_0^{+,s}) \cup W^u(P_0^{+,s})|_{\M^+_0} = \left\{(b,p,w,q) \in \M^+_0: \mathcal{H}_{0}^+(w,q) \equiv \mathcal{H}_{0}^{+,s} := \mathcal{H}_{0}^+(w_0^{+,s},0)\right\}.
\feq
In the upcoming analysis, we will be especially interested in the case of 2 critical points $P_0^{+,s}, P_0^{+,c} \in \mathcal{M}^+_0$, therefore we investigate this situation on some more detail. First, we introduce $\D^{SN}$ and $\sigma \geq 0$ by setting,
\beq
\label{defsig}
\D(\sigma^2) = \D^{SN} - \A \sigma^2 = \frac{\C^2}{4 \A} - \A \sigma^2 > 0: \sigma = \sqrt{\frac{\D - \D^{SN}}{\A}},
\feq
so that the solutions of (\ref{critM+0}) are given by $\W = \W^{SN} \pm \sigma = - \frac{\C}{2 \A} \pm \sigma$. We rewrite \eqref{eq:SRSplus} in terms of $(a,\A,\C,\D)$
\begin{eqnarray}
w_{XX} = - \mathcal{A} + \left(\D + \left(a+\frac14\right)\A \right) w + \mathcal{C} w \sqrt{a + \frac14 - \frac{1}{w}}.
\label{eq:SRSplus_}
\end{eqnarray}
Clearly, $\sigma = 0$ corresponds to the degenerate saddle-node case in which $P_0^{+,s}$ and $P_0^{+,c}$ merge,
\beq
\label{defPSNwSN}
P_0^{SN} = (b_+(w_0^{SN}),0,w_0^{SN},0) \; \; {\rm with} \; \; w_0^{SN} = \frac{4 \A^2}{(1 + 4a)\A^2 - \C^2}, \quad (1 + 4a)\A^2 - \C^2 \neq 0 \, ,
\feq
where we note that $w_0^{SN}$ satisfies (\ref{condw0}) for $0 < \C^2 < \A^2$ (independent of $a$).
In fact, we can consider the unfolding of the saddle-node bifurcation by the additional assumption
that $0 < \sigma \ll 1$,
\beq
\label{unfSN}
w_0^{+,j} = w_0^{SN} \pm w_1^{SN} \sigma + \mathcal{O}(\sigma^2) =
w_0^{SN} \pm 2 \W^{SN} (w_0^{SN})^2 \sigma + \mathcal{O}(\sigma^2),
\feq
$(j=1,2)$, where the $+$-case represents the saddle $P_0^{+,s}$ and the $-$-case the center $P_0^{+,c}$: $w_0^{+,c} < w_0^{+,s}$  -- see Fig. \ref{fig:slowM+0}. In this parameter region, the slow reduced system \eqref{eq:SRSplus} features a homoclinic orbit $ (w_{\rm hom},q_{\rm hom})$ to $P_0^{+,s}$ and a family of periodic solutions around the center point $P_0^{+,c}$ (Fig. \ref{fig:slowM+0}).

\begin{remark}
\label{rem:parameters}
We conclude from (\ref{eq:SRSplus_}) that the reduced slow flow on $\mathcal{M}^+_0$ is fully determined by the values of $(a,\A, \B, \D)$. Clearly, the (linear) mapping $(\Psi, \Phi, \Omega, \Theta) \mapsto (\A, \B, \D)$ has a kernel: we can vary one of the parameters -- for instance $\Phi$ -- and determine $(\Psi, \Omega, \Theta)$ such that this does not have an effect on the reduced flow (\ref{eq:SRSplus_}) on $\mathcal{M}^+_0$  (by choosing $(\Psi(\Phi), \Theta(\Phi), \Omega(\Phi))$ such that $(\A, \B, \D)$ are kept at a chosen value). We will make use of this possibility extensively in section \ref{sec:locfront}.
\end{remark}

\subsection{Critical points and homogeneous background states}
\label{sec:crithom}

\begin{figure}[t]
	\centering
	\begin{subfigure}[b]{0.32\textwidth}
		\centering
		\includegraphics[width=\textwidth]{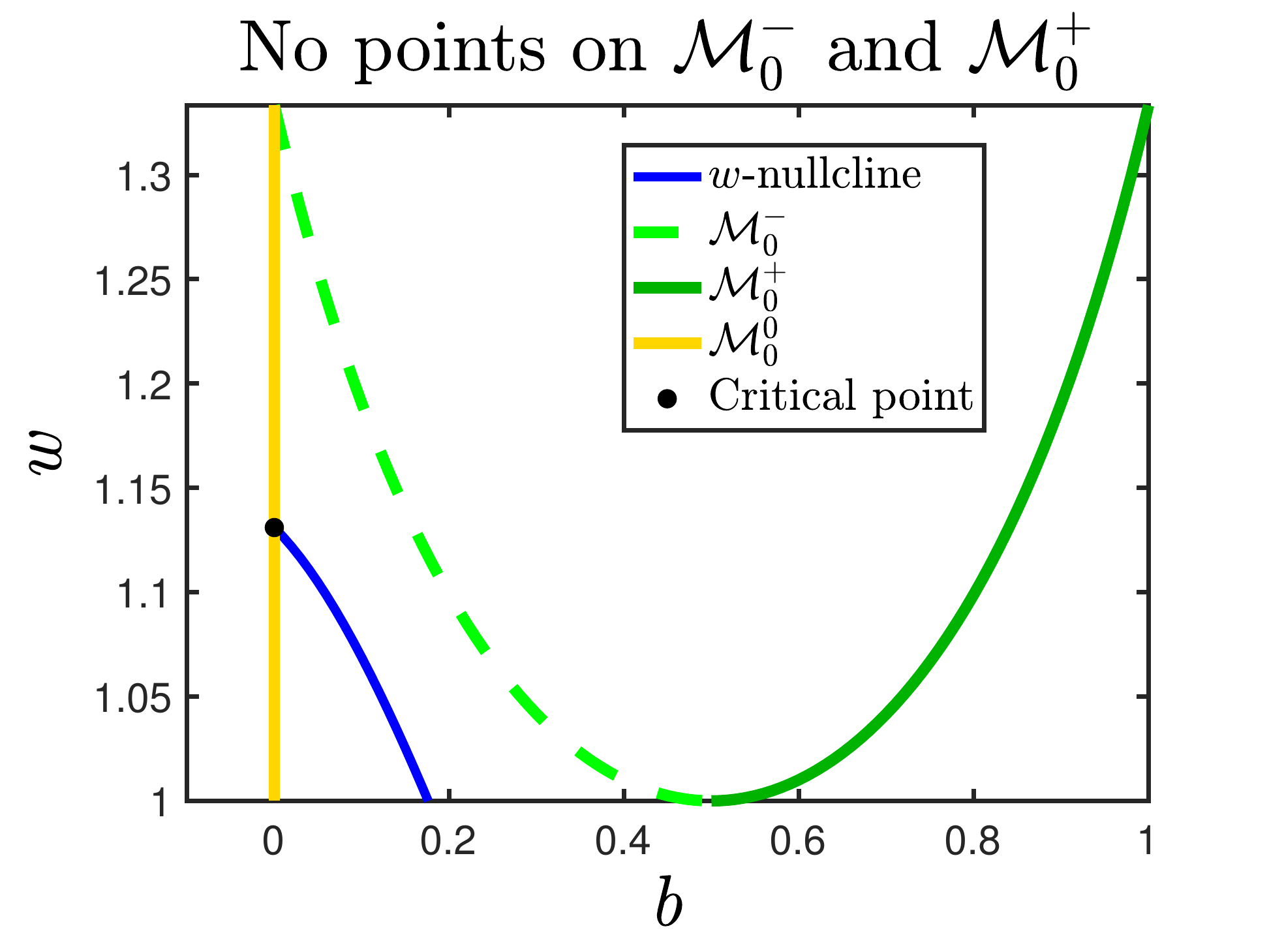}
		\caption{\label{fig:nullcline_none} No intersections of the $w$-nullcline with either $\mathcal{M}^{-}_0$ or $\mathcal{M}^{+}_0$ ($ a = 0.75, \Psi =   0.1131, \Omega = 0.0369, \Phi = 0.1$, $\Theta = 0.2131$).}		
	\end{subfigure}%
	\hspace{.2cm}
	\begin{subfigure}[b]{0.32\textwidth}
		\centering
		\includegraphics[width=\textwidth]{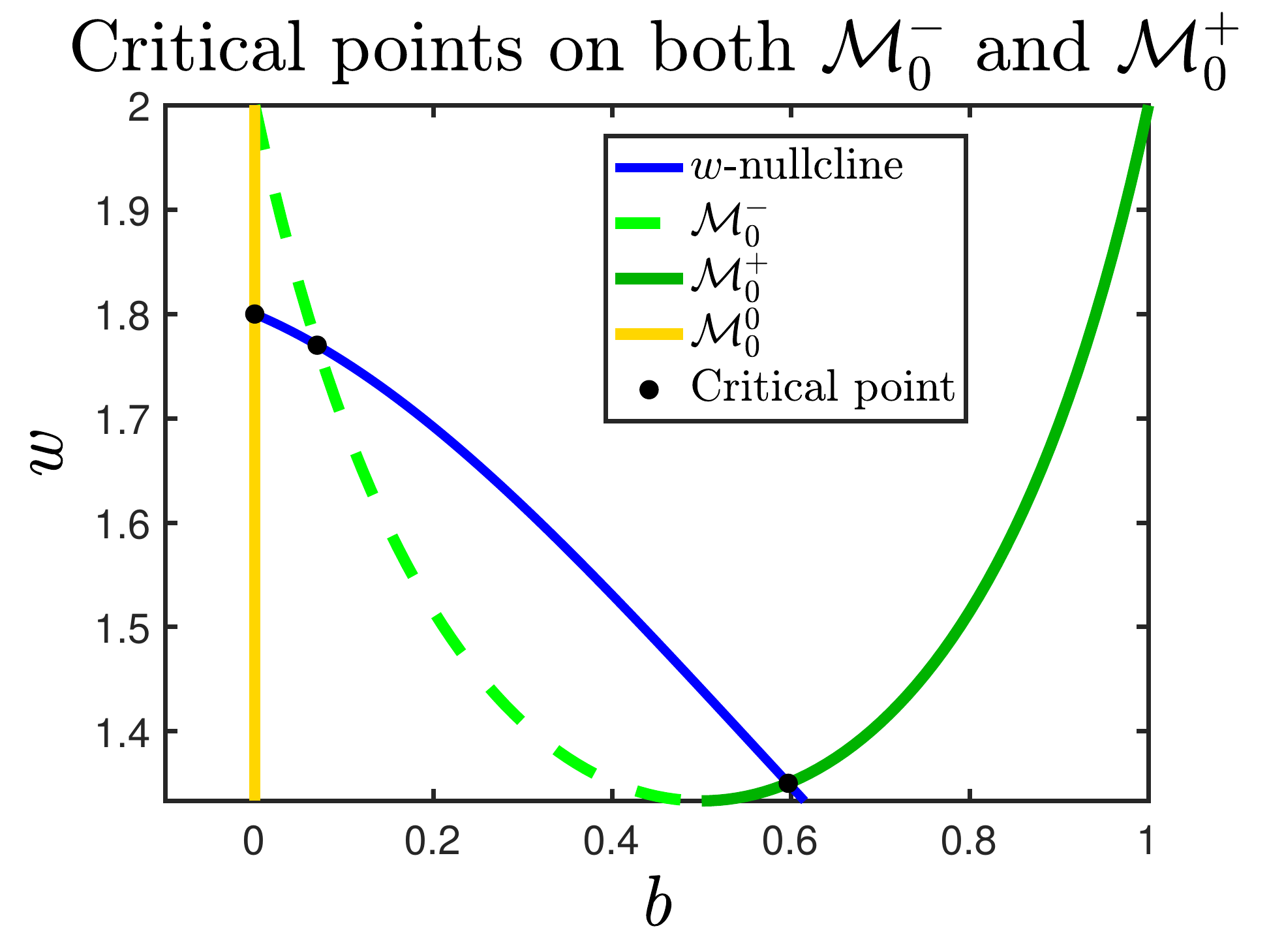}
		\caption{\label{fig:nullcline_one} A unique intersection of the $w$-nullcline with both $\mathcal{M}^{-}_0$ and $\mathcal{M}^{+}_0$ ($ a = 0.1, \Psi = 1.9, \Omega = 0.1, \Phi = 0.3$, $ \Theta = 0.5$).}
	\end{subfigure}%
	\hspace{.2cm}
\begin{subfigure}[b]{0.32\textwidth}
	\centering
	\includegraphics[width=\textwidth]{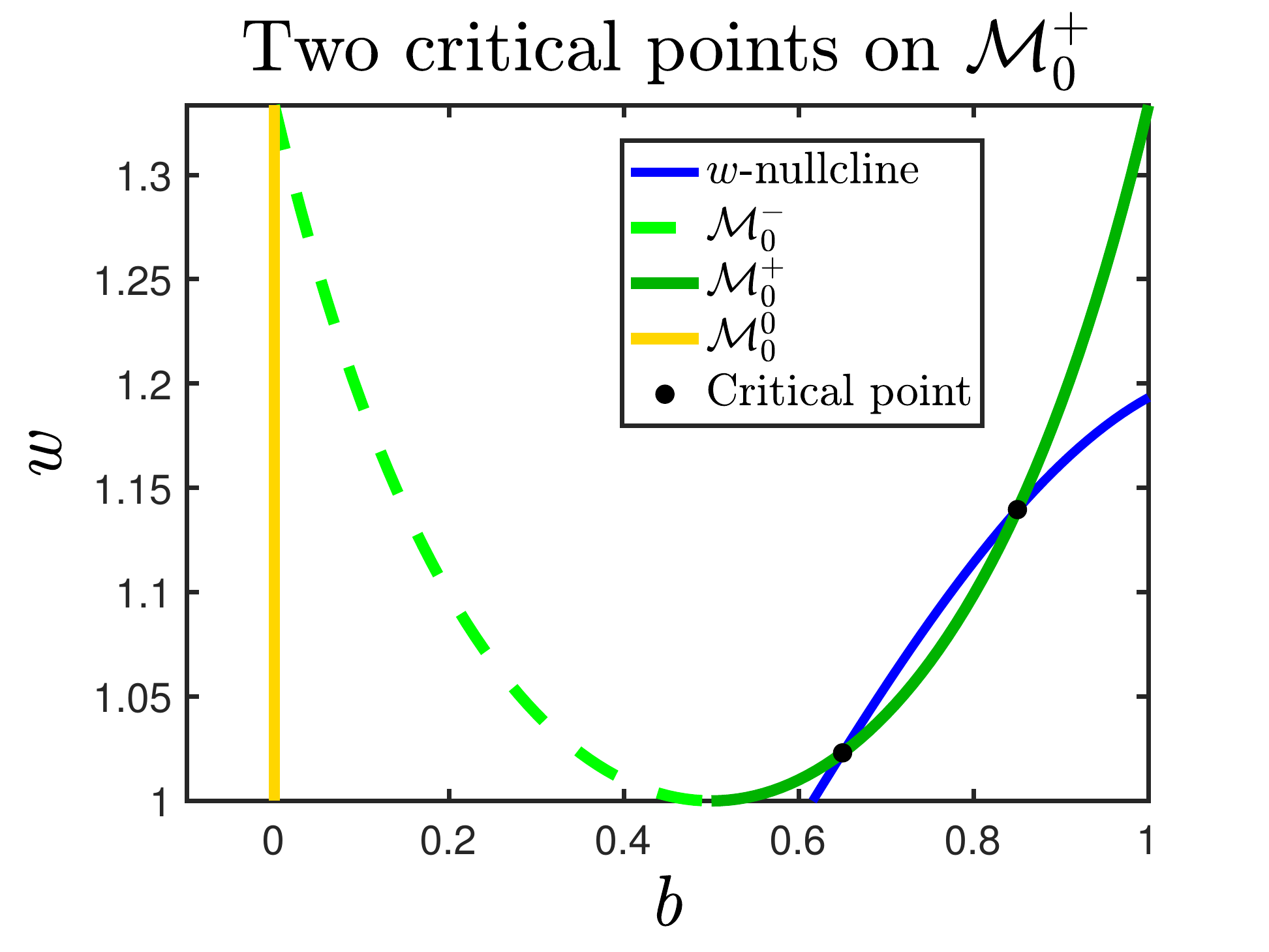}
	\caption{\label{fig:nullcline_two} Two intersections of the $w$-nullcline with $\mathcal{M}^{+}_0$ and none with $\mathcal{M}^{-}_0$ ($ a = 0.75, \Psi = 0.2983, \Omega =-0.4517, \Phi = 0.5$, $\Theta =  0.2017$). }
\end{subfigure}%

\caption{Various relative configurations of the nullclines (\ref{eqbackgroundstates}) and the associated critical points for $w \in \mathcal{U}_a$ (\ref{condw0}).}
\label{fig:nullclines}
\end{figure}

Since the critical points $P^j = (b_j,p^j,w^j,q^j)$ of the full $\eps \neq 0$ system (\ref{eq:fast}) must have $p^j=q^j=0$, their $(b,w)$ coordinates are determined by the intersections of the $b$- and $w$-nullclines,
\beq
\label{eqbackgroundstates}
wb^3 - wb^2 + (1-aw)b = 0, \; \;
-\Psi +\left[\Phi + \Omega b + \Theta b^2\right]w = 0,
\feq
where we recall that the $b$-nullcline determines the slow manifolds $\M^0_0$ and $\M^{\pm}_0$ -- see Fig. \ref{fig:nullclines}. Hence, all critical points $P^j$ must correspond to critical points of the slow reduced flows on either one of the (unperturbed) slow manifolds $\M^0_0$, $\M^+_0$ or $\M^-_0$. This immediately implies that $P^1 = P^0_0 = (0,0,\Psi/\Phi,0) \in \M^0_0$. The (potential) critical points on $\M^-_0$ can be determined completely analogously to $P_0^{+,j} \in \M^+_0$ in section \ref{sec:slowred} -- the only difference is that the term $+\C \W$ in (\ref{critM+0}) must be replaced by $-\C \W$. Thus, we conclude that there are two additional critical points $P^2$ and $P^3$ if $\C^2 - 4 \A \D \geq 0$ (and that $P^1 = P^0_0$ is the unique critical point if $\C^2 - 4 \A \D \leq 0$). Moreover, if $\C^2 - 4 \A \D \geq 0$ then,
\\
$\bullet$  if $\D < 0$, then $P^2 = P^{-}_0 \in \M^-_0$ and $P^3 = P^{+}_0 \in \M^+_0$;
\\
$\bullet$  if $\D > 0$ and $\C > 0$, then $P^{2} = P^{-,1}$, $P^{3} = P^{-,2}$ and both $P^{-,j} \in \M^-_0$;
\\
$\bullet$  if $\D > 0$ and $\C < 0$, then $P^{2} = P^{+,1}$, $P^{3} = P^{+,2}$ and both $P^{+,j} \in \M^+_0$.
\\ \\
Naturally, the critical points $P^j$ correspond to homogeneous background states $(B(x,t),W(x,t)) \equiv (\bar{B}^j,\bar{W}^j)$ of the full PDE (\ref{eq:scaled}). In this paper, we focus on the {\it existence} of patterns in (\ref{eq:scaled}) and do not consider the stability of these patterns (which is the subject of work in progress). However, there is a strong relation between the local character of critical points $P^j$ in the spatial system (\ref{eq:fast}) and their (in)stability as homogeneous background pattern in (\ref{eq:scaled}) -- see for instance \cite{Dreview}. Therefore, we may immediately conclude,
\\
$\bullet$ the bare soil state $(\bar{B},\bar{W}) = (0,\Psi/\Phi)$ is stable as solution of (\ref{eq:scaled}) for $\Psi/\Phi < 1/a$, i.e. as long as $(0,\Psi/\Phi)$ corresponds to a critical point on the normally hyperbolic part of $\M^0_0$ (\ref{condw0});
\\
$\bullet$ background states $(\bar{B},\bar{W})$ that correspond to critical points on $\M^-_0$ are unstable;
\\
$\bullet$ a background state $(\bar{B},\bar{W})$ that corresponds to a center point on $\M^+_0$ is unstable;
\\
$\bullet$ a background state $(\bar{B},\bar{W})$ that corresponds to a saddle point on $\M^+_0$ is stable as solution of (\ref{eq:scaled}) if one additional (technical) condition on the parameters of (\ref{eq:scaled}) is satisfied.
\\ \\
Of course this motivates our choice to study homoclinic and heteroclinic connections between the saddle points on $\M^0_0$ and $\M^+_0$ in this paper.
\begin{remark}
\label{rem:Turing}
The singular perturbation point of view also immediately provides insight in the possible occurrence of a Turing bifurcation in (\ref{eq:scaled}). In the setting of (\ref{eq:fast}) -- with $c=0$ -- a Turing bifurcation corresponds to a reversible $1:1$ resonance Hopf bifurcation \cite{HaragusIooss2011}, i.e. the case of a critical point with 2 colliding pairs of purely imaginary eigenvalues. By the slow/fast nature of the flow of (\ref{eq:fast}), such a critical point cannot lay inside one of the $3$ possible reduced slow manifolds $\M^0_0$, $\M^-_0$ or $\M^+_0$ (critical points not asymptotically close to the boundaries must have 2 $\OO(\eps)$ and 2 $\OO(1)$ eigenvalues). Thus, critical points that may undergo a Turing/reversible $1:1$ Hopf bifurcation have to be asymptotically close to the edge of $\M^+_\eps$ where it approaches $\M^-_\eps$ (where we note that we a priori do not claim that $\M^-_\eps$ persists). Indeed, the bifurcation appears in that region -- although we refrain from going into the details. See Fig. \ref{fig:numerics_turing} for a thus found spatially periodic Turing pattern in (\ref{eq:scaled}).
\end{remark}
\begin{remark}
By directly focusing on \eqref{eqbackgroundstates} -- and thus by not following the path indicated by the singularly perturbed structure of (\ref{eq:fast}) -- the uniform vegetation background states can also be computed in a more straightforward way: assuming $ b \neq 0$, yields $ w = -\frac{1}{b^2-b-a} $, which implies that $(\Theta + \Psi) b^2 + (\Omega - \Psi) b + (\Phi - a \Psi) =0$. Hence it follows (for $(\Omega - \Psi)^2 - 4 (\Theta + \Psi)(\Phi - a \Psi) \geq 0$) that
\[
(b_{1,2},w_{1,2}) = \left(
\frac{-(\Omega - \Psi) \pm \sqrt{(\Omega - \Psi)^2 - 4 (\Theta + \Psi)(\Phi - a \Psi)}}{2 (\Theta + \Psi)}, -\frac{1}{b_{1,2}^2-b_{1,2}-a}
\right).
\]
\end{remark}

\subsection{The slow flows of the $\eps \neq 0$ system}
\label{sec:slowfull}

Condition (\ref{condw0}) was chosen such that the points $(0,0,w,q) \in \M^0_0$ and $(b_{+}(w),0,w,q) \in \M^+_0$ are saddles for the fast reduced limit problem (\ref{eq:fastblue}) (so that the associated background states may be stable as trivial, homogeneous, patterns of (\ref{eq:scaled}) -- section \ref{sec:crithom}). Thus, where (\ref{condw0}) holds, $\M^0_0$ and $\M^+_0$ are normally hyperbolic and they thus persist as $\M^0_\eps$ and $\M^+_\eps$ for $\eps \neq 0$ \cite{Jones1995,Kaper1999}. Clearly, $\M^0_0$ is also invariant under the flow of the full system (\ref{eq:fast}): $\M^0_\eps = \M^0_0$. Moreover, the flow on $\M^0_\eps$ is only a slight -- $\mathcal{O}(\eps)$ -- (linear) perturbation of the unperturbed flow (\ref{eq:SRSzero}) on $\M^0_0$ -- due to the (asymmetric) $-\eps c q$ term. As a consequence, only the orientation of the (un)stable manifolds $W^{s,u}(P^0_0)|_{\M^0_\eps} = \ell^{s,u}_\eps$ undergoes an $\mathcal{O}(\eps)$ change w.r.t. $\ell^{s,u}_0$ (\ref{deflsu}).
\\ \\
The situation is very different for $\M^+_\eps$. A direct perturbation analysis yields,
\beq
\label{M+eps}
\mathcal{M}^{+}_{\eps} = \left\{ (b,p,w,q) \in \mathbb{R}^4: b=b_+(w) +  \eps c q  b_1(w) + \mathcal{O}(\varepsilon^2), \;
p = \varepsilon  q p_1(w) + \mathcal{O}(\varepsilon^2)\right\},
\end{eqnarray}
with
\beq
\label{p1b1}
p_1(w) = \frac{1}{2w^2\sqrt{a + \frac14 - \frac{1}{w}}}, \; \; b_1(w) = \frac{p_1(w)}{2 w b_+(w) \sqrt{a + \frac14 - \frac{1}{w}}}
\feq
and $b_+(w)$ as defined in (\ref{b0bpm}). Since we only consider situations in which there are critical points (of the full flow) $P^{+,j}$ on $\M^+_0$, and thus on $\M^+_\eps$, we know (and use) that $\mathcal{M}_\eps^{+}$ is determined uniquely. The slow flow on $\M^+_\eps$ is given by
\begin{eqnarray}
w_{XX} = - \mathcal{A} + \left(\mathcal{B}+ a\Theta \right) w + \mathcal{C} w \sqrt{a + \frac14 - \frac{1}{w}} + \eps c q \rho_1(w) + \mathcal{O}(\eps^2),
\label{eq:slow+eps}
\end{eqnarray}
(cf. (\ref{eq:SRSplus})), with
\beq
\label{rho1}
\rho_1(w) = \left(\Omega + 2 b_+(w) \Theta\right) w b_1(w) - 1 = \left( \C + 2 \Theta \sqrt{a + \frac14 - \frac{1}{w}} \, \right)w b_1(w) - 1.
\feq
Thus, for $c \neq 0$ the flow on $\M^+_\eps$ is a perturbed integrable planar system with `nonlinear friction term' $\eps c q \rho_1(w)$. In the case that there is only one critical point $P^{+,s}$ of saddle type on $\M^+_0$ -- and thus on $\M^+_0$ -- the impact of this term is asymptotically small. The situation is comparable to that of the flow on $\M^0_\eps$ w.r.t. the flow on $\M^0_0$. The stable and unstable manifolds of $P^{+,s}$ restricted to the slow manifolds remain close:  $W^{u,s}(P^{+,s})|_{\M^+_\eps}$ is $\mathcal{O}(\eps)$ close to $W^{u,s}(P^{+,s})|_{\M^+_0}$ (for $\mathcal{O}(1)$ values of $(w,q)$) and the span $W^{u,s}(P^{+,s}) \cup W^{u,s}(P^{+,s})|_{\M^+_\eps}$ has becomes slightly asymmetric -- cf. (\ref{manifoldsP+s}). This is drastically different in the case that there are 2 critical points $P^{+,c}$ -- the center -- and $P^{+,s}$ -- the saddle -- on $\M^+_\eps$. We deduce by classical dynamical system techniques -- such as the Melnikov method (see for instance \cite{GuckHolm}) -- the following (bifurcation) properties of (\ref{eq:slow+eps}), and thus of (\ref{eq:fastblue}).
\begin{theorem}
\label{T:slowflow}
Let parameters $(a, \Psi, \Phi, \Omega, \Theta)$ of (\ref{eq:fast}) be such that there is a center $P^{+,c}= (b_+(w^{+,c}),0,w^{+,c},0)$ and a saddle $P^{+,s}= (b_+(w^{+,s}),0,w^{+,s},0)$ on $\M^+_\eps$ and assume that the unperturbed homoclinic orbit $(w_{{\rm hom},0}(X),q_{{\rm hom},0}(X))$ to $P^{+,s}$ of (\ref{eq:slowblue}) on $\M^+_0$ lies entirely in the $w$-region in which both $\M^0_0$ and $\M^+_0$ are normally hyperbolic. More explicitly, assume that,
\beq
\label{asshomM+eps}
\C^2 - 4 \A \D > 0, \C < 0, \D > 0 \; \; {\rm and} \; \; \frac{4}{1+4a} < \underline{w}_{h,0} < w^{+,c} < w^{+,s} < \frac{1}{a}
\feq
(\ref{defABC}), (\ref{defWD}), (\ref{condw0}), where $(\underline{w}_{h,0},0)$ is the intersection of $(w_{{\rm hom},0}(X),q_{{\rm hom},0}(X))$ with the $w$-axis  -- see Fig. \ref{fig:slowM+0}. Then, for all $c \neq 0$ (but $\mathcal{O}(1)$ w.r.t. $\eps$) and $\eps$ sufficiently small,
\\
$\bullet$ there is a co-dimension 1 manifold $\R_{\rm Hopf} = \R_{\rm Hopf}(a, \Psi, \Phi, \Omega, \Theta)$ such that a periodic solution (dis)appears in (\ref{eq:slow+eps}) -- and thus in (\ref{eq:fast}) -- for parameters $(a, \Psi, \Phi, \Omega, \Theta)$ that cross through $\R_{\rm Hopf}$; moreover, $\R_{\rm Hopf}$ is at leading order (in $\eps$) determined by $\rho_1(w^{+,c}) = 0$ (\ref{rho1});
\\
$\bullet$ there is a co-dimension 1 manifold $\R_{\rm hom} = \R_{\rm hom}(a, \Psi, \Phi, \Omega, \Theta)$ such that for $(a, \Psi, \Phi, \Omega, \Theta) \in \R_{\rm hom}$, the unperturbed homoclinic solution $(w_{{\rm hom},0}(X),q_{{\rm hom},0}(X))$ on $\M^+_0$ persists as homoclinic solution to $P^{+,s}$ of (\ref{eq:slow+eps})/(\ref{eq:fast}); moreover, $\R_{\rm hom}$ is at leading order determined by,
\beq
\label{persslowhom}
\Delta H_{\rm hom} = c \int_{\underline{w}_{h,0}}^{w^{+,s}} \rho_1(w)
\sqrt{2 \mathcal{H}_{0}^{+,s} - 2 \mathcal{A}w + \left(\mathcal{B}+a\Theta\right)w^2 + 2 \C \mathcal{J}_{0}^+(w)} \, dw = 0.
\feq
with $\mathcal{H}_{0}^{+,s}$, $\mathcal{J}_{0}^+(w)$ as defined in (\ref{manifoldsP+s}), (\ref{defJ+0}).
\\
$\bullet$ there is an open region $\mathcal{S}_{\rm per}$ in $(a, \Psi, \Phi, \Omega, \Theta)$-space -- with $\R_{\rm Hopf} \cup \R_{\rm hom} \subset \partial \mathcal{S}_{\rm per}$ -- such that for all $(a, \Psi, \Phi, \Omega, \Theta) \in \mathcal{S}_{\rm per}$, one of the (restricted) periodic solutions $(w_{p,0}(X),q_{p,0}(X))$ of the integrable flow (\ref{eq:SRSplus}) on $\M^+_0$ persists as a periodic solution $(b_{p,\eps}(X),p_{p,\eps}(X),w_{p,\eps}(X),q_{p,\eps}(X))$ of (\ref{eq:slow+eps})/(\ref{eq:fast}) on $\M^+_\eps$; the stability of the periodic orbit on $\M^+_\eps$ is determined by (the sign of) $c$.
\\
The flow on $\M^+_\eps$ is reversible for $c = 0$: there always is a one-parameter family of periodic solutions on $\M^+_\eps$ enclosed by a homoclinic loop if (\ref{asshomM+eps}) holds, i.e. the phase portrait remains as in the $\eps = 0$ case of Fig. \ref{fig:slowM+0}, it is not necessary to restrict parameters $(a, \Psi, \Phi, \Omega, \Theta)$ to $\mathcal{S}_{\rm per}$ or to $\R_{\rm hom}$ for $c=0$.
\end{theorem}
{\bf Proof.}
A periodic solution $(w_{p,0}(X),q_{p,0}(X))$ of the unperturbed flow (\ref{eq:SRSplus}) on $\M^+_0$ is described by the value $\mathcal{H}_{p,0}^+$ of the Hamiltonian $\mathcal{H}_{0}^+(w,q)$ (\ref{eq:HamiltonianSRS}), where necessarily $\mathcal{H}_{p,0}^+ \in (\mathcal{H}_{0}^{+,c}, \mathcal{H}_{0}^{+,s})$ -- with $\mathcal{H}_{0}^{+,c} < \mathcal{H}_{0}^{+,s}$ the values of $\mathcal{H}_{0}^+(w,q)$ at the center $P^{+,c}_0$, resp. saddle $P^{+,s}_0$ (cf. (\ref{manifoldsP+s})). We define $L_{p,0} = L_{p,0}(\mathcal{H}_{p,0}^+)$ as the period -- or wave length -- of $(w_{p,0}(X),q_{p,0}(X))$ and $\underline{w}_{p,0} = \underline{w}_{p,0}(\mathcal{H}_{p,0}^+)$ and $\overline{w}_{p,0} = \overline{w}_{p,0}(\mathcal{H}_{p,0}^+)$ as the minimal and maximal values of $w_{p,0}(X)$, i.e
$\underline{w}_{p,0} \leq w_{p,0}(X) \leq \overline{w}_{p,0}$ -- see Fig. \ref{fig:slowM+0}.
\\ \\
Hamiltonian $\mathcal{H}_{0}^+(w,q)$ (\ref{eq:HamiltonianSRS}) becomes a slowly varying function
in the perturbed system (\ref{eq:slow+eps}),
\[
\frac{d \mathcal{H}_{0}^+}{dX}(w,q) = \eps c q^2 \rho_1(w) + \mathcal{O}(\eps^2).
\]
Thus, unperturbed periodic solution $(w_{p,0}(X),q_{p,0}(X))$ on $\M^+_0$ persists as periodic solution $(w_{p,\eps}(X),q_{p,\eps}(X))$ of (\ref{eq:slow+eps}) on $\M^+_\eps$ -- with $|L_{p,\eps} - L_{p,0}|, |\overline{w}_{p,\eps} - \overline{w}_{p,0}| = \mathcal{O}(\eps)$ and, by definition, $\underline{w}_{p,\eps} = \underline{w}_{p,0}$ -- if,
\[
\int_0^{L_{p,\eps}} \frac{d \mathcal{H}_{0}^+}{dX}(w_{p,\eps}(X),q_{p,\eps}(X))\, dX =
\eps c \int_0^{L_{p,\eps}} (q_{p,\eps}(X))^2 \rho_1(w_{p,\eps}(X))\, dX  + \mathcal{O}(\eps^2) = 0.
\]
The approximation of $(w_{p,\eps}(X),q_{p,\eps}(X))$ by $(w_{p,0}(X),q_{p,0}(X))$ yields, together with (\ref{eq:HamiltonianSRS}),
\[
\begin{array}{rcl}
\int_0^{L_{p,\eps}} (q_{p,\eps}(X)^2 \rho_1(w_{p,\eps}(X))\, dX
&=&
\int_0^{L_{p,0}} q_{p,0}(X)^2 \rho_1(w_{p,0}(X))\, dX + \mathcal{O}(\eps)
\\
&=&
2 \int_{\underline{w}_{p,0}}^{\overline{w}_{p,0}} \rho_1(w)
\sqrt{2 \mathcal{H}_{p,0}^+ - 2 \mathcal{A}w + \left(\mathcal{B}+a\Theta\right)w^2 + 2 \C \mathcal{J}_{0}^+(w)} \, dw
+ \mathcal{O}(\eps).
\end{array}
\]
Thus, unperturbed periodic solution/pattern $(w_{p,0}(X),q_{p,0}(X))$ persists as periodic solution on $\M^+_{\eps}$ for parameter combinations such that,
\beq
\label{exslowper}
\Delta H(\mathcal{H}_{p,0}^+) =
c \int_{\underline{w}_{p,0}(\mathcal{H}_{p,0}^+)}^{\overline{w}_{p,0}(\mathcal{H}_{p,0}^+)} \rho_1(w)
\sqrt{2 \mathcal{H}_{p,0}^+ - 2 \mathcal{A}w + \left(\mathcal{B}+a\Theta\right)w^2 + 2 \C \mathcal{J}_{0}^+(w)} \, dw = 0.
\feq
Note that this expression does not depend on the speed $c$ -- see however Remark \ref{rem:higherorder} -- but that (the sign of) $c$ indeed determines the stability of $(w_{p,\eps}(X),q_{p,\eps}(X))$ on $\M^+_\eps$. For given $\mathcal{H}_{p,0}^+ \in (\mathcal{H}_{0}^{+,c}, \mathcal{H}_{0}^{+,s})$, condition (\ref{exslowper}) determines a co-dimension 1 manifold $\R_{\rm per}(\mathcal{H}_{p,0}^+)$ in
$(a, \Psi, \Phi, \Omega, \Theta)$-space for which a periodic orbit $(b_{p,\eps}(X),p_{p,\eps}(X),w_{p,\eps}(X),q_{p,\eps}(X))$ on $\M^+_\eps$ exists. Clearly $\mathcal{S}_{\rm per} \subset \cup_{\mathcal{H}_{p,0}^+ \in (\mathcal{H}_{0}^{+,c}, \mathcal{H}_{0}^{+,s})} \R_{\rm per}(\mathcal{H}_{p,0}^+)$. Moreover, $\underline{w}_{p,0}(\mathcal{H}_{p,0}^+) \uparrow w^{+,c}$ and $\overline{w}_{p,0}(\mathcal{H}_{p,0}^+) \downarrow w^{+,c}$ as $\mathcal{H}_{p,0}^+ \downarrow \mathcal{H}_{0}^{+,c}$, so that (\ref{exslowper}) indeed reduces to $\rho_1(w^{+,c}) = 0$ as $\mathcal{H}_{p,0}^+ \downarrow \mathcal{H}_{0}^{+,c}$: $\R_{\rm Hopf} = \R_{\rm per}(\mathcal{H}_{0}^{+,c})$. Note that $\rho_1(w) \to -\infty$ as $w \downarrow 4/(1+4a)$ -- recall that $\C < 0$ -- and that
\[
\rho_1\left(\frac{1}{a}\right) = a^2 \left(\Omega + 2 \Theta\right) - 1 = - \left(1 - a^2 \C \right) + a^2 \Theta
\]
can be made positive by choosing $\Theta$ sufficiently large: $\rho_1(w)$ must change sign for $\Theta$ not too small (in fact, it can be shown by straightforward analysis of (\ref{rho1}) that $\rho_1(w)$ may change sign twice (at most)). It thus follows that $\R_{\rm Hopf} \neq \emptyset$ and consequentially that $\mathcal{S}_{\rm per}$ is nonempty. Since $\underline{w}_{p,0}(\mathcal{H}_{p,0}^+) \downarrow \underline{w}_{h,0}$ and $\overline{w}_{p,0}(\mathcal{H}_{p,0}^+) \uparrow w^{+,s}$ as $\mathcal{H}_{p,0}^+ \uparrow \mathcal{H}_{0}^{+,s}$, it follows that $\Delta H(\mathcal{H}_{p,0}^+) \to \Delta H_{\rm hom}$ and thus that $\R_{\rm hom} = \R_{\rm per}(\mathcal{H}_{0}^{+,s})$, which also can be shown to be non-empty -- see Lemma \ref{lem:posigmasmal}.
\hfill $\Box$
\\ \\
Of course, Theorem \ref{T:slowflow} has a direct interpretation in terms of traveling waves in the full PDE (\ref{eq:scaled}),
\begin{corollary}
\label{cor:slowflow}
Let the conditions formulated in Theorem \ref{T:slowflow} hold, then for all $c \in \mathbb{R}$ $\mathcal{O}(1)$ w.r.t. $\eps$,
\\
$\bullet$ there is a traveling spatially periodic wave (train) solution $(B_{p,\eps}(\eps(x-ct)), W_{p,\eps}(\eps(x-ct))$ of (\ref{eq:scaled}) for $(a, \Psi, \Phi, \Omega, \Theta) \in \mathcal{S}_{\rm per}$;
\\
$\bullet$ there is a traveling pulse $(B_{{\rm hom},\eps}(\eps(x-ct)), W_{{\rm hom},\eps}(\eps(x-ct))$ in (\ref{eq:scaled}) -- homoclinic to the background state $(\bar{B}^{+,s}, \bar{W}^{+,s}) = (b_+(w^{+,s}),w^{+,s})$ -- for $(a, \Psi, \Phi, \Omega, \Theta) \in \R_{\rm hom}$.
\end{corollary}
It is possible to (locally) get full analytical control over the set $\mathcal{S}_{\rm per}$ and its boundary manifolds $\R_{\rm Hopf}$ and $\R_{\rm hom}$ in $(a, \Psi, \Phi, \Omega, \Theta)$-space by considering the unfolding of the saddle-node bifurcation on $\M^+_\eps$ as in section \ref{sec:slowred} (cf. (\ref{unfSN})).
\begin{lemma}
\label{lem:posigmasmal}
Let the conditions formulated in Theorem \ref{T:slowflow} hold, introduce $\sigma> 0$ as in (\ref{defsig}) and consider $\sigma$ sufficiently small (but still $\mathcal{O}(1)$ w.r.t. $\eps$). Then, system (\ref{eq:slow+eps})/(\ref{eq:fast}) has a periodic solution $(b_{p,\eps}(X),p_{p,\eps}(X),w_{p,\eps}(X),q_{p,\eps}(X))$ on $\M^+_\eps$ for all $(a, \Psi, \Phi, \Omega, \Theta)$ such that,
\beq
\label{SperSN}
\frac57 \sigma w^{SN}_1 \rho_1'(w^{SN}_0) + \mathcal{O}(\sigma^2) < \rho_1(w^{SN}_0) < \sigma w^{SN}_1 \rho_1'(w^{SN}_0) + \mathcal{O}(\sigma^2),
\feq
where $\rho_1(w)$, $\sigma$, $w^{SN}_0$ and $w^{SN}_1$ are explicitly given in terms of the parameters $(a, \Psi, \Phi, \Omega, \Theta)$ in (\ref{rho1}), (\ref{defsig}), (\ref{unfSN}) (with (\ref{defABC}),(\ref{defWD})): $\mathcal{S}_{\rm per}$ is given by (\ref{SperSN}) and its boundaries $\R_{\rm Hopf}$ and $\R_{\rm hom}$ by the upper, respectively lower, boundary of (\ref{SperSN}).
\end{lemma}
{\bf Proof.} For $\D$ $\mathcal{O}(\sigma^2)$ close to $\D^{SN}$ (\ref{defsig}), the unperturbed flow (\ref{eq:SRSplus}) on $\M^+_0$ can be given locally, i.e. in an $\mathcal{O}(\sigma)$ neighborhood of the critical points $P^{+,c} = (b_+(w_0^{+,c}), 0, w_0^{+,c}, 0)$ and $P^{+,s} = (b_+(w_0^{+,s}), 0, w_0^{+,s}, 0)$ with $w_0^{+,c} = w_0^{+,1} < w_0^{+,2} = w_0^{+,s}$ (\ref{unfSN}), be given by its quadratic approximation,
\[
w_{XX} = \tilde{\alpha} (w - w_0^{+,c})(w - w_0^{+,s}) + \mathcal{O}(\sigma^3) = \tilde{\alpha}
\left( (w-w^{SN}_0)^2 - \sigma^2 (w^{SN}_1)^2 \right) + \mathcal{O}(\sigma^3),
\]
(\ref{unfSN}), where $\tilde{\alpha} > 0$ is the second derivative of the right-hand side of (\ref{eq:SRSplus}) evaluated
at $w_0^{SN}$. Thus, the integral $\mathcal{H}_{0}^+$ (\ref{eq:HamiltonianSRS}) can locally be given by,
\beq
\label{HslowSN}
\mathcal{H}_{0}^+(w,q) = \frac12 q^2 - \tilde{\alpha} \left( \frac13 (w-w^{SN}_0)^3 - \sigma^2 (w^{SN}_1)^2 w \right)
+ \mathcal{O}(\sigma^4).
\feq
Direct evaluation yields that the stable/unstable manifolds of $P^{+,s}$ (restricted to $\M^+_0$)
are given by,
\beq
\label{H+sSlowSN}
\mathcal{H}_{0}^+(w,q) =  \mathcal{H}_{0}^{+,s} = \tilde{\alpha}
\sigma^2 (w^{SN}_1)^2 \left(w^{SN}_0 + \frac23 \sigma w^{SN}_1 \right)
+ \mathcal{O}(\sigma^4)
\feq
(cf. (\ref{manifoldsP+s})), which implies that the (second) intersection with the $w$-axis of the homoclinic orbit connected
to $P^{+,s}$ (in $\M^+_0$) is given by,
\beq
\label{underlinewhSN}
\underline{w}_{h,0} = w^{SN}_0 - 2 \sigma w^{SN}_1 + \mathcal{O}(\sigma^2) \left( < w_0^{+,c} =
w^{SN}_0 - \sigma w^{SN}_1 + \mathcal{O}(\sigma^2) \right)
\feq
(cf. Theorem \ref{T:slowflow})). Now, we consider parameter combinations such that $\rho_1(w)$ has a zero $\mathcal{O}(\sigma)$ close to $w^{SN}_0$, i.e. we set $\rho_1(w) = \tilde{\beta}\left( w - (w^{SN}_0 + \sigma \mu)\right) + \mathcal{O}(\sigma^2)$, where $\sigma \mu$ represents the position of the zero and $\tilde{\beta} = \rho_1'(w^{SN}_0)$. Hence, the condition $\Delta H_{\rm hom} = 0$ (\ref{persslowhom}) -- that determines the manifold $\R_{\rm hom}$ -- is at leading order (in $\sigma$) given by,
\beq
\label{persslowhomSN1}
c \tilde{\beta} \int_{\underline{w}_{h,0}}^{w^{+,s}} \left( w - (w^{SN}_0 + \sigma \mu)\right) \sqrt{2 \mathcal{H}_{0}^{+,s} +
2 \tilde{\alpha} \left( \frac13 (w-w^{SN}_0)^3 - \sigma^2 (w^{SN}_1)^2 w \right)} \, dw = 0
\feq
(\ref{H+sSlowSN}). Introducing $\omega$ by $w = w^{SN}_0 + \sigma \omega$ and using (\ref{unfSN}), (\ref{underlinewhSN}), we reduce (\ref{persslowhomSN1}) to,
\[
\tilde{\beta} \sigma^3 \sqrt{\frac23 \tilde{\alpha} \sigma} \sqrt{\sigma} \int_{-2 w^{SN}_1}^{w^{SN}_1} (\omega - \mu) \sqrt{\omega^3 - 3(w^{SN}_1)^2 \omega + 2 (w^{SN}_1)^3} \, d\omega =0.
\]
Thus, the homoclinic orbit to $P^{+,s}$ (in $\M^+_0$) persists for $\mu$ such that,
\[
\int_{-2 w^{SN}_1}^{w^{SN}_1} (\omega - \mu) (\omega - w^{SN}_1)\sqrt{\omega + 2 w^{SN}_1} \, d\omega =0
\]
(at leading order in $\sigma$ (and in $\eps$)). Straightforward integration yields that $\mu = \mu_{\rm hom} = -\frac57 \w^{SN}_1 + \mathcal{O}(\sigma)$, i.e. that on $\R_{\rm hom}$, the zero of $\rho_1(w)$ must be at $w^{SN}_0 - \frac57 \sigma w^{SN}_1 + \mathcal{O}(\sigma^2) > w_0^{+,c} = w^{SN}_0 - \sigma w^{SN}_1 + \mathcal{O}(\sigma^2)$.
\\ \\
We conclude that for $\sigma$ (and $\eps$) sufficiently small, the boundaries $\R_{\rm Hopf}$ and $\R_{\rm hom}$ of the domain $\mathcal{S}_{\rm per}$ are given by $\rho_1(w^{SN}_0 - \sigma w^{SN}_1 + \mathcal{O}(\sigma^2)) = 0$ (first bullet of Theorem \ref{T:slowflow}), respectively $\rho_1(w^{SN}_0 - \frac57 \sigma w^{SN}_1 + \mathcal{O}(\sigma^2)) = 0$ -- which is equivalent to the boundaries of (\ref{SperSN}) by Taylor expansion (in $\sigma$). Finally, we notice that for parameter values between $\R_{\rm Hopf}$ and $\R_{\rm hom}$, i.e. for which (\ref{SperSN}) holds, one of the periodic orbits between the center point and the homoclinic loop must persist -- in other words, for parameter combinations that satisfy (\ref{SperSN}), $\Delta H(\mathcal{H}_{p,0}^+) = 0$ (\ref{exslowper}) for certain $\mathcal{H}_{p,0}^+ \in (\mathcal{H}_{0}^{+,c},\mathcal{H}_{0}^{+,s})$. \hfill $\Box$

\begin{remark}
\label{rem:BogdanovTakens}
Lemma \ref{lem:posigmasmal} `rediscovers' the periodic solutions associated to a Bogdanov-Takens bifurcation. In \ref{Appendix:BogdanovTakens} we present a brief embedding of our result into the normal form approach to the Bogdanov-Takens bifurcation scenario.
\end{remark}
\begin{remark}
\label{rem:higherorder}
A higher order perturbation analysis yields that the $\mathcal{O}(\eps)$ corrections to $\R_{\rm Hopf}$ and $\R_{\rm hom}$ -- and thus to $\mathcal{S}_{\rm per}$ -- explicitly depend on $c$.
\end{remark}
\begin{remark}
\label{rem:otherslowflowperiodics}
Of course one can also establish the persistence of periodic orbits of the slow reduced flow -- as in Theorem \ref{T:slowflow} -- under the assumption that there is only one critical point $P^{+,c}$ of center type on $\M^+_0$, instead of focusing on the present case in which the reduced slow flow (\ref{eq:slowblue}) has a homoclinic orbit on $\M^+_0$ (Theorem \ref{T:slowflow}). Since we decided to focus on situations in which there is a saddle point on $\M^+_0$ -- that is potentially stable as homogeneous background state in (\ref{eq:scaled}) (section \ref{sec:crithom}) -- we do not consider this possibility here. Note however that the analysis of this case is essentially the same as presented here. See also Remark \ref{rem:moreonslowperiodics}.
\end{remark}

\section{Localized front patterns}
\label{sec:locfront}

In this section we use the slow-fast geometry of the phase space associated to (\ref{eq:fast}) to establish a remarkably rich variety of localized vegetation patterns (potentially) exhibited by model (\ref{eq:scaled}). First, we consider various kinds of traveling and stationary `invasion fronts' that connect the bare soil state to a uniform or an `oscillating' vegetation state and their associated bifurcation structures (sections \ref{sec:primtrav1front}, \ref{sec:frontsbyM+eps} and \ref{sec:stand1fronts}), next we study stationary homoclinic 2-front spot and gap patterns (section \ref{sec:stand2fronts}) and finally spatially periodic multi-front (spot/gap) patterns (section \ref{sec:perpatt}). As starting point, we need to control the intersection of $W^u(P^0)$ and $W^s(\M^+_\eps)$.

\begin{remark}
\label{rem:symmetry}
We start by considering localized patterns that correspond to orbits in $W^u(P^0)$, i.e. patterns that approach the bare soil state $(\bar{B},\bar{W}) = (0,\Psi/\Phi)$ of (\ref{eq:scaled}) as $x \to -\infty$. In fact, the upcoming results on 1-fronts are all on orbits in (\ref{eq:fast}) that connect $P^0 \in \M^0_\eps$ either to a critical point or to a persisting periodic orbit in $\M^+_\eps$ (Theorem \ref{T:slowflow}): all constructed 1-fronts originate from the uniform bare soil state. The existence of 1-front patterns that approach $(\bar{B},\bar{W}) = (0,\Psi/\Phi)$ as $t \to + \infty$ is embedded in these results through the application of the symmetry (\ref{eq:symmc}).
\end{remark}


\subsection{$W^u(P^0) \cap W^s(\M^+_\eps)$ and its touch down points on $\M^+_\eps$}
\label{sec:touchdown}

A (traveling) front pattern between the bare soil state $(0,\Psi/\Phi)$ and a (potentially stable) uniform vegetation state $(\bar{B}, \bar{W})$ of (\ref{eq:scaled}) corresponds to a heteroclinic solution $\gamma_{h}(\xi) = (w_{h}(\xi),p_{h}(\xi),b_{h}(\xi),q_{h}(\xi))$ of (\ref{eq:fast}) between the critical points $P^0=P^0_0 = (0,0,\Psi/\Phi,0) \in \M^0_\eps$ and $P^{+,s} = (b_+(w^{+,s}),0,w^{+,s},0) \in \M^+_\eps$ -- see section \ref{sec:crithom}. We know by Fenichels second Theorem that, by the normal hyperbolicity $\M^0_0$ and $\M^+_0$, their stable and unstable manifolds $W^{s,u}(\M^0_0)$ and $W^{s,u}(\M^+_0)$ persist as $W^{s,u}(\M^0_\eps)$ and $W^{s,u}(\M^+_\eps)$ for $\eps \neq 0$ as $w \in (1/(a + 1/4), 1/a)$ (\ref{condw0}), \cite{Jones1995,Kaper1999}. Thus, $\gamma_{h}(\xi) \subset W^u(P^0) \cap W^s(P^{+,s}) \subset W^u(\M^0_\eps) \cap W^s(\M^+_\eps)$ -- where we note that the manifolds $W^u(P^0)$ and $W^s(P^{+,s})$ are 2-dimensional, while $W^u(\M^0_\eps)$ and $W^s(\M^+_\eps)$ are 3-dimensional (and that the intersections take place in a 4-dimensional space).
\\ \\
We know by (\ref{WsWu0}) that $W^u(\M^0_0)$ and $W^s(\M^+_0)$ intersect transversely -- and thus that $W^u(\M^0_0) \cap W^s(\M^+_0)$ is 2-dimensional. Since $W^{u}(\M^0_\eps)$ and $W^{s}(\M^+_\eps)$ are $C^1$-$\mathcal{O}(\eps)$ close to $W^{u}(\M^0_0)$ and $W^{s}(\M^+_0)$, it immediately follows that $W^u(\M^0_\eps)$ and $W^s(\M^+_\eps)$ also intersect transversely, that $W^u(\M^0_\eps) \cap W^s(\M^+_\eps)$ is 2-dimensional and at leading order (in $\eps$) given by (\ref{WsWu0}). Since $W^u(P^0) \subset W^u(\M^0_\eps)$, $W^u(P^0)\cap W^s(\M^+_\eps)$ is a 1-dimensional subset of $W^u(\M^0_\eps) \cap W^s(\M^+_\eps)$ -- i.e. an orbit -- that follows $W^{u}(P^0)|_{\M^0_\eps} = \ell^{u}_\eps$ (\ref{deflsu}) exponentially close until its $w$-component reaches $w_h^+(c)$ (\ref{wafoc}) at which it `takes off' from $\M^0_\eps$ to follow the fast flow along the `parabolic' manifold given by (\ref{WsWu0}), all at leading order in $\eps$ -- see sections \ref{sec:fast}, \ref{sec:slowfull}. Since $w,q$ only vary slowly (\ref{eq:fast}), the $(w,q)$-components of the orbit $W^u(P^0)\cap W^s(\M^+_\eps)$ remain constant at leading order during its fast jump: it `touches down' on $\M^+_\eps$ with (at leading order) the same $(w,q)$-coordinates
(Remark \ref{rem:FenIII}). Therefore, we define the touch down curve $\T_{\rm down}(c) \subset \M^+_\eps$ as the set of touch down points of the orbits $W^u(P^0)\cap W^s(\M^+_\eps)$ that take off from $\M^0_\eps$ exponentially close to the intersection $\ell^u_\eps \cap \{w=w^+_h(c)\}$ (\ref{wafoc}), parameterized by $c$; it is at leading order (in $\eps$) given by,
\beq
\label{Tdown}
\T_{\rm down}(c) = \left\{\left(b_+(w^+_h(c)),0,w^+_h(c),\sqrt{\Phi}\left(w^+_h(c) - \frac{\Psi}{\Phi}\right)\right)\right\}
\feq
In terms of the projected $(w,q)$-coordinates by which the dynamics on $\M^+_\eps$ are described (\ref{eq:slow+eps}), $\T_{\rm down}(c)$ describes a smooth 1-dimensional manifold $\I_{\rm down} = \{(w_{\rm down}(c), q_{\rm down}(c))\}$ parameterized by $c$ with boundaries (its endpoints): the family of base points of the Fenichel fibers of $W^u(P^0) \cap W^s(\M^+_\eps)$ on $\M^+_\eps$ -- Remark \ref{rem:FenIII}; at leading order in $\eps$, $\I_{\rm down}$ is a straight interval with endpoints determined by the bounds (\ref{condw0}) on $w = w^+_h(c)$.
\begin{lemma}
\label{lem:Idown}
At leading order in $\eps$, $\I_{\rm down} = \left\{\left(w^+_h(c),\sqrt{\Phi}\left(w^+_h(c) - \frac{\Psi}{\Phi}\right)\right), c \in [-\frac{1}{\sqrt{2(1+4a)}},\frac{1}{\sqrt{2a}}]\right\}$ . The map $[-\frac{1}{\sqrt{2(1+4a)}},\frac{1}{\sqrt{2a}}] \to \I_{\rm down}$ is bijective and
\[
w^+_h\left(-\frac{1}{\sqrt{2(1+4a)}}\right) = \frac{4}{1+4a}, \; w^+_h\left( 0 \right) = \frac{9}{2+9a}, \; w^+_h\left(\frac{1}{\sqrt{2a}}\right) = \frac{1}{a}.
\]
\end{lemma}
Expression (\ref{wafoc}) a priori does not exclude the possibility that $w^+_h$ has several extremums as function of $c$, in fact $\frac{d}{dc} w^+_h(-\frac{1}{\sqrt{2(1+4a)}}) = 0$. The proof -- derivation -- of this lemma thus requires some careful, but straightforward, analysis. We refrain from going into the details here.
\\ \\
We conclude this section by noticing that heteroclinic connections $\gamma_{h}(\xi)$ between $P^0 \in \M^0_\eps$ and $P^{+,s} \in \M^+_\eps$ directly correspond to intersections $\I_{\rm down} \cap W^s(P^{+,s})|_{\M^+_\eps}$ (Remark \ref{rem:FenIII}). However, the coordinates of this intersection determine $c$ (through $\I_{\rm down}$), while $W^s(P^{+,s})|_{\M^+_\eps}$ also varies as function of $c$. Moreover, by the perturbed integrable nature of the flow on $\M^+_\eps$ (\ref{eq:slow+eps}), there can a priori be (countably) many intersections $\I_{\rm down} \cap W^s(P^{+,s})|_{\M^+_\eps}$. Thus, the analysis is more subtle and richer than (perhaps) expected -- as we shall see in the upcoming sections.

\begin{remark}
\label{rem:FenIII}
We (for instance) refer to \cite{DGK01} for a more careful treatment of `take off' and `touch down' points/manifolds. In fact, these points/manifolds correspond to base points of Fenichel fibers (that persist under perturbation by Fenichel's third Theorem \cite{Jones1995,Kaper1999}). By construction/definition, an orbit that touches down at a certain (touch down) point on a slow manifold is asymptotic to the orbit of the slow flow that has this point as initial condition. Therefore, if an orbit touches down on a stable manifold of a critical point on the slow manifold, it necessarily is asymptotic to this critical point.
\end{remark}

\subsection{Traveling 1-front patterns -- primary orbits}
\label{sec:primtrav1front}

\begin{figure}[t]
\centering
	\begin{minipage}{.45\textwidth}
		\centering
		\includegraphics[width =\linewidth]{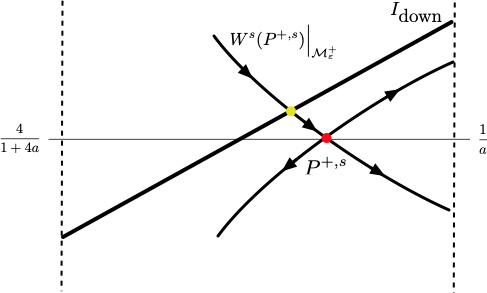}
	\end{minipage}%
	\hspace{1cm}
	\begin{minipage}{0.45\textwidth}
		\includegraphics[width=\linewidth]{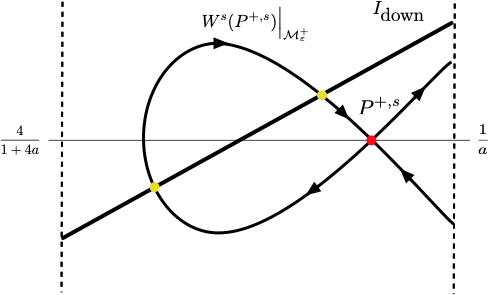}
	\end{minipage}
\caption{Sketches of the intersections of $\I_{\rm down}$ and $W^s(P^{+,s})|_{\M^+_0}$ in $\M^+_0$, i.e. the leading order configurations as described by the (integrable) slow reduced flow (\ref{eq:SRSplus}), in the 2 cases considered in Theorem \ref{T:primaryfronts}: there is one critical point $P^{+,s}$ of saddle type on $\M^+_0$ or there is a center $P^{+,c}$ and a saddle $P^{+,s}$ on $\M^+_0$.}
\label{fig:primhets}
\end{figure}

Our first result -- on the existence of {\it primary} heteroclinic orbits -- can be described in terms of the slow reduced flow on $\M^+_0$, or more precise, on intersections of the touch down manifold $\I_{\rm down}$ and the restricted stable manifold $W^s(P^{+,s})|_{\M^+_0} \subset \{ \mathcal{H}_{0}^+(w,q) = \mathcal{H}_{0}^{+,s}\}$ (\ref{manifoldsP+s}) of the reduced slow flow (\ref{eq:SRSplus}) on $\M^+_0$. However, it is a priori unclear whether such intersections may exist and how many of such intersections may occur: the many parameters of system (\ref{eq:fast}) have a `nontrivial' effect on $\I_{\rm down}$ and  $W^s(P^{+,s})|_{\M^+_0}$ and thus on their relative positions. To obtain a better insight in this, we `freeze' the flow of (\ref{eq:SRSplus}) by fixing $a, \A, \C, \D$ at certain values. Since $\B + a \Theta = \D + (a + \frac14) \A$ (\ref{defWD}), this indeed fixes all coefficients of the reduced slow flow (\ref{eq:SRSplus}) on $\M^+_0$. At the same time, this leaves a 1-parameter freedom in the parameters $\Phi, \Psi, \Omega, \Theta$. Defining,
\beq
\label{defchi}
\chi = \frac{1}{a}\left(\frac14 \A - \frac12 \C + \D\right),
\feq
we see that for all $\Phi$, the choices
\beq
\label{PsiasPhi}
\Psi = \frac{1}{a}\Phi - \chi, \; \Theta = \A - \frac{1}{a}\Phi + \chi, \; \Omega = \C - \A + \frac{1}{a}\Phi - \chi
\feq
yield identical slow reduced flows (\ref{eq:SRSplus}). On the other hand, the (leading order) interval $\I_{\rm down}$ clearly varies as function of $\Phi$,
\beq
\label{IdownPhi}
\I_{\rm down}(\Phi) = \left\{ q = \sqrt{\Phi} \left(w- \left( \frac{1}{a} - \frac{\chi}{\Phi}\right) \right), \; w \in \left(\frac{4}{1 + 4a}, \frac{1}{a}\right) \right\}.
\feq
Note that for $\chi > 0$, the intersection of $\I_{\rm down}$ with the $w$-axis can be varied between the critical $w$ values $4/(1 + 4a)$ and $1/a$ by increasing $\Phi$ from $a(1+4a) \chi$ to $\infty$. In fact, $\chi > 0$ necessarily holds in case there are 2 critical points on $\M^+_\eps$ (since in that case $\A, \D > 0$, $\C< 0$), while $\chi$ can also chosen to be positive in the case that there is only 1 critical point on $\M^+_\eps$. Thus, by choosing $\Psi, \Omega, \Theta$ as in (\ref{PsiasPhi}) and varying $\Phi$ we can control $\I_{\rm down} \cap W^s(P^{+,s})|_{\M^+_0}$ .
\begin{theorem}
\label{T:primaryfronts}
Let $P^{+,s} = (b_+(w^{+,s}),0,w^{+,s},0) \in \M^+_0$ be a critical point of (\ref{eq:fast}) that is a saddle point for the slow reduced flow (\ref{eq:SRSplus}) on $\M^+_0$, and consider the touch down manifold $\I_{\rm down}$ at leading order given in Lemma \ref{lem:Idown} and the restricted stable manifold $W^s(P^{+,s})|_{\M^+_0}$ of the reduced slow flow  (\ref{eq:SRSplus}). If there is a non-degenerate intersection point $(\bar{w}_{{\rm prim},0},\bar{q}_{{\rm prim},0}) \in  \I_{\rm down} \cap W^s(P^{+,s})|_{\M^+_0}$, then, for $\eps$ sufficiently small, there exists for $c = c_{\rm prim}$ a primary heteroclinic orbit $\gamma_{\rm prim}(\xi) = (w_{\rm prim}(\xi),p_{\rm prim}(\xi),b_{\rm prim}(\xi),q_{\rm prim}(\xi)) \subset W^u(P^0) \cap W^s(P^{+,s})$ of (\ref{eq:fast}) connecting $P^0 \in \M^0_\eps$ to $P^{+,s} \in \M^+_\eps$ -- where $c_{\rm prim} = c_{{\rm prim}, 0} + \mathcal{O}(\eps)$ and $c_{{\rm prim},0}$ is the unique solution of $w^+_h(c) = \bar{w}_{{\rm prim}, 0}$ (\ref{wafoc}). Departing from $P^0$ (and at leading order in $\eps$), $\gamma_{\rm prim}(\xi)$ first follows $\ell^u_0 \subset \M^0_0$ (\ref{deflsu}) until it reaches the take off point $(0,0,\bar{w}_{\rm prim},\bar{q}_{\rm prim})$ from which it jumps off from $\M^0_0$ and follows the fast flow along $W^u(\mathcal{M}_0^{0}) \cap W^{s}(\mathcal{M}_0^{+})$ (\ref{WsWu0}) to touch down on $\M^+_0$ at $(b_+(\bar{w}_{\rm prim}),0,\bar{w}_{\rm prim},\bar{q}_{\rm prim}) \in W^s(P^{+,s})|_{\M^+_0}$; from there, it follows $W^s(P^{+,s})|_{\M^+_0}$ towards $P^{+,s}$. Moreover,
\\
$\bullet$ if $P^{+,s}$ is the only critical point on $\M^+_\eps$, i.e. if $\C^2 - 4 \A \D > 0$, $\D < 0$, $\E > 0$ (\ref{saceM+0}), there is an open region $\mathcal{S}_{\rm s-prim}^{1}$ in $(a,\Psi,\Phi,\Omega,\Theta)$ parameter space for which $\I_{\rm down}$ and  $W^s(P^{+,s})|_{\M^+_0}$ intersect transversely; however, there is at most one intersection $(\bar{w}_{\rm prim},\bar{q}_{\rm prim}) \in  \I_{\rm down} \cap W^s(P^{+,s})|_{\M^+_0}$ and thus at most one primary heteroclinic orbit $\gamma_{\rm prim}(\xi)$; in fact, this is the only possible heteroclinic orbit between $P^0$ and $P^{+,s}$;
\\
$\bullet$ if there are two critical points on $\M^+_\eps$, the center $P^{+,c}$ and saddle $P^{+,s}$, i.e. if $\C^2 - 4 \A \D > 0$, $\C < 0$, $\D > 0$, then there are open regions $\mathcal{S}_{\rm cs-prim}^{1}$, respectively $\mathcal{S}_{\rm cs-prim}^{2}$, in $(a,\Psi,\Phi,\Omega,\Theta)$ parameter space for which $\I_{\rm down}$ and $W^s(P^{+,s})|_{\M^+_0}$ have 1, resp. 2, (transversal) intersections, so that there can be (up to) 2 distinct primary heteroclinic orbit $\gamma_{\rm prim}^j(\xi)$ that travel with different speeds, i.e. $c_{\rm prim}^2 < c_{\rm prim}^1$.
\\
A primary heteroclinic orbit $\gamma_{\rm prim}(\xi) = (w_{\rm prim}(\xi),p_{\rm prim}(\xi),b_{\rm prim}(\xi),q_{\rm prim}(\xi))$ corresponds to a (localized, traveling, invasion) 1-front pattern $(B(x,t), W(x,t)) = (b_{\rm prim}(x-c_{\rm prim} t), w_{\rm prim}(x-c_{\rm prim} t))$ in PDE (\ref{eq:scaled}) that connects the bare soil state $(\bar{B},\bar{W})=(0,\Psi/\Phi)$ to the uniform vegetation state $(\bar{B},\bar{W}) = (b_+(w^{+,s}),w^{+,s})$.
\end{theorem}
In the case of 2 critical points on $\M^+_\eps$, we shall see that the primary orbits may only be the first of many `higher order' heteroclinic orbits -- see section \ref{sec:frontsbyM+eps}. We refer to Fig. \ref{fig:primhets} for sketches of the constructions in $\M^0_\eps$ that yield the primary heteroclinic orbits $\gamma_{\rm prim}(\xi)$ and to Figs. \ref{fig:basicpatterns}a, \ref{fig:numerics_front_ehud_both} and \ref{fig:numerics_1_front_traveling} for the associated -- numerically obtained -- primary 1-front patterns in (\ref{eq:scaled})-- see especially Fig. \ref{fig:numerics_1_front_ehud_phase_space} in which the the slow-fast-slow structure of a (numerically obtained) heteroclinic front solutions of (\ref{eq:scaled}) is exhibited by its projection in the 3-dimensional $(b,w,q)$-subspace of the 4-dimensional phase space associated to (\ref{eq:fast}).
\\

{\bf Proof.} The existence of the heteroclinic orbit $\gamma_{\rm prim}(\xi)$ follows by construction -- Remark \ref{rem:FenIII} -- from an intersection of $\I_{\rm down}$ and $W^s(P^{+,s})|_{\M^+_\eps}$. Thus, we first need to show that a (non-degenerate) intersection $\I_{\rm down} \cap W^s(P^{+,s})|_{\M^+_0}$ implies an intersection $\I_{\rm down} \cap W^s(P^{+,s})|_{\M^+_\eps}$. More precise, since $W^s(P^{+,s})|_{\M^+_\eps}$ varies with $c$, i.e. since $W^s(P^{+,s})|_{\M^+_\eps} = W^s(P^{+,s})|_{\M^+_\eps}(c)$, we need to determine $c^*$ such that $W^s(P^{+,s})|_{\M^+_\eps}(c^*)$ intersects $\I_{\rm down} = \{(\bar{w}_{\rm down}(c), \bar{q}_{\rm down}(c))\}$ exactly at $(\bar{w}_{\rm down}(c^*), \bar{q}_{\rm down}(c^*))$.
\\ \\
By the assumption that $(\bar{w}_{{\rm prim},0},\bar{q}_{{\rm prim},0}) \in \I_{\rm down} \cap W^s(P^{+,s})|_{\M^+_0}$ is a non-degenerate intersection point, we know that the intersection is transversal, and thus that $W^s(P^{+,s})|_{\M^+_\eps}(\tilde{c})$ -- i.e. $W^s(P^{+,s})|_{\M^+_\eps}$ for (\ref{eq:slow+eps}) with $c = \tilde{c}$ -- also intersects $\I_{\rm down}$ transversally as $\tilde{c}$ is varied around $c_{{\rm prim},0}$ in an $\mathcal{O}(1)$ fashion. Thus, for $\tilde{c}$ sufficiently (but $\mathcal{O}(1)$) close to $c_{{\rm prim}, 0}$, $\I_{\rm down} \cap W^s(P^0)|_{\M^+_\eps}(\tilde{c}) = (\bar{w}_{\rm down}(c_i),\bar{q}_{\rm down}(c_i))$ determines a curve $c_i = c_i(\tilde{c})$ by $c_i = \bar{w}_{\rm down}(c_i)$. Since the flows of (\ref{eq:SRSplus}) and (\ref{eq:slow+eps}) are $\mathcal{O}(\eps)$ close, we know that $\|(\bar{w}_{{\rm prim},0},\bar{q}_{{\rm prim},0}) - (\bar{w}_{\rm down}(c_i),\bar{q}_{\rm down}(c_i))\| = \mathcal{O}(\eps)$, which implies that $c_{i}(\tilde{c}) = c_{{\rm prim},0} + \mathcal{O}(\eps)$. Hence, the $\mathcal{O}(1)$ variation of $\tilde{c}$ through $c_{{\rm prim},0}$ yields at leading order (in $\eps$) a horizontal line $c_{i}(\tilde{c}) \equiv c_{{\rm prim},0}$: there must be a unique intersection $c_{i}(\tilde{c}^*) = \tilde{c}^* $, and thus, by construction, $\I_{\rm down} \cap W^s(P^{+,s})|_{\M^+_\eps}(\tilde{c}^*) = (\bar{w}_{\rm down}(\tilde{c}^*), \bar{q}_{\rm down}(\tilde{c}^*))$: $\tilde{c}^* = c_{\rm prim}$.
\\ \\
If $P^{+,s}$ is the only critical point on $\M^+_\eps$ -- i.e. if $\C^2 - 4 \A \D > 0$, $\D < 0$, $\E > 0$ -- we freeze the flow of (\ref{eq:SRSplus}) with $\A, \C, \D$ such that $\chi > 0$ (\ref{defchi}) and define $\Phi = \Phi^{+,s}$ such that $\Psi(\Phi)/\Phi = 1/a - \chi/\Phi = w^{+,s}$, the $w$-coordinate of the saddle $P^{+,s}$ on $\M^+_0$ -- see (\ref{PsiasPhi}), (\ref{IdownPhi}). Since $q$ is an increasing function of $w$ on $\I_{\rm down}$ and $W^s(P^{+,s})|_{\M^+_\eps}$ is decreasing near $P^{+,s}$ -- see Fig \ref{fig:primhets}a -- it follows that there must be a transversal intersection $\I_{\rm down} \cap W^s(P^{+,s})|_{\M^+_0}$ for values of $\Phi$ in an (open) interval around $\Phi^{+,s}$. Transversality implies that the intersection persists under varying $\A, \C, \D$ around their initially frozen values, which establishes the existence of the open region $\mathcal{S}_{\rm s-prim}^{1}$ in $(a,\Psi,\Phi,\Omega,\Theta)$-space for which $\I_{\rm down}$ and $W^s(P^{+,s})|_{\M^+_0}$ intersect. Moreover, the manifold $W^s(P^{+,s})|_{\M^+_0} \subset \{ \mathcal{H}_{0}^+(w,q) = \mathcal{H}_{0}^{+,s}\}$ (\ref{manifoldsP+s}) is given by a (strictly) decreasing function $q^{+,s}|_{\M^+_0}(w)$ for all $w \in (4/(1 + 4a),1/a)$ since it cannot have extremums: zeroes of $\frac{d}{dw} q^{+,s}|_{\M^+_0}(w)$ correspond to zeroes of $\frac{\partial}{\partial w} \mathcal{H}_{0}^+(w,q)$ (\ref{eq:HamiltonianSRS}) and thus to critical points of (\ref{eq:SRSplus}). By assumption, there are no critical points besides $P^{+,s}$, which yields that there indeed can be maximally one intersection $\I_{\rm down} \cap W^s(P^{+,s})|_{\M^+_0}$.
\\ \\
To control the case with a center $P^{+,c}$ and saddle $P^{+,s}$ on $\M^+_0$, we again consider the unfolded saddle-node case of Lemma \ref{lem:posigmasmal} and define $\Phi = \Phi^{+,c}$ such that $\Psi(\Phi)/\Phi = 1/a - \chi/\Phi = w^{+,c}$, the $w$-coordinate of the center $P^{+,c}$. The level set $\{ \mathcal{H}_{0}^+(w,q) = \mathcal{H}_{0}^{+,s}\}$ forms a small (w.r.t. the unfolding parameter $\sigma$) homoclinic loop around $P^{+,c}$ that intersects $\I_{\rm down}$ (transversally) in two points $(\bar{w}_{\rm prim}^{j,0},\bar{q}_{\rm prim}^{j,0})$, $j = 1,2$ -- see Fig. \ref{fig:primhets}(b). By varying $\Phi$ around $\Phi = \Phi^{+,c}$ and $\A, \C, \D$ around their initially frozen values, we find the open region $\mathcal{S}_{\rm cs-prim}^{2}$ in $(a,\Psi,\Phi,\Omega,\Theta)$-space for which both elements of the intersection $\I_{\rm down} \cap W^s(P^{+,s})|_{\M^+_0}$ persist: for $(a,\Psi,\Phi,\Omega,\Theta) \in \mathcal{S}_{\rm cs-prim}^{2}$, (\ref{eq:fast}) has 2 (distinct) primary heteroclinic orbits $\gamma_{\rm prim}^j(\xi)$, $j=1,2$, that correspond to 1-front patterns traveling with speeds $c_{\rm prim}^1 \neq c_{\rm prim}^2$ -- where $c_{\rm prim}^{j,0}$ is the unique solution of $w^+_h(c) = \bar{w}_{\rm prim}^{j,0}$. Finally, we note that the existence of the open set $\mathcal{S}_{\rm cs-prim}^{1}$ follows by considering $\I_{\rm down} \cap W^s(P^{+,s})|_{\M^+_0}$ for values of $\Phi > \Phi^{+,s}$ (as defined above). \hfill $\Box$

\subsection{Traveling 1-front patterns by the perturbed integrable flow on $\M^+_\eps$}
\label{sec:frontsbyM+eps}

\begin{figure}[t]
  \centering
	\includegraphics[width=12cm]{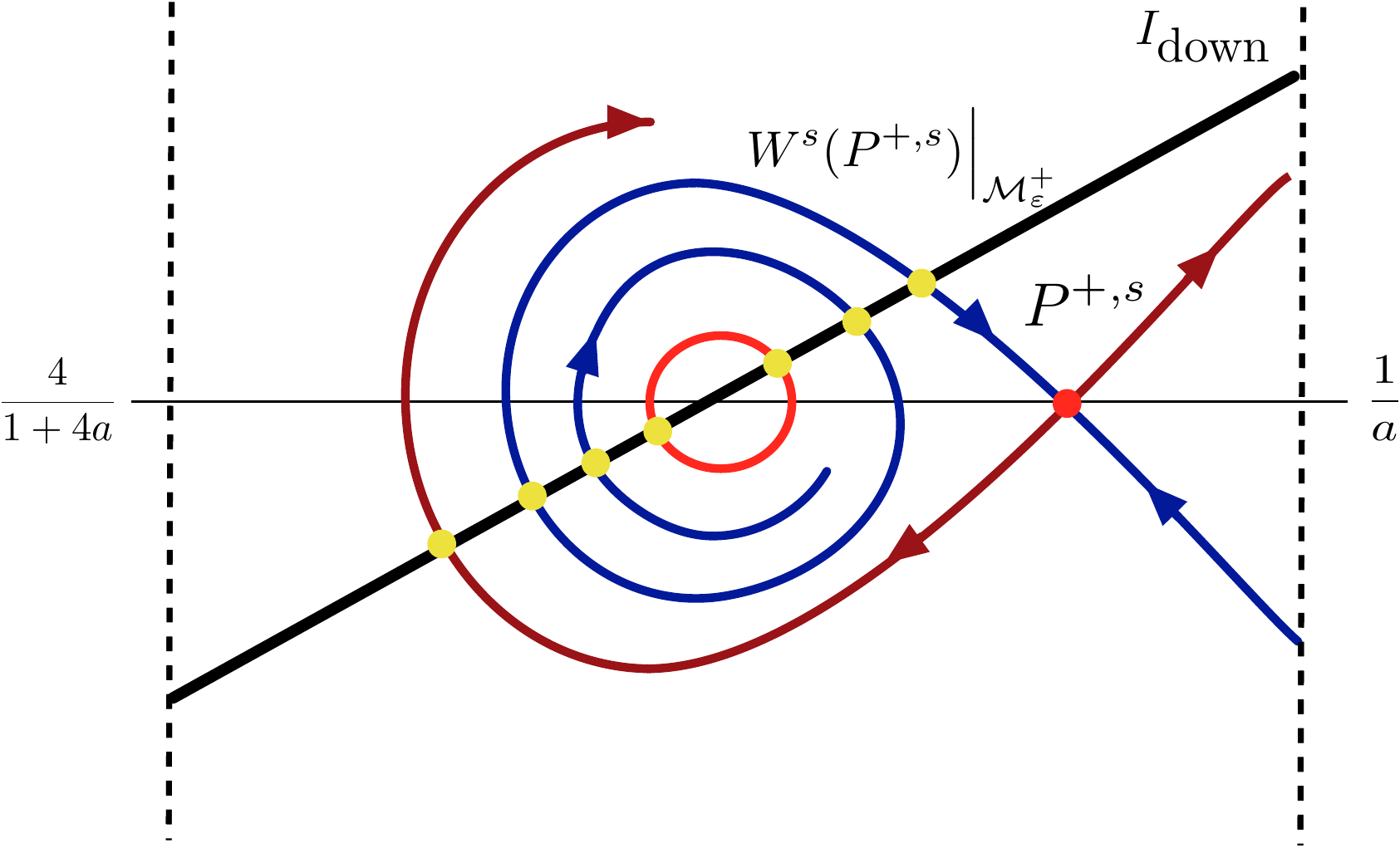}
     \caption{A sketch of the flow (\ref{eq:slow+eps}) on $\M^+_\eps$ for $(a,\Psi,\Phi,\Omega,\Theta) \in \mathcal{S}_{\rm h-p}$ (Theorem \ref{T:frontstoperiodic}) and $c = c_{\rm prim}^1$ (Theorem \ref{T:primaryfronts}) in the case that (\ref{DeltaHpos}). Since $W^s(P^{+,s})|_{\M^+_\eps}$ `wraps around' the persistent periodic solution $(b_{p,\eps}(X),p_{p,\eps}(X),w_{p,\eps}(X),q_{p,\eps}(X))$ (Theorem \ref{T:frontstoperiodic}) -- in backwards time -- and since $\I_{\rm down}$ intersects this orbit in 2 points (by assumption), there are two countable sets of intersections $\I_{\rm down} \cap W^s(P^{+,s})|_{\M^+_\eps}$.}
 	\label{fig:counthets}
\end{figure}

As in Theorem \ref{T:slowflow}, we assume throughout this section that there is a center $P^{+,c}= (b_+(w^{+,c}),0,w^{+,c},0)$ and a saddle $P^{+,s}= (b_+(w^{+,s}),0,w^{+,s},0)$ on $\M^+_\eps$ and -- for simplicity -- that the unperturbed homoclinic orbit $(w_{{\rm hom},0}(X),q_{{\rm hom},0}(X))$ to $P^{+,s}$ of (\ref{eq:SRSplus}) on $\M^0_0$ -- that is a subset of $W^s(P^{+,s})|_{\M^+_0} \subset \{ \mathcal{H}_{0}^+(w,q) = \mathcal{H}_{0}^{+,s}\}$ -- lies entirely in the $w$-region in which both $\M^0_0$ and $\M^+_0$ are normally hyperbolic, i.e. we assume that (\ref{asshomM+eps}) holds.
\\ \\
The homoclinic orbit $(w_{{\rm hom},0}(X),q_{{\rm hom},0}(X))$ of (\ref{eq:SRSplus}) typically breaks open under the perturbed flow of (\ref{eq:slow+eps}), and $W^s(P^{+,s})|_{\M^+_\eps}$ either spirals inwards in backwards `time', i.e. as $\xi \to -\infty$, or not. In the former case, there will be (typically many) further intersections $\I_{\rm down} \cap W^s(P^{+,s})|_{\M^+_0}$ -- see Fig. \ref{fig:counthets}. Of course, this is determined by the sign of $\Delta H_{\rm hom}$ (\ref{persslowhom}): if
\beq
\label{DeltaHpos}
\Delta H_{\rm hom} = c \int_{\underline{w}_{h,0}}^{w^{+,s}} \rho_1(w)
\sqrt{2 \mathcal{H}_{0}^{+,s} - 2 \mathcal{A}w + \left(\mathcal{B}+a\Theta\right)w^2 + 2 \C \mathcal{J}_{0}^+(w)} \, dw > 0
\feq
(at leading order in $\eps$), we may expect further heteroclinic connections $\gamma_{h, j}$ in (\ref{eq:fast}) connecting $P^0 \in \M^0_\eps$ to $P^{+,s} \in \M^+_\eps$ beyond the primary orbits $\gamma_{\rm prim}(\xi)$ established in Theorem \ref{T:primaryfronts}. In fact, it follows directly that $\gamma_{\rm prim}^1(\xi)$ and $\gamma_{\rm prim}^2(\xi)$ are the only heteroclinic orbits between $P^0$ and $P^{+,s}$ if (\ref{DeltaHpos}) does not hold. If (\ref{DeltaHpos}) does hold, the (spiraling part of) $W^s(P^{+,s})|_{\M^+_0}$ clearly must limit -- for $\xi \to -\infty$ -- on either the center $P^{+,c}$ or, if $(a,\Psi,\Phi,\Omega,\Theta) \in \mathcal{S}_{\rm per}$, on the persistent periodic solution $(b_{p,\eps}(X),p_{p,\eps}(X),w_{p,\eps}(X),q_{p,\eps}(X)) \subset \M^+_\eps$ (Theorem \ref{T:slowflow}). Therefore, we first formulate a result on the existence of heteroclinic connections between $P^{0} \in \M^0_\eps$ and $(b_{p,\eps}(X),p_{p,\eps}(X),w_{p,\eps}(X),q_{p,\eps}(X)) \subset \M^+_\eps$. Like in Theorem \ref{T:primaryfronts}, this can be done in terms of the unperturbed flow in $\M^+_0$.
\begin{theorem}
\label{T:frontstoperiodic}
Assume that (\ref{asshomM+eps}) holds and that $(a,\Psi,\Phi,\Omega,\Theta) \in \mathcal{S}_{\rm per}$. Let $(w_{p,0}(X),q_{p,0}(X)) \subset \{ \mathcal{H}_{0}^+(w,q) = \mathcal{H}_{p,0}^+ \}$ with $\mathcal{H}_{p,0}^+ \in (\mathcal{H}_{0}^{+,c}, \mathcal{H}_{0}^{+,s})$ (\ref{eq:HamiltonianSRS}) be the periodic solution of (\ref{eq:SRSplus}) that persists (on $\M^+_\eps$) as periodic solution $(b_{p,\eps}(X), p_{p,\eps}(X), w_{p,\eps}(X), q_{p,\eps}(X))$ on $\M^+_\eps$ of (\ref{eq:fast}). Then there is an open set $\mathcal{S}_{\rm h-p} \subset \mathcal{S}_{\rm per} \cap \mathcal{S}_{\rm cs-prim}^{2}$ -- with $\mathcal{S}_{\rm cs-prim}^{2}$ defined in Theorem \ref{T:primaryfronts} -- such that there are 2 (non-degenerate) intersection points $(\bar{w}_{\rm h-p}^j,\bar{q}_{\rm h-p}^j) \in \I_{\rm down} \cap \{ \mathcal{H}_{0}^+(w,q) = \mathcal{H}_{p,0}^+ \}$, $j=1,2$, that correspond -- for $\eps$ sufficiently small -- to 2 distinct heteroclinic orbits $\gamma_{\rm h-p}^j(\xi)= (b_{\rm h-p}^j(\xi),p_{\rm h-p}^j(\xi),w_{\rm h-p}^j(\xi),q_{\rm h-p}^j(\xi))$ of (\ref{eq:fast}) -- in which $c = c_{\rm h-p}^j$ -- between the critical point $P^0 \in \M^0_\eps$ and the periodic orbit $(b_{p,\eps}(X),p_{p,\eps}(X),w_{p,\eps}(X),q_{p,\eps}(X)) \subset \M^+_\eps$; at leading order in $\eps$, $c_{\rm h-p}^j$ is determined by $w_h^+(c) = \bar{w}_{\rm h-p}^j$, with $c_{\rm prim}^2 < c_{\rm h-p}^2 < c_{\rm h-p}^1 <c_{\rm prim}^1$ (Theorem \ref{T:primaryfronts}).
\\
The orbits $\gamma_{\rm h-p}^j(\xi)$ correspond traveling 1-front patterns $(B(x,t), W(x,t)) = (b_{\rm h-p}^j(x-c_{\rm h-p}^j t),w_{\rm h-p}^j(x-c_{\rm h-p}^j t))$ in PDE (\ref{eq:scaled}) that connect the bare soil state $(\bar{B},\bar{W})=(0,\Psi/\Phi)$ to the traveling wave train $(B_{p,\eps}(\eps(x-c_{\rm h-p}^j t)), W_{p,\eps}(\eps(x-c_{\rm h-p}^j t))$ of Corollary \ref{cor:slowflow}.
\end{theorem}
Notice that this result is independent of condition (\ref{DeltaHpos}), i.e. Theorem \ref{T:frontstoperiodic} holds independent of the sign of $\Delta H_{\rm hom}$. Moreover, we could formulate similar limiting result concerning heteroclinic 1-front connections between $P^0 \in \M^0_\eps$ and $P^{+,c} \in \M^+_\eps$ for $(a,\Psi,\Phi,\Omega,\Theta)$ on a certain co-dimension 1 manifold. Since the background state associated to $P^{+,c}$ cannot be stable -- section \ref{sec:crithom} -- we refrain from going into the details.
\\ \\
{\bf Proof.} The proof goes exactly along the lines of that of Theorem \ref{T:primaryfronts}. \hfill $\Box$
\\ \\
Theorem \ref{T:frontstoperiodic} provides the foundation for a result on the existence of multiple -- in fact countably many -- distinct traveling 1-front connections between $P^0 \in \M^0_\eps$ and $P^{+,j} \in \M^0_\eps$ for an open set in parameter space -- see also the sketches in Figs. \ref{fig:sketchessection1}a and \ref{fig:sketchessection1}b.
\begin{theorem}
\label{T:count1fronts}
Assume that the conditions of Theorem \ref{T:frontstoperiodic} hold and let $(a,\Psi,\Phi,\Omega,\Theta) \in \mathcal{S}_{\rm h-p}$. If $c_{\rm h-p}^j$ and $c_{\rm prim}^j$ have the same sign (for either $j=1$ or 2) and if (\ref{DeltaHpos}) holds for $c$ of this sign, then -- for $\eps$ sufficiently small -- there are countably many distinct heteroclinic orbits $\gamma_{h}^{j,k}(\xi)= (b_{h}^{j,k}(\xi),p_{h}^{j,k}(\xi),w_{h}^{j,k}(\xi),q_{h}^{j,k}(\xi))$, $k \geq 0$, of (\ref{eq:fast}) with $c = c_{h}^{j,k}$ connecting $P^0 \in \M^0_\eps$ to $P^{+,s} \in \M^+_\eps$. Moreover, $\gamma_{h}^{j,0}(\xi) = \gamma_{\rm prim}^j(\xi)$, $|c_{h}^{1,k+1} - c_{h}^{1,k}| = \mathcal{O}(\eps)$, and,
\[
\begin{array}{lll}
j=1: & c_{\rm h-p}^1 < ... < c_{h}^{1,k} < ... < c_{h}^{1,1} < c_{h}^{1,0} = c_{\rm prim}^1, & c_{h}^{1,k} \downarrow c_{\rm h-p}^1 \; \; {\rm for} \; \; k \to \infty,
\\
j=2: & c_{\rm prim}^2 = c_{h}^{2,0} < c_{h}^{2,1} < ... < c_{h}^{2,k} < ... < c_{\rm h-p}^2, & c_{h}^{2,k} \uparrow c_{\rm h-p}^2 \; \; {\rm for} \; \; k \to \infty.
\end{array}
\]
Each orbit $\gamma_{h}^{j,k}(\xi)$ corresponds to a (localized, traveling, invasion) 1-front pattern $(B(x,t), W(x,t)) = (b_{h}^{k,j}(x-c_{h}^{k,j}t), w_{h}^{k,j}(x-c_{h}^{k,j}t))$ in PDE (\ref{eq:scaled}) that connects the bare soil state $(\bar{B},\bar{W})=(0,\Psi/\Phi)$ to the uniform vegetation state $(\bar{B},\bar{W}) = (b_+(w^{+,s}),w^{+,s})$.
\end{theorem}

\begin{figure}[t]
	\centering
	\includegraphics[width=15cm]{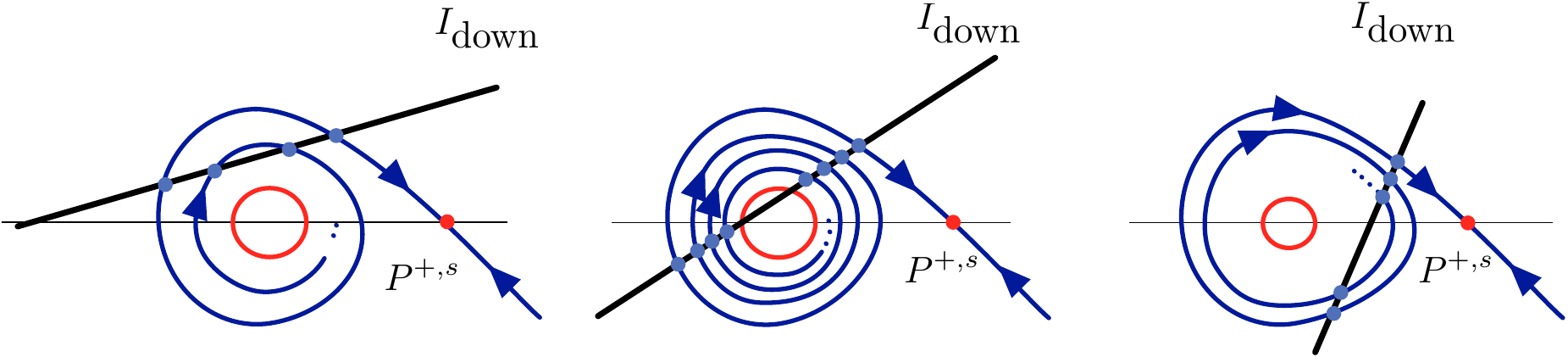}
	\caption{Sketches of 3 relative configurations of $\I_{\rm down}(\Phi)$ and $W^s(P^{+,s})|_{\M^+_0}$ for increasing $\Phi$ that represent 3 distinct stages in the bifurcation scenario of Fig. \ref{fig:cPhi2}.}
	\label{fig:spiralaroundperiodic-3config}
\end{figure}

\begin{figure}[t]
	\centering
	\includegraphics[width=15cm]{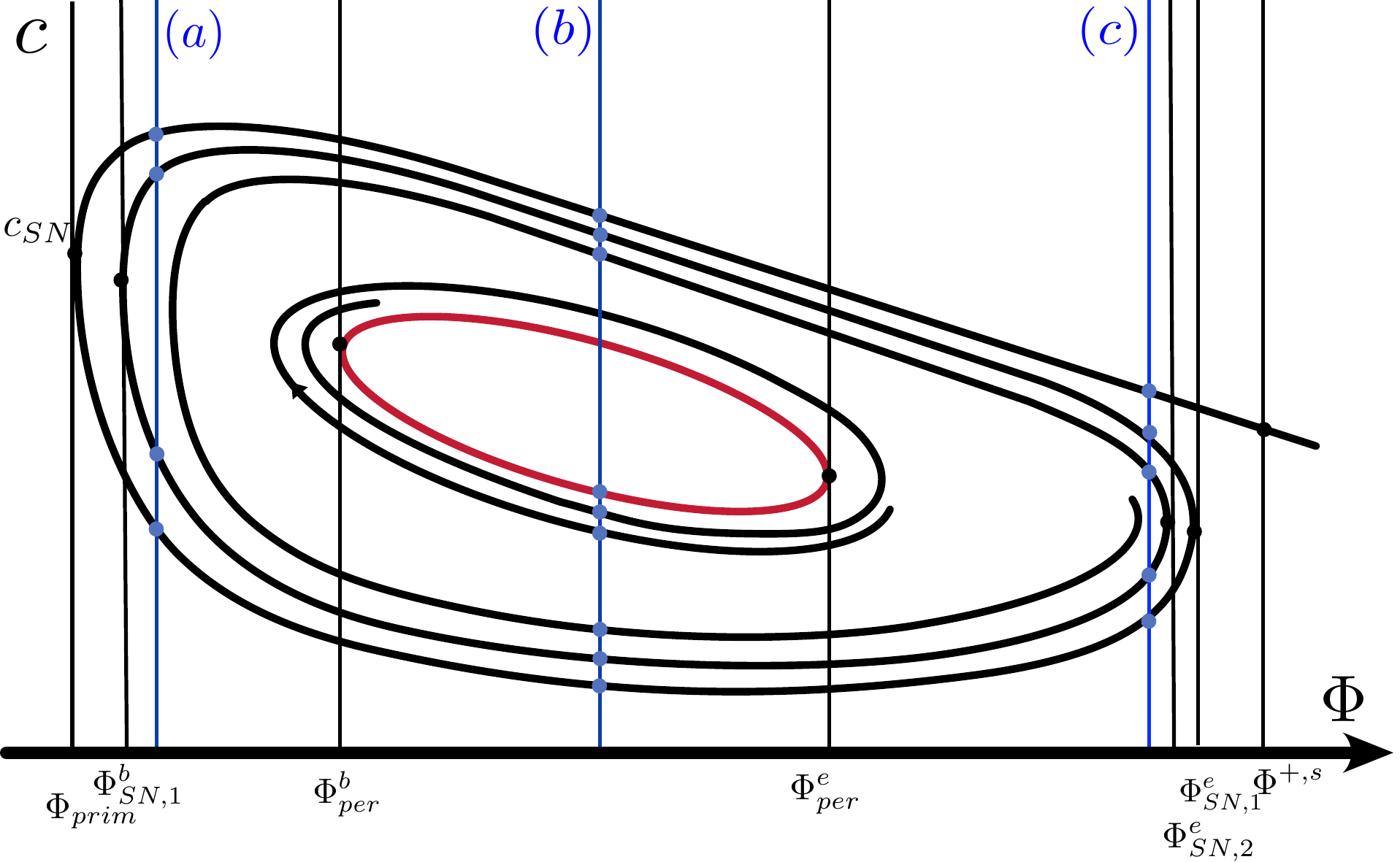}
	\caption{A sketch of the bifurcation scenario as function of $\Phi$ representing the appearance in a saddle-node bifurcation of the 2 primary heteroclinic 1-front orbits $\gamma_{\rm prim}^j(\xi)$, $j=1,2$, of Theorem \ref{T:primaryfronts}, followed by further saddle-node bifurcations leading to the situation governed by Theorem \ref{T:count1fronts} in which countably many 1-front orbits exist; these orbits subsequently disappear in another cascade of saddle-node bifurcations eventually leaving only one (primary) 1-front orbit behind (Theorem \ref{T:primaryfronts}). The relative configurations of $\I_{\rm down}(\Phi)$ and $W^s(P^{+,s})|_{\M^+_0}$ sketched in Fig. \ref{fig:spiralaroundperiodic-3config} occur at the $\Phi$-values indicated by the vertical (a), (b) and (c) lines.}
	\label{fig:cPhi2}
\end{figure}

As in the proofs of Theorems \ref{T:slowflow} and \ref{T:primaryfronts}, we can verify that there indeed are open regions in $(a,\Psi,\Phi,\Omega,\Theta)$-space for which $c_{\rm h-p}^j$ and $c_{\rm prim}^j$ have the same sign (for either $j=1, 2$ or for both) and such that (\ref{DeltaHpos}) holds, by considering the unfolded saddle-node case of Lemma \ref{lem:posigmasmal}. In fact, we know from Lemma \ref{lem:Idown} that $c$ changes sign as the $w$-coordinate of the intersection point on $\I_{\rm down} \cap W^s(P^{+,s})|_{\M^+_0}$ passes through $9/(2+9a)$. Thus (and for instance), all 4 values $c_{\rm h-p}^j$ and $c_{\rm prim}^j$, $j=1,2$, must have the same sign as the entire homoclinic orbit spanned by $W^s(P^{+,s})|_{\M^+_0}$ either is to the left or to the right of $w=9/(2+9a)$ -- more precise, if either $(\underline{w}_{h,0}, w^{+,s}) \subset (4/(1+4a),9/(2+9a))$ or $(\underline{w}_{h,0}, w^{+,s}) \subset (9/(2+9a),1/a)$ (cf. (\ref{asshomM+eps})). Note that it follows from (\ref{defPSNwSN}) that $w^{SN}_0 = 9/(2+9a)$ implies that $\C^2 = \A^2/9$ (independent of $a$), so that we can indeed move the homoclinic loop associated to the unfolded saddle-node -- i.e. $\sigma \ll 1$ as in (\ref{unfSN}) -- through $w = 9/(2+9a)$ by increasing $\C^2 \in (0,\A^2)$ through $\A^2/9$. On the other hand, it is certainly also possible that $c_{\rm h-p}^j$ and $c_{\rm prim}^j$ do not have the same sign. Hence, apart from the PDE point of view -- from which it is natural to consider stationary patterns -- this gives us an additional motivation to study the sign-changing stationary case $c=0$ in more detail, as we will briefly do in Remark \ref{rem:counthetsc=0} in section \ref{sec:stand1fronts}.
\\ \\
{\bf Proof.} We only consider the case $j =1$, i.e. we assume that $c_{\rm h-p}^1$ and $c_{\rm prim}^1$ have the same sign and that (\ref{DeltaHpos}) holds for $c = c_{\rm h-p}^1, c_{\rm prim}^1$. The proof for $j=2$ goes exactly along the same lines.
\\ \\
For $c=c_{\rm prim}^1$, $W^s(P^{+,s})|_{\M^+_\eps}$ by assumption spirals inwards in backwards `time' and `wraps around' the (perturbed) periodic orbit $(w_{p,\eps}(X),q_{p,\eps}(X))$ on (the projection of) $\M^+_\eps$ -- see Fig. \ref{fig:counthets}. Since $(a,\Psi,\Phi,\Omega,\Theta) \in \mathcal{S}_{\rm h-p}$, $W^s(P^{+,s})|_{\M^+_\eps}$ must intersect $\I_{\rm down}$ countably many times. We define $(\bar{w}_{i}^{1,1}, \bar{q}_{i}^{1,1})$ as the next intersection of $W^s(P^{+,s})|_{\M^+_\eps}(c_{\rm prim}^1)$ with $\I_{\rm down}$ beyond the 2 primary intersection points: it is the first non-primary intersection point and has $\bar{q}_{i}^{1,1} > 0$. As before, $(\bar{w}_{i}^{1,1}, \bar{q}_{i}^{1,1}) \in \I_{\rm down} = \{(\bar{w}_{\rm down}(c), \bar{q}_{\rm down}(c))\}$ determines the value $c_{i}^{1,1}$ through $\bar{w}_{\rm down}(c) = \bar{w}_{i}^{1,1}$ -- where we know that $c_{i}^{1,1} < c_{\rm prim}^1$ since the $w$-component of $\I_{\rm down}$ is a monotonically increasing function of $q$ (Lemma \ref{lem:Idown}). Since the perturbation term in (\ref{eq:slow+eps}) is $\mathcal{O}(\eps)$, it follows that $\|(\bar{w}_{\rm prim}^1,\bar{q}_{\rm prim}^1) - (\bar{w}_{i}^{1,1}, \bar{q}_{i}^{1,1})\| =  \mathcal{O}(\eps)$ and thus that $c_{\rm prim}^1 - c_{i}^{1,1} = \mathcal{O}(\eps)$. An $\mathcal{O}(\eps)$ change in $c$ yields an $\mathcal{O}(\eps^2)$ change in the flow of (\ref{eq:slow+eps}), hence for all $\tilde{c}$ $\mathcal{O}(\eps)$ close to $c_{\rm prim}^1$, the first non-primary intersection of $W^s(P^{+,s})|_{\M^+_\eps}(\tilde{c})$ and $\I_{\rm down}$ -- denoted by $(\bar{w}_{i}^{1,1}(\tilde{c}), \bar{q}_{i}^{1,1}(\tilde{c}))$ -- must be $\mathcal{O}(\eps^2)$ close to $(\bar{w}_{i}^{1,1}, \bar{q}_{i}^{1,1}) \in \I_{\rm down} \cap W^s(P^{+,s})|_{\M^+_\eps}(c_{\rm prim}^1)$. Thus, the speed $c_{i}^{1,1}(\tilde{c})$ associated to this intersection -- by $\bar{w}_{\rm down}(c) = \bar{w}_{i}^{1,1}(\tilde{c})$ -- must also be $\mathcal{O}(\eps^2)$ close to $c_{i}^{1,1}$. The situation is therefore similar to that in the proof of Theorem \ref{T:primaryfronts}: an $\mathcal{O}(\eps)$ variation of $\tilde{c}$ around $c_{i}^{1,1}$ in (\ref{eq:slow+eps}) yields only an $\mathcal{O}(\eps^2)$ change in the $c$-coordinate associated the first non-primary intersection $\I_{\rm down} \cap W^s(P^{+,s})|_{\M^+_\eps}(\tilde{c})$ so that there must be an unique $\tilde{c} = c_{h}^{1,1}$ such that $\I_{\rm down} \cap W^s(P^{+,s})|_{\M^+_\eps}(c_{h}^{1,1}) = (\bar{w}_{\rm down}(c_{h}^{1,1}), \bar{q}_{\rm down}(c_{h}^{1,1}))$. This establishes the existence of the first non-primary heteroclinic 1-front orbit $\gamma_{h}^{1,1}(\xi)$ for $c = c_{h}^{1,1}$ in (\ref{eq:fast}).
\\ \\
We can now iteratively consider the first intersection in backwards `time' -- denoted by $(\bar{w}_{i}^{1,2}, \bar{q}_{i}^{1,2})$ -- of $W^s(P^{+,s})|_{\M^+_\eps}(c_{h}^{1,1})$ with $\I_{\rm down}$ beyond $(\bar{w}_{\rm down}(c_{h}^{1,1}), \bar{q}_{\rm down}(c_{h}^{1,1}))$ with $\bar{q}_{i}^{1,2} > 0$ -- so that the speed $c_{i}^{1,2}$ associated to this intersection is $\mathcal{O}(\eps)$ close to $c_{h}^{1,1}$. Completely analogous to the above arguments, we deduce the existence of an unique $c=c_{h}^{1,2}$ such that $\I_{\rm down} \cap W^s(P^{+,s})|_{\M^+_\eps}(c_{h}^{1,2}) = (\bar{w}_{\rm down}(c_{h}^{1,2}), \bar{q}_{\rm down}(c_{h}^{1,2}))$, which establishes the existence of the next non-primary 1-front orbit $\gamma_{h}^{1,2}(\xi)$ of (\ref{eq:fast}) with $0 < c_{h}^{1,1} - c_{h}^{1,2} = \mathcal{O}(\eps)$. Next, we construct $\gamma_{h}^{1,3}(\xi)$ in (\ref{eq:fast}) with $0 < c_{h}^{1,2} - c_{h}^{1,3} = \mathcal{O}(\eps)$ through the intersection $\I_{\rm down} \cap W^s(P^{+,s})|_{\M^+_\eps}(c_{h}^{1,3}) = (\bar{w}_{\rm down}(c_{h}^{1,3}), \bar{q}_{\rm down}(c_{h}^{1,3}))$, etc.
\\ \\
Theorem \ref{T:slowflow} holds independent of $c$, which implies that $W^s(P^{+,s})|_{\M^+_\eps}(c)$ wraps around the periodic orbit periodic orbit $(w_{p,\eps}(X),q_{p,\eps}(X))$ of (\ref{eq:slow+eps}) (in backwards `time') for all $c$ with the same sign as $c_{\rm h-p}^1$ and $c_{\rm prim}^1$ (cf. Fig. \ref{fig:counthets}). Thus, there must be countably many heteroclinic orbits $\gamma_{h}^{1,k}(\xi)$ -- every $W^s(P^{+,s})|_{\M^+_\eps}(c)$ intersects $\I_{\rm down}$ countably many times -- and the associated speeds $c_{h}^{1,k}$ must all be between $c_{\rm h-p}^1$ and $c_{\rm prim}^1$. Moreover, the decreasing sequence $\{c_{h}^{1,k}\}_{k=1}^\infty$ must have a limit that cannot differ from $c_{\rm h-p}^1$: $c_{h}^{1,k} \downarrow c_{\rm h-p}^1$ as $k \to \infty$.
\hfill $\Box$
\\ \\
By establishing the existence of countably many distinct heteroclinic connections between $P^0$ and $P^{+,s}$, Theorem \ref{T:count1fronts} in a sense considers (one of) the most complex case(s), which is quite far removed from situations in which there are no such connections. To obtain insight in the bifurcations that occur `in between', we can again freeze the reduced slow flow and vary $\I_{\rm down} = \I_{\rm down}(\Phi)$ by increasing $\Phi$ (from $a(1+4a)\chi$ to $\infty$ (\ref{defchi}), (\ref{PsiasPhi}), (\ref{PsiasPhi})). We consider the most simple case and assume that the homoclinic orbit of the frozen flow lies entirely in the $w$-region $(4/(1+4a),9/(2+9a))$ -- so that all $c_{h}^{j,k}$'s of Theorem \ref{T:count1fronts} are positive -- and that $\I_{\rm down}(\Phi) \cap W^s(P^{+,s})|_{\M^+_0} = \emptyset$ at $\Phi = a(1+4a)\chi$ (this can easily be achieved by the unfolded saddle-node approach). As $\Phi$ increases, $\I_{\rm down}(\Phi)$ becomes steeper and the intersection $\I_{\rm down}(\Phi) \cap \{ q =0 \}$ moves over the entire interval determined by (\ref{condw0}), i.e. from $4/(1+4a)$ to $1/a$. Thus, $\I_{\rm down}(\Phi)$ moves through the homoclinic loop spanned by $W^s(P^{+,s})|_{\M^+_0}$ and through the enclosed persistent periodic orbit established by Theorem \ref{T:slowflow}. We tune the parameters such that during the passage of the latter, (\ref{DeltaHpos}) holds and $(a,\Psi,\Phi,\Omega,\Theta) \in \mathcal{S}_{\rm h-p}$, i.e. that Theorem \ref{T:count1fronts} can be applied -- which is also possible. It should be noticed that although the reduced flow (\ref{eq:SRSplus}) is frozen, this is not the case for the perturbed flow (\ref{eq:slow+eps}), since $\rho_1(w)$ varies with $\Theta$ (\ref{rho1}) and $\Theta = \Theta(\Phi)$ (\ref{PsiasPhi}), thus the persistent periodic orbit established in Theorem \ref{T:slowflow} is not frozen, but also varies with $\Phi$ -- this is represented in the sketches of Fig. \ref{fig:spiralaroundperiodic-3config} by the decreasing size of the limiting periodic orbit on $\M^+_\eps$.
\\ \\
Fig. \ref{fig:spiralaroundperiodic-3config} exhibits sketches of 3 configurations of $\I_{\rm down}(\Phi)$ and $W^s(P^{+,s})|_{\M^+_0}$ for increasing $\Phi$, the associated bifurcation scenario is sketched in Fig. \ref{fig:cPhi2}. In Fig. \ref{fig:spiralaroundperiodic-3config}a, $\Phi$ has already passed through the first bifurcation value $\Phi_{\rm prim}$ at which the first 2 primary heteroclinic orbits $\gamma_{\rm prim}^j(\xi)$, $j=1,2$ of Theorem \ref{T:primaryfronts} are created, and through a second one, $\Phi^b_{SN,1}$ -- $\mathcal{O}(\eps)$ close $\Phi_{\rm prim}$ -- at which the first 2 secondary orbits appear. This bifurcation is followed by countably subsequent saddle-node bifurcations until $\Phi$ reaches $\Phi^b_{\rm per}$ at which the 2 limiting heteroclinic orbits between $P^0$ and the persistent periodic solution of Theorem \ref{T:frontstoperiodic} appear and we enter into the realm of Theorem \ref{T:count1fronts}. These orbits next disappear at $\Phi^e_{\rm per}$, Fig. \ref{fig:spiralaroundperiodic-3config}(b) is similar to Fig. \ref{fig:counthets} and represents the 2 countable families of heteroclinic orbits that exist for $\Phi \in (\Phi^b_{\rm per},\Phi^e_{\rm per})$ (Theorem \ref{T:count1fronts}). All these orbits step-by-step disappear in pairs as $\Phi$ is increased further: Fig. \ref{fig:spiralaroundperiodic-3config}c shows the situation with only 5 left -- 4 of these will disappear just before $\Phi$ reaches $\Phi^{+,s}$ at which $\I_{\rm down}(\Phi)$ passes through $P^{+,s}$.
\\ \\
We refrain from giving all rigorous details on which the above sketched scenario is based -- this is in essence a matter of following the lines set out in the proofs of the preceding results. Moreover, we also refrain from working out all possible alternative bifurcation scenarios that may occur -- there are many (sub)cases to consider, some more simple, others more complex than that of Fig. \ref{fig:cPhi2}. Nevertheless, we do briefly come back to this in the upcoming section -- where we consider stationary, sign-changing, case $c=0$ case.
\begin{remark}
\label{rem:moreonslowperiodics}
As in Remark \ref{rem:otherslowflowperiodics}, we note that a result like Theorem \ref{T:frontstoperiodic} on the existence of heteroclinic connections between $P^0$ and a periodic orbit on $\M^+_\eps$ can also be established under the assumption that there is only one critical point $P^{+,c}$ of center type on $\M^+_0$. Similar remarks can be made about the upcoming Theorems \ref{T:statfrontstoperiodic}, \ref{T:higherorder2fronts} and Corollary \ref{cor:higherorderperiodics}. We note -- also as in Remark \ref{rem:otherslowflowperiodics} -- that the analysis of these additional cases is essentially the same as already presented.
\end{remark}

\subsection{Stationary 1-front patterns}
\label{sec:stand1fronts}

\begin{figure}[t]
	\centering
	\includegraphics[width=0.8\linewidth]{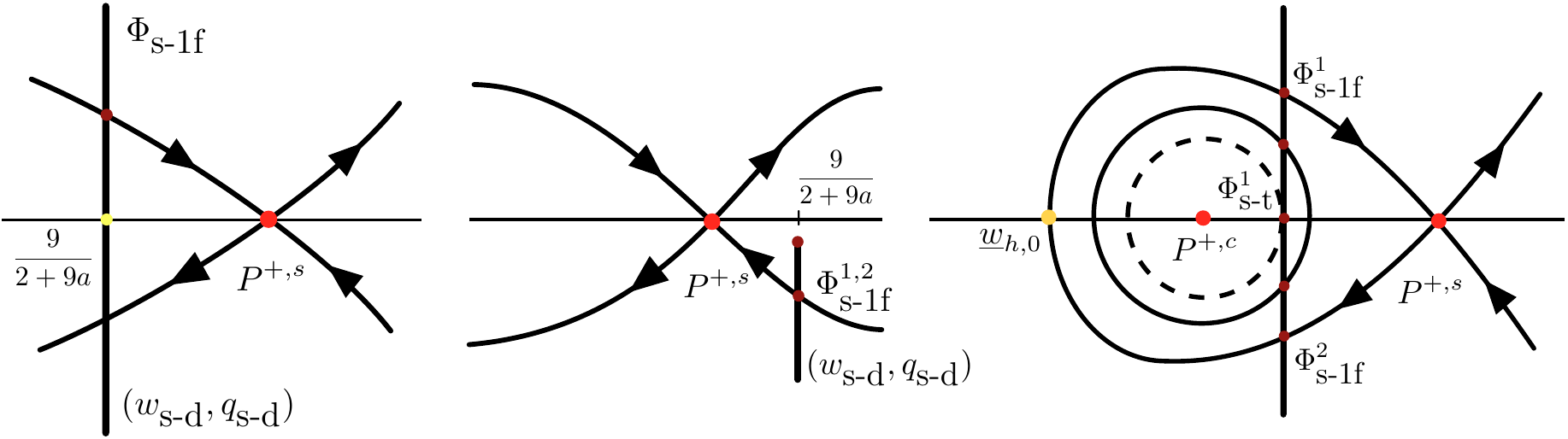}
	\caption{Sketches of 3 configurations of the line $\J_{\rm s-d}$ (\ref{defJs-d}) and $W^s(P^{+,s})|_{\M^+_0}$ as consider in Theorem \ref{T:stand1fronts}: (a) $P^{+,s}$ is the only critical point on $\M^+_\eps$ and $\chi > 0$; (b) $P^{+,s}$ is the only critical point on $\M^+_\eps$ and $\chi < 0$; (c) a center and a saddle on $\M^+_\eps$ with $\underline{w}_{h,0} < 9/(2 + 9a) < w^{+,s}$.}
	\label{fig:3configurationsline-}
\end{figure}

In this section, we construct stationary heteroclinic 1-front patterns that are similar to those constructed in Theorems \ref{T:primaryfronts} and \ref{T:frontstoperiodic}. We immediately note that if the reduced flow on $\M^+_0$ has an unperturbed homoclinic loop $W^s(P^{+,s})|_{\M^+_0} \cap W^s(P^{+,s})|_{\M^+_0}$ -- as in Fig. \ref{fig:slowM+0} -- that it persists as homoclinic solution of (\ref{eq:fast}) on $\M^+_\eps$ for $\eps \neq 0$ -- since (\ref{eq:fast}) with $c=0$ is a reversible system (see also Theorem \ref{T:slowflow}). Thus, we a priori deduce that there cannot be any further non-primary heteroclinic 1-front connections between $P^0 \in \M^0_\eps$ and $P^{+,s} \in \M^+_\eps$ as those of Theorem \ref{T:count1fronts} for $c \neq 0$ (see however also Remark \ref{rem:counthetsc=0} for a result similar to Theorem \ref{T:count1fronts}). In the subsequent sections, we will proceed to construct homoclinic and periodic multi-front patterns -- i.e. solutions of (\ref{eq:fast}) that jump up and down between $\M^0_\eps$ and $\M^+_\eps$ -- and show that there is a richness in these kinds of patterns similar to that of Theorem \ref{T:count1fronts}.
\\ \\
As in the previous sections, we approach the bifurcation analysis by freezing the flow on $\M^+_\eps$. Thus, we choose $\Psi, \Omega, \Theta$ as in (\ref{PsiasPhi}) and vary $\Phi$. We know by Lemma \ref{lem:Idown} and (\ref{IdownPhi}) that for $c=0$, the touch down point of $W^u(P^0) \cap W^s(\M^+_\eps)$ is represented by a vertical line/half line $\J_{s-d}$ in the $(w,q)$-plane
\beq
\J_{\rm s-d}= \left\{J_{\rm s-d}(\Phi) = (w_{\rm s-d}(\Phi),q_{\rm s-d}(\Phi)) = \left(\frac{9}{2 + 9a}, \sqrt{\Phi} \left(\frac{\chi}{\Phi} - \frac{2}{a(2 + 9a)}\right)\right), \Phi > 0\right\}
\label{defJs-d}
\feq
(at leading order in $\eps$). Clearly, a 1-front connection between $P^0$ and $P^{+,s}$ corresponds to those values of $\Phi$ for which $J_{\rm s-d}(\Phi) \in W^s(P^{+,s})|_{\M^+_0}$.
\begin{theorem}
\label{T:stand1fronts}
Let $c=0$ and let $\eps$ be sufficiently small. Then, there is a co-dimension 1 set $\R_{\rm s-1f}$ in $(a,\Phi,\Psi,\Omega,\Theta)$-space for which (stationary) 1-front heteroclinic orbits $\gamma_{\rm s-1f}(\xi) \subset W^u(P^0) \cap W^s(P^{+,s})$ exists in (\ref{eq:fast}). More precise, let $\Psi, \Omega, \Theta$ as in (\ref{PsiasPhi}), then:
\\
{\bf (A)} If $P^{+,s} = (w^{+,s},0)$ is the only critical point on (the projection of) $\M^+_\eps$, i.e. if $\C^2 - 4 \A \D > 0$, $\D < 0$, $\E > 0$ (\ref{saceM+0}), then,
\\
$\bullet$ if $\chi > 0$, then there is a unique value $\Phi_{\rm s-1f}$ such that there is a 1-front heteroclinic orbit $\gamma_{\rm s-1f}(\xi) \subset W^u(P^0) \cap W^s(P^{+,s})$ in (\ref{eq:fast});
\\
$\bullet$ if $\chi < 0$, $w^{+,s} < w_{\rm s-d} = 9/(2+9a)$ (\ref{defJs-d}) and there are $\Phi$ such that $\mathcal{H}_{0}^+(w_{\rm s-d}(\Phi),q_{\rm s-d}(\Phi)) < \mathcal{H}_{0}^{+,s}$ (\ref{eq:HamiltonianFRS}), (\ref{manifoldsP+s}), then there are 2 values $\Phi_{\rm s-1f}^j$, $j=1,2$ for which 1-front heteroclinic orbits $\gamma_{\rm s-1f}(\xi) \subset W^u(P^0) \cap W^s(P^{+,s})$ exist in (\ref{eq:fast});
\\
$\bullet$ if $\chi < 0$ and either one of the above additional conditions does not hold, then there is no such stationary 1-front orbit.
\\
{\bf (B)} If there are two critical points on $\M^+_\eps$, the center $P^{+,c}$ and saddle $P^{+,s}$, i.e. if $\C^2 - 4 \A \D > 0$, $\C < 0$, $\D > 0$, then,
\\
$\bullet$ if $\underline{w}_{h,0} < 9/(2 + 9a) < w^{+,s}$ -- with $\underline{w}_{h,0}$ as defined in Theorem \ref{T:slowflow} -- then there are 2 values $\Phi_{\rm s-1f}^1 < \Phi_{\rm s-1f}^2$ for which 1-front heteroclinic orbits $\gamma_{\rm s-1f}(\xi) \subset W^u(P^0) \cap W^s(P^{+,s})$ exist in (\ref{eq:fast});
\\
$\bullet$ if $9/(2 + 9a) >  w^{+,s}$ there is a unique value $\Phi_{\rm s-1f}$ such that a 1-front heteroclinic orbit $\gamma_{\rm s-1f}(\xi) \subset W^u(P^0) \cap W^s(P^{+,s})$ exists in (\ref{eq:fast});
\\
Every heteroclinic orbit $\gamma_{\rm s-1f}(\xi) = (w_{\rm s-1f}(\xi),p_{\rm s-1f}(\xi),b_{\rm s-1f}(\xi),q_{\rm s-1f}(\xi))$ corresponds to a stationary 1-front pattern $(B(x,t), W(x,t)) = (b_{\rm s-1f}(x), w_{\rm s-1f}(x))$ in PDE (\ref{eq:scaled}) that connects the bare soil state $(0,\Psi/\Phi)$ to the uniform vegetation state $(b_+(w^{+,s}),w^{+,s})$.
\end{theorem}
We refer to Fig. \ref{fig:numerics_1_front_standing} for an example of a numerical simulation of (\ref{eq:scaled}) exhibiting a stationary 1-front pattern. Moreover, we notice that -- by symmetry (\ref{eq:symmc}) of (\ref{eq:fast}) with $c=0$ -- the heteroclinic orbit $\gamma_{\rm s-1f}(\xi) \subset W^u(P^0) \cap W^s(P^{+,s})$ has a counterpart $\subset W^u(P^{+,s}) \cap W^s(P^0)$, i.e. an orbit from $P^{+,s}$ to $P^0$. Together, these obits form a heteroclinic cycle between the saddles $P^{+,s}$ to $P^0$.
\\ \\
{\bf Proof.} The result follows directly by studying the possible intersections of $W^s(P^{+,s})|_{\M^+_\eps}$ and the (vertical) line $\J_{\rm s-d}$ in combination with the observation that the range of $q_{\rm s-d}(\Phi)$ is $\mathbb{R}$ for $\chi > 0$, while it's bounded from above by a negative number for $\chi < 0$ (\ref{defJs-d}). See Fig. \ref{fig:3configurationsline-}.
\\ \\
Since all periodic orbits on $\M^+_\eps$ persist for $c = 0$ (Theorem \ref{T:slowflow}), we also `automatically' obtain a result similar to Theorem \ref{T:frontstoperiodic} on the existence of heteroclinic connections $\gamma_{\rm s-p}^j(\xi)= (b_{\rm s-p}^j(\xi),p_{\rm s-p}^j(\xi),w_{\rm s-p}^j(\xi),q_{\rm s-p}^j(\xi))$ of (\ref{eq:fast}) between the critical point $P^0 \in \M^0_\eps$ and one of the periodic orbits $\gamma_{p,\eps}(X)=(b_{p,\eps}(X),p_{p,\eps}(X),w_{p,\eps}(X),q_{p,\eps}(X)) \subset \M^+_\eps$ determined by $\mathcal{H}_{p,0}^+ \in (\mathcal{H}_{0}^{+,c}, \mathcal{H}_{0}^{+,s})$ (\ref{eq:HamiltonianSRS}) -- note that this orbit is $\mathcal{O}(\eps^2)$ close to the level set $\mathcal{H}_{0}^+(w,q) = \mathcal{H}_{p,0}^+$ (cf. (\ref{eq:slow+eps}) with $c=0$).
\begin{theorem}
\label{T:statfrontstoperiodic}
Let $c=0$, $\Psi, \Omega, \Theta$ as in (\ref{PsiasPhi}) such that (\ref{asshomM+eps}) holds, let $(b_{p,\eps}(X),p_{p,\eps}(X),w_{p,\eps}(X),q_{p,\eps}(X)) \subset \M^+_\eps$ be a periodic solution of (\ref{eq:fast}) determined by $\mathcal{H} \in (\mathcal{H}_{0}^{+,c}, \mathcal{H}_{0}^{+,s})$ and let $\eps$ be sufficiently small. Assume that $\underline{w}_{h,0} < 9/(2 + 9a) < w^{+,s}$ and define
\begin{equation}
\label{defPhist}
\Phi_{\rm s-t} = \frac12 a (2 + 9a) \chi \in (\Phi_{\rm s-1f}^1,\Phi_{\rm s-1f}^2), \;
\mathcal{H}_{\rm s-t}^+ = \mathcal{H}_{0}^+(w_{\rm s-d}(\Phi_{\rm s-t}),q_{\rm s-d}(\Phi_{\rm s-t})) = \mathcal{H}_{0}^+(\frac{9}{2+9a},0)
\end{equation}
with $\Phi_{\rm s-1f}^{1,2}$ as defined in Theorem \ref{T:stand1fronts}, $(w_{\rm s-d}(\Phi),q_{\rm s-d}(\Phi))$ as in (\ref{defJs-d}) and $\mathcal{H}_{\rm s-t}^+ > \mathcal{H}_{0}^{+,c}$ (unless $P^{+,c} = (\frac{9}{2+9a},0)$ (restricted to $\M^+_0$) -- see Fig. \ref{fig:3configurationsline-}c). For all $\mathcal{H} \in (\mathcal{H}_{\rm s-t}^+, \mathcal{H}_{0}^{+,s})$, there are 2 values $\Phi_{\rm p-1f}^{1,2} = \Phi_{\rm p-1f}^{1,2}(\mathcal{H}) $ -- with $\Phi_{\rm s-1f}^{1} < \Phi_{\rm p-1f}^{1} < \Phi_{\rm s-t}< \Phi_{\rm p-1f}^{2} < \Phi_{\rm s-1f}^{2}$ -- that determine 2 distinct heteroclinic orbits $\gamma_{\rm s-p}^j(\xi; \mathcal{H})$ of (\ref{eq:fast}) between the critical point $P^0 \in \M^0_\eps$ and the periodic orbit $(b_{p,\eps}(X),p_{p,\eps}(X),w_{p,\eps}(X),q_{p,\eps}(X)) \subset \M^+_\eps$ (and $\gamma_{\rm s-p}^j(\xi; \mathcal{H})= (b_{\rm s-p}^j(\xi; \mathcal{H}),p_{\rm s-p}^j(\xi; \mathcal{H}),w_{\rm s-p}^j(\xi; \mathcal{H}),q_{\rm s-p}^j(\xi; \mathcal{H}))$). The orbits $\gamma_{\rm s-p}^j(\xi; \mathcal{H})$ ($j=1,2$) correspond to stationary 1-front patterns $(B(x), W(x)) = (b_{\rm s-p}^j(x; \mathcal{H}),w_{\rm s-p}^j(x; \mathcal{H}))$ in PDE (\ref{eq:scaled}) that connect the bare soil state $(0,\Psi/\Phi)$ to the spatially periodic pattern $(b_{p,\eps}(X),w_{p,\eps}(X))$.
\end{theorem}
\begin{remark}
\label{rem:counthetsc=0}
Together, Theorems \ref{T:stand1fronts} and \ref{T:statfrontstoperiodic} provide the possibility to establish a result similar to that of Theorem \ref{T:count1fronts} in the case that $c_{\rm h-p}^j$ and $c_{\rm prim}^j$ do not have the same signs. Assume we have -- for a certain parameter combination $(a,\Psi,\Phi,\Omega,\Theta)$ such that (\ref{asshomM+eps}) holds -- that $c_{\rm h-p}^1 < 0 < c_{\rm prim}^1$. This implies that the point $J_{\rm s-d}(\Phi)$ on $\mathcal{T}_{\rm down}$ -- Lemma \ref{lem:Idown} -- must lie between the intersections of $\mathcal{T}_{\rm down}$ with the unperturbed homoclinic orbit that determines $c_{\rm prim}^1 > 0$ (Theorem \ref{T:stand1fronts}) and the persisting periodic orbit that determines $c_{\rm h-p}^1 < 0$ (Theorem \ref{T:statfrontstoperiodic}). Thus, $J_{\rm s-d}(\Phi)$ determines a level set $\mathcal{H}_{0}^+(w,q) = \mathcal{H}_{0}^+(w_{\rm s-d}(\Phi),q_{\rm s-d}(\Phi)) = \bar{\mathcal{H}} \in (\mathcal{H}_{0}^{+,c}, \mathcal{H}_{0}^{+,c})$ and we know by Theorem \ref{T:statfrontstoperiodic} that $\Phi = \Phi^2_{\rm p-2f}(\bar{\mathcal{H}})$: $(w_{\rm s-d}(\Phi),q_{\rm s-d}(\Phi))$ is the touchdown point of a heteroclinic orbit between the bare soil state and the persisting (stationary) periodic orbit determined by the level set $\mathcal{H}_{0}^+(w,q) = \bar{\mathcal{H}}$. It then follows by arguments similar to those in the proof of Theorem \ref{T:count1fronts} that there are countably many $c$-values $0 < ... < c_{h}^{1,k} < ... < c_{h}^{1,1} < c_{h}^{1,0} = c_{\rm prim}^1$ with $c_{h}^{1,k} \downarrow 0$ for $k \to \infty$ for which non-primary heteroclinic connections between $P^0$ and $P^{+,s}$ exist -- as in Theorem \ref{T:count1fronts}. The main difference with Theorem \ref{T:count1fronts} is that $c = 0$ determines a stationary orbit and not an attracting one: for $c$ slightly above $c = 0$, the unstable manifold $W^s(P^{+,s})|_{\M^+_\eps}$ only spirals inwards very weakly (in backward time). As a consequence, the number of intersections $W^s(P^{+,s})|_{\M^+_\eps} \cap \mathcal{T}_{\rm down}$ with $\mathcal{H}_{0}^+(w,q) > \bar{\mathcal{H}}$ increases (without bound) as $c \downarrow 0$.
\end{remark}

\subsection{Stationary homoclinic 2-front patterns: vegetation spots and gaps}
\label{sec:stand2fronts}

\begin{figure}[t]
	\includegraphics[width= \textwidth]{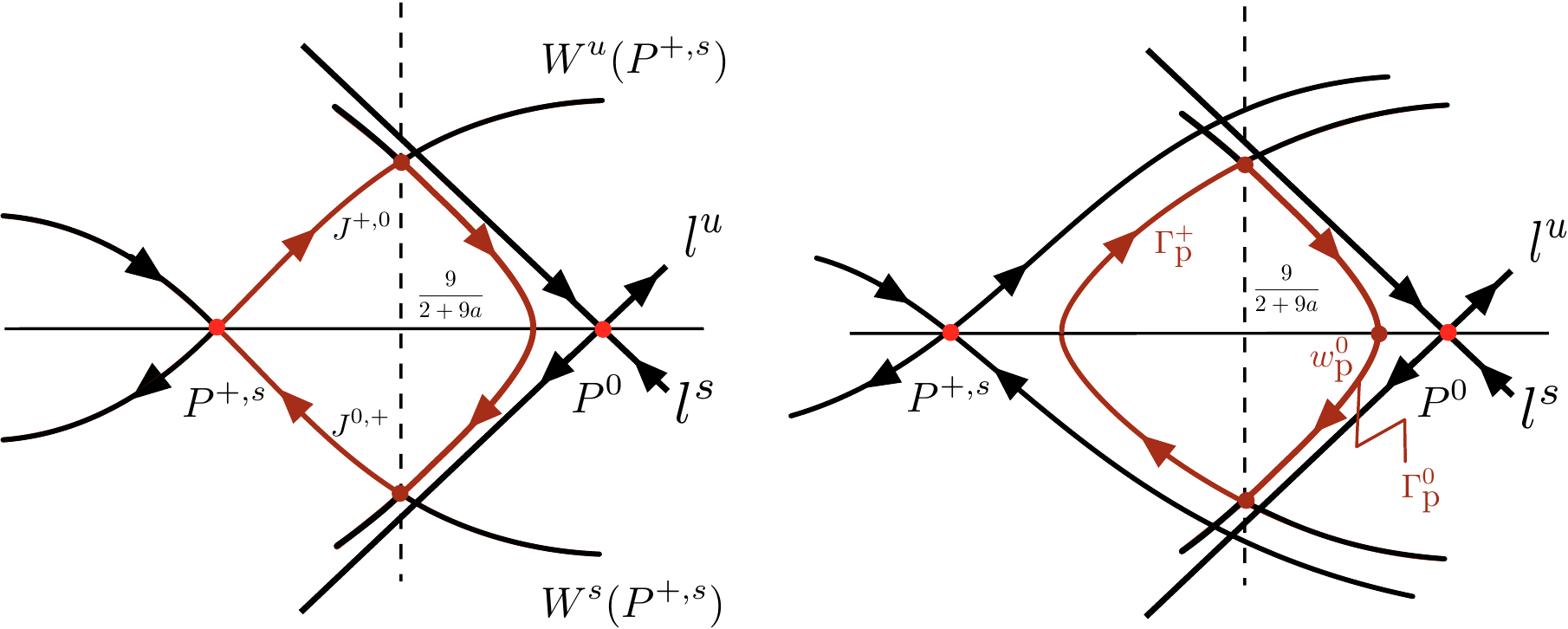}
	\caption{Sketches of 2 skeleton structures of multi-front patterns: (a) The stationary homoclinic 2-front gap pattern of Theorem \ref{T:standing2fgaps}. (b) The spatially periodic multi-front spot/gap pattern of Theorem \ref{T:periodics} as (projected) closed orbit.}
	\label{fig:gapsandperiodics}
\end{figure}

In this section we construct stationary 2-front patterns that correspond to vegetation spots or vegetation gaps -- the latter sometimes also interpreted as fairy circles. These patterns are observed in nature and appear as stable patterns in simulations of (\ref{eq:meron})/(\ref{eq:scaled}) -- see \cite{Zelnik2015} and Figs. \ref{fig:numerics_2_front_standing} and \ref{fig:numerics_fairy_circle}. The patterns/orbits to be constructed are symmetric with respect to the reversibility symmetry in (\ref{eq:scaled}) that persists as (\ref{eq:symmc}) into (\ref{eq:fast}) -- with $c = 0$. As a consequence, we may expect that these patterns are generic, in the sense that they exist in open regions within parameter space -- see for instance \cite{Dreview}. Notice that this is unlike the stationary -- but non-symmetric -- 1-front patterns of the previous section, that only exist for $(a, \Phi, \Psi, \Omega, \Theta) \in \R_{\rm s-1f}$, an explicitly determined co-dimension 1 manifold (Theorem \ref{T:stand1fronts}).
\\
\begin{figure}[t]
	\centering
	\includegraphics[width=0.8\textwidth]{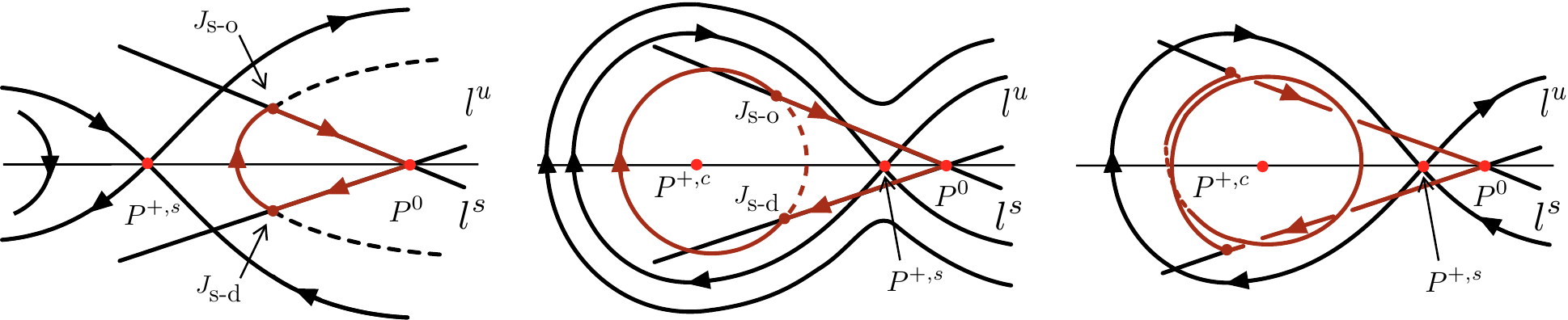}
	\caption{Sketches of 3 (projected) `skeleton structures' of stationary 2-fronts homoclinic to $P^0$: (a) (b) Two examples of the 2 skeleton structures $\Gamma_{\mbox{s-2f}}$ as considered in Theorem \ref{T:standing2fspots}; (c) The extended skeleton $\Gamma_{\rm s-2f}^{\rm ext}$ of Theorem \ref{T:higherorder2fronts} and an associated higher order 2-front homoclinic with an additional full extra `spatial oscillation' during its passage along $\mathcal{M}^+_{\varepsilon}$.}
	\label{fig:3skeletons}
\end{figure}
\\
We first consider the (localized) spots: localized vegetated regions embedded within bare soil. Thus, these spots correspond to solutions of (\ref{eq:fast}) that are homoclinic to the bare soil state $P^0$. Singularly perturbed models of the type (\ref{eq:fast}), can have homoclinic (pulse) solutions of various types. The localized vegetation (spot) patterns constructed in \cite{SewaltD17,BastiaansenD19} in the context of the extended Klausmeier model -- also called generalized Klausmeier-Gray-Scott model \cite{vdSteltetal13} -- make a fast excursion away from the slow manifold that contains the critical point associated to the bare soil state following a homoclinic solution of the fast reduced system. As a consequence these spots are `narrow', their size scales with $\eps$. Although such pulses are also exhibited by the present model -- see Fig. \ref{fig:numerics_fairy_circle_and_1_pulse}b and Remark \ref{rem:2ftofasthom} -- we focus here on 2-front patterns, i.e. orbits homoclinic to $P^0$ that `jump' from $\M^0_\eps$ to $\M^+_\eps$, follow the slow flow on $\M^+_\eps$ over an $\OO(1)$ distance and jump back again -- by its second fast reduced heteroclinic front -- to $\M^0_\eps$. Since the extended Klausmeier model only has one slow manifold \cite{SewaltD17,BastiaansenD19}, such orbits cannot exist in that model. Moreover, these patterns have a well-defined width that does not decrease to $0$ as $\eps \downarrow 0$, a property that appears to be natural in observed ecosystems \cite{Zelnik2015}.
\\ \\
The construction of the most simple -- primary, cf. section \ref{sec:primtrav1front} -- singular `skeleton structure' $\Gamma_{\rm s-2f} \subset \mathbb{R}^4$ -- the phase space of (\ref{eq:fast}) -- of a stationary homoclinic 2-front orbit $\gamma_{\rm s-2f}(\xi)$ to $P^0$ of (\ref{eq:fast}) is relatively straightforward (but somewhat involved/technical). Since the homoclinic orbit $\gamma_{\rm s-2f}(\xi) \subset W^u(P^0)$, it follows $\ell^u$ along $\M^0_\eps$, takes off from $\M^0_\eps$ by following the fast reduced flow and touches down on $\M^+_\eps$ near an element of $\I_{\rm down}$ -- see section \ref{sec:touchdown}. In fact, since we consider stationary spots, the touch down point is near a point $J_{\rm s-d} \in \J_{\rm s-d}$ (\ref{defJs-d}) -- where we note that $\gamma_{\rm s-2f}(\xi)$ cannot exactly touch down on  $\I_{\rm down}$/$\J_{\rm s-d}$ since $\gamma_{\rm s-2f}(\xi) \notin \W^s(\M^+_\eps)$ -- see the proof of (upcoming) Theorem \ref{T:standing2fspots}. The touch down point $J_{\rm s-d}$ determines (at leading order) a level set $\HH_{\rm s-d}$ of the Hamiltonian $\mathcal{H}_{0}^+(w,q)$ (\ref{eq:HamiltonianSRS}) of the slow reduced flow (\ref{eq:SRSplus}) on $\M^+_\eps$: $\HH_{\rm s-d} = \HH_{\rm s-d}(\Phi) = \mathcal{H}_{0}^+(w_{\rm s-d}(\Phi),q_{\rm s-d}(\Phi))$. As long as it remains (exponentially) close to $\M^+_\eps$, the homoclinic orbit-to-be-constructed remains asymptotically close to the level set $\mathcal{H}_{0}^+(w,q) = \HH_{\rm s-d}$. This construction provides the first half of skeleton $\Gamma_{\rm s-2f}$, the second part follows by the symmetry (\ref{eq:symmc}) -- with $c=0$. Completely analogous to $\I_{\rm down}$, one can define $\I_{\rm off}$ as the points on $\M^+_\eps$ that determine the evolution of orbits in $W^s(P^0) \cap W^u(\M^+_\eps)$ after their jump from $\ell^s \subset \M^0_\eps$ through the fast field {\it in backwards time}. In fact, it follows by the symmetry (\ref{eq:symmc}) that $\I_{\rm off}$ and its stationary counterpart $\J_{\rm s-o}$ correspond exactly to the reflections of $\I_{\rm down}$ and $\J_{\rm s-d}$ with respect to the $q$-axis -- where we consider $\I_{\rm down}$/$\J_{\rm s-d}$ and $\I_{\rm down}$/$\J_{\rm s-d}$ within the (projected) 2-dimensional representation of $\M^+_\eps$ as in Lemma \ref{lem:Idown} and in (\ref{defJs-d}). Thus, for a given $\Phi$, the take off point $J_{\rm s-d}$ is given by $J_{\rm s-o}(\Phi) = (w_{\rm s-d}(\Phi),-q_{\rm s-d}(\Phi))$; this point also lies on the level set $\HH_{\rm s-d}$ -- since $\mathcal{H}_{0}^+(w,q)$ (of course) also is symmetric in $q \to -q$.
\\ \\
We define the region $\SSS_{s-2f}$ in $(a,\Psi,\Phi,\Omega,\Theta)$-space for which the point $J_{\rm s-d}$ (and thus $J_{\rm s-o}$) can be constructed (as above) {\it and} there is a solution of the slow reduced flow (\ref{eq:SRSplus}) on $\M^+_\eps$ that connects $J_{\rm s-d}$ to $J_{\rm s-o}$ -- so that $\Gamma_{\rm s-2f}$ indeed exists as closed singular `loop'. Obviously, $\SSS_{s-2f} \neq \emptyset$ -- see also Fig. \ref{fig:3skeletons} -- however, the fact that both $J_{\rm s-d}, J_{\rm s-o} \in \{\mathcal{H}_{0}^+(w,q) = \HH_{\rm s-d} \}$ does not necessarily imply that $(a,\Psi,\Phi,\Omega,\Theta) \in \SSS_{s-2f}$ for all values for which these points exist on $\M^+_\eps$. For instance, in the case that the saddle $P^{+,s}$ is the only critical point on $\M^+_\eps$, $J_{\rm s-d}$ and $J_{\rm s-o}$ are not connected by a solution of (\ref{eq:SRSplus}) if $\HH_{\rm s-d} > \mathcal{H}_{0}^{+,s}$ -- the value of $\mathcal{H}_{0}^+(w,q)$ at $P^{+,s}$ -- see Fig. \ref{fig:3skeletons}a. Moreover, if $\HH_{\rm s-d} < \mathcal{H}_{0}^{+,s}$ there is an additional condition on $\Phi$ that is determined by the relative positions of $w^{+,s}$ (the $w$-coordinate of $P^{+,s}$), $9/(2 + 9a)$ (the $w$-coordinate of $J_{\rm s-d}/J_{\rm s-o}$) and $\Psi/\Phi$ (the $w$-coordinate of the bare soil state associated to $P^0$). Here, we refrain from working out the full `bookkeeping' details by which (the boundary of) $\SSS_{s-2f}$ is determined -- see also a further brief discussion following Theorem \ref{T:standing2fspots}. We refer to Fig. \ref{fig:3skeletons}a for a case $w^{+,s} < 9/(2 + 9a) < \Psi/\Phi$ (and implicitly $\chi >0$) for which $J_{\rm s-d}$ can only be connected to $J_{\rm s-o}$ if $\Phi < \Phi_{\rm s-t} = \frac12 a (2 + 9a) \chi$ (\ref{defPhist}) -- since we need that $q_{\rm s-d}(\Phi) < 0$. Notice that the sketch in Fig. \ref{fig:3skeletons}a in principle also covers a (sub)case of the situation with two critical points $P^{+,c}$ and $P^{+,s}$ on $\M^+_\eps$ and that Fig. \ref{fig:3skeletons}b considers the case $9/(2 + 9a) < w^{+,s} < \Psi/\Phi$ for this situation. Clearly, there are no further restrictions on $\Phi$ if $\HH_{\rm s-d} < \mathcal{H}_{0}^{+,s}$ if there are two critical points $P^{+,c}$ and $P^{+,s}$ on $\M^+_\eps$, since the orbits on the level set associated to $\HH_{\rm s-d}(\Phi)$ are periodic, while one again has to impose $\Phi < \Phi_{\rm s-t}$ to have a connection between $J_{\rm s-d}$ and $J_{\rm s-o}$ for level sets outside the homoclinic loop, i.e. for $\HH_{\rm s-d} > \mathcal{H}_{0}^{+,s}$. Finally, we note that the skeleton structure $\Gamma_{\rm s-2f}$ can in principle also be constructed for $(a,\Psi,\Phi,\Omega,\Theta)$ such that there is no critical point on $\M^+_\eps$, or only one critical point that is not a saddle but a center (in the limit $\eps \to 0$).
\\ \\
Summarizing, the (open) region $\SSS_{s-2f}$ is defined such that for parameter combinations $(a,\Psi,\Phi,\Omega,\Theta) \in \SSS_{s-2f}$, a singular skeleton $\Gamma_{\rm s-2f} \subset \mathbb{R}^4$ can be constructed as above. In the limit $\eps \to 0$, $\Gamma_{\rm s-2f}$ is spanned by a piece of $\ell^u \subset \M^0_0$ from $P^0$ up to the ($\eps \to 0$ limit of the) take off point from $\M^0_0$ (that has the same $(w,q)$-coordinates as $J_{s-d}$ in the limit $\eps \to 0$), the jump through the fast field along (a piece of) $W^u(\M^0_0) \cap W^s(\M^+_0)$ (\ref{WsWu0}) towards the $\eps \to 0$ limit of the (projected) touch down point $J_{s-d} \in \M^+_0$, the connection along $\M^+_0$ to $J_{s-d}$ by the slow reduced flow (on the level set $\{\mathcal{H}_{0}^+(w,q) = \HH_{\rm s-d} \}$) up to (the $\eps \to 0$ limit of) the take off point $J_{\rm s-o}$, followed by a fast jump backwards along (a piece of) $W^s(\M^0_0) \cap W^u(\M^+_0)$ to $\M^0_0$ and a final piece of $\ell^s$ (up to $P^0$) -- see Figs. \ref{fig:3skeletons}a and \ref{fig:3skeletons}b for 2 sketches of projections of $\Gamma_{\rm s-2f}$ on $\M^+_0$ that skip both jumps through the fast field. The proof of the persistence of $\Gamma_{\rm s-2f}$ for $\eps \neq 0$ relies heavily on the reversibility symmetry of (\ref{eq:fast}) with $c=0$ (\ref{eq:symmc}).
\begin{theorem}
\label{T:standing2fspots}
Let $(a,\Psi,\Phi,\Omega,\Theta) \in \SSS_{s-2f}$ and $\Gamma_{\rm s-2f} \subset \mathbb{R}^4$ be the singular skeleton constructed above. Then, there is for $\eps > 0$ sufficiently small a symmetric homoclinic 2-front orbit $\gamma_{\rm s-2f}(\xi) = (b_{\rm s-2f}(\xi),p_{\rm s-2f}(\xi),w_{\rm s-2f}(\xi),q_{\rm s-2f}(\xi)) \subset W^u(P^0) \cap W^s(P^0)$ of (\ref{eq:fast}) with $c=0$ that merges with $\Gamma_{\rm s-2f}$ as $\eps \downarrow 0$. The associated stationary pattern $(B(x,t), W(x,t)) = (b_{\rm s-2f}(x),w_{\rm s-2f}(x))$ in (\ref{eq:scaled}) represents a stationary localized vegetation spot embedded in bare soil.
\end{theorem}
\noindent
We refer to Figs. \ref{fig:basicpatterns}b and \ref{fig:numerics_2_front_standing} for numerical observations of these 2-front spot patterns. In Fig. \ref{fig:projection_pulse_2front_manifolds}, the projection of the 2-front orbit on the $(w,b)$-plane is given; it clearly shows the slow-fast-slow-fast-slow nature of the pattern: it first follows $\M^0_\eps$ (slowly), jumps to $\M^+_\eps$, follows the slow flow on $\M^+_\eps$, jumps back to $\M^0_\eps$ and slowly returns to $P^0$.
\\ \\
To get some insight in the boundaries of $\SSS_{s-2f}$ -- and thus is the bifurcations of $\gamma_{\rm s-2f}(\xi)$ -- we can (as usual) `freeze' the flow on $\M^+_0$ by choosing $\Psi, \Omega, \Theta$ as in (\ref{PsiasPhi}) and vary $\Phi$. We need to be aware though that $\ell^{u,s}$ and $P^0$ do vary with $\Phi$ (i.e. they are not frozen). In the situation sketched in Fig. \ref{fig:3skeletons}a -- thus with $w^{+,s} < 9/(2 + 9a) < \Psi/\Phi$ and $\chi >0$ -- we see that the distance between the 2 fronts of $\gamma_{\rm s-2f}(\xi)$ approaches $\infty$ as $\Phi \downarrow \Phi_{\rm s-1f} < \Phi_{\rm s-t}$ as defined in Theorems \ref{T:stand1fronts} and \ref{T:statfrontstoperiodic}: $\gamma_{\rm s-2f}(\xi)$ obtains the character of the superposition of the 1-front heteroclinic orbit $\gamma_{\rm s-1f}(\xi)$ of Theorem \ref{T:stand1fronts} between $P^0$  and $P^{+,s}$ and its symmetrical counterpart -- (\ref{eq:symmc}) -- connecting $P^{+,s}$ back to $P^{0}$. The other boundary corresponds to $\Phi \uparrow \Phi_{\rm s-t}$: the distance traveled along $\M^+_\eps$ decreases to $0$ and $\gamma_{\rm s-2f}(\xi)$ detaches from $\M^+_\eps$. However, since the angle between $\ell^u$ and $\ell^s$ is determined by $\sqrt{\Phi}$ (\ref{deflsu}), this can only happen as also  $\Psi/\Phi \downarrow 9/(2 + 9a)$, i.e. as the `projected triangle' of Fig. \ref{fig:3skeletons}a that represents the skeleton $\Gamma_{\rm s-2f}$ entirely contracts to a point -- see also Remark \ref{rem:2ftofasthom}. The bifurcational structure associated to the situation sketched in Fig. \ref{fig:3skeletons}b is quite different: as $\Phi$ decreases towards $\Phi_{\rm s-1f}^1 < \Phi_{\rm s-t}$, $\gamma_{\rm s-2f}(\xi)$ does not merge with a (superposition of 2) 1-front(s) $\gamma_{\rm s-1f}(\xi)$ of Theorem \ref{T:stand1fronts}. In fact, $\gamma_{\rm s-2f}(\xi)$ does not bifurcate at all, the distance between the 2 fronts of $\gamma_{\rm s-2f}(\xi)$ remains bounded as $\Phi$ passes through $\Phi_{\rm s-1f}^1$ (the main difference between the cases $\Phi > \Phi_{\rm s-1f}^1$ and $\Phi < \Phi_{\rm s-1f}^1$ is the sign of $\HH_{\rm s-d} - \mathcal{H}_{0}^{+,s}$: $\gamma_{\rm s-2f}(\xi)$ follows an orbit of (\ref{eq:SRSplus}) outside the homoclinic loop for $\Phi < \Phi_{\rm s-1f}^1$). Moreover, as $\Phi$ increases towards $\Phi_{\rm s-t}$, the distance $\gamma_{\rm s-1f}(\xi)$ travels along $\M^+_\eps$ also does not go to 0, in fact $\gamma_{\rm s-1f}(\xi)$ (almost) follows the entire periodic orbit of the level set $\{\mathcal{H}_{0}^+(w,q) = \HH_{\rm s-d}(\Phi_{\rm s-t}) = \HH_{\rm s-t}\}$, the critical/limiting orbit of Theorem \ref{T:statfrontstoperiodic}. Again, this can only happen if also the projection of $P^0$ on $\M^+_\eps$ merges with this periodic orbit. We refrain from going any further into the details of these -- and other -- bifurcations of $\gamma_{\rm s-2f}(\xi)$.
\\ \\
{\bf Proof of Theorem \ref{T:standing2fspots}.} The proof follows the geometrical approach developed in \cite{DH96,DGK01}, equivalently the more analytical approach of \cite{JKK96} could be employed. The construction of $\gamma_{\rm s-2f}(\xi)$ is based on the `intermediate' orbit $\gamma_{\rm i-1f}(\xi) \subset W^u(P^0) \cap W^s(\M^{+}_\eps)$, the heteroclinic connection between $P^0$ and $\M^+_\eps$ that touches down on $J_{s-d} \in \M^+_\eps$; $\gamma_{\rm i-1f}(\xi)$ follows the slow flow along $\M^+_\eps$ and is thus asymptotically close to the skeleton $\Gamma_{\rm s-2f}$ up to the take off point $J_{\rm s-o}$ ($\gamma_{\rm i-1f}(\xi) \subset W^s(\M^{+}_\eps)$ and thus cannot take off from $\M^+_\eps$). The homoclinic orbit $\gamma_{\rm s-2f}(\xi)$ is constructed as a symmetric orbit -- i.e. an orbit that passes through the plane $\{p=q=0\}$ at its `midpoint' -- that is exponentially close to $\gamma_{\rm i-1f}(\xi)$ up to the point it takes off from $\M^+_\eps$.
\\ \\
Since $b_-(w) < 1/2 < b_+(w)$ (\ref{b0bpm}), $\gamma_{\rm i-1f}(\xi)$ intersects the hyperplane $\{b = 1/2\}$ transversally in the point $P_{\rm i-1f}$ -- by definition. We define for some sufficiently small $\sigma$ (independent of $\eps$), the (bounded) 1-dimensional curve $C_{\rm i-1f}^\sigma \subset \{b = 1/2\}$ as the (first, transversal) intersection of $W^u(P^0)$ and $\{b = 1/2\}$ that is at a distance of maximal $\sigma$ away from $P_{\rm i-1f}$; in other words, $C_{\rm i-1f}^\sigma =  W^u(P^0) \cap \{b = 1/2\} \cap \{|(b,p,w,q) - P_{\rm i-1f}| \leq \sigma\}$. By choosing $\gamma(0) \in C_{\rm i-1f}^\sigma$, the curve $C_{\rm i-1f}^\sigma$ provides a parametrization of orbits $\gamma(\xi)$ in $W^u(P^0)$ near $\gamma_{\rm i-1f}(\xi)$ . In fact, the saddle structure of the fast flow around $\M^+_\eps$ cuts $W^u(P^0)$, and thus $C_{\rm i-1f}^\sigma$, exactly in two along $\gamma_{\rm i-1f}(\xi) \subset W^s(\M^{+}_\eps)$: orbits $\gamma(\xi) \subset W^u(P^0)$ with (by definition) $\gamma(0) \in C_{\rm i-1f}^{\sigma, {\rm r}}$ cross through the plane $\{p=0\}$ near $\M^+_\eps$ so that their $b$-coordinate changes direction; the $b$-coordinates of $\gamma(\xi)$'s with $\gamma(0)$ in $C_{\rm i-1f}^{\sigma}\setminus C_{\rm i-1f}^{\sigma, {\rm r}}$ do not change direction, these $\gamma(\xi)$'s pass along $\M^+_\eps$ without the possibility of returning to $\M^0_\eps$. Thus, we may uniquely parameterize orbits $\gamma(\xi) \subset W^u(P^0)$ that pass through $\{p=0\}$ near by $\M^+_\eps$ that are $\sigma$-close to $\gamma_{\rm i-1f}(\xi)$ by the distance $d$ between their initial point $\gamma(0) \in C_{\rm i-1f}^{\sigma, {\rm r}}$ and $P_{\rm i-1f}$ ($\in \partial C_{\rm i-1f}^{\sigma, {\rm r}}$) -- where we have implicitly used the fact that $W^u(\M^0_\eps) \supset W^u(P^0)$ is $C^1-\OO(\eps)$ close to its $\eps \to 0$ limit $W^u(\M^0_0)$ \cite{Jones1995,Kaper1999}. We denote these $\gamma(\xi)$'s by $\gamma(\xi; d)$.
\\ \\
Since the flow on $\M^+_\eps$ is $\OO(\eps)$ slow, orbits $\gamma(\xi; d)$ with $\gamma(0) \in C_{\rm i-1f}^{\sigma, {\rm r}}$ can only follow $\M^+_\eps$ over an $\OO(1)$ distance (w.r.t. $\eps$) for $d$ exponentially small; in fact it is necessary that $d = \OO({\rm exp}(-\lambda_{f,+}(9/(2 + 9a))/\eps)$, where $\lambda_{f,+}(w_0) = \sqrt{w_0b_+(w_0) + 2aw_0 - 2}$, the unstable eigenvalue of the reduced fast flow (\ref{eq:fastblue}) -- with $c=0$ -- associated to the critical point $(b_+(w_0),0)$, and $w_0 = 9/(2+9a)$ is the (leading order) $w$-coordinate of $J_{s-d}$ (\ref{defJs-d}). As $d$ decreases, the `time' (i.e. distance) $\gamma(\xi; d)$ remains exponentially close to $\M^+_\eps$ increases monotonically (as follows from a direct perturbation analysis -- for instance along the lines of \cite{CarterSandstede15}). Equivalently, the distance between $J_{s-d}$ and the (projected) point on $\M^+_\eps$ at which $\gamma(\xi; d)$ crosses through the $\{p=0\}$-plane also increases monotonically with decreasing $d$ -- where we note that this point marks the transition between $\gamma(\xi; d)$ approaching the `fast saddle' $\M^+_\eps$ exponentially close to $W^s(\M^+_\eps)$ and moving away from $\M^+_\eps$ exponentially close to $W^u(\M^+_\eps)$. Clearly, for $d$ `too large' $\gamma(\xi; d)$ passes through $\{p=0\}$ before the slow flow on $\M^+_\eps$ -- and thus $\gamma(\xi)$ itself (since it is exponentially close to $\M^+_\eps$) -- passes through $\{q=0\}$ -- Fig. \ref{fig:3skeletons}. However, by decreasing $d$, we can delay the passage of $\gamma(\xi; d)$ through $\{p=0\}$ until after it passes through $\{q=0\}$. It follows that there must be a $d_*$ such that the associated orbit $\gamma(\xi; d_*)$ passes through $\{p=0\}$ and $\{q=0\}$ simultaneously: $\gamma(\xi; d_*) \cap \{p=q=0\} \neq \emptyset$.
\\ \\
The orbit $\gamma_{\rm s-2f}(\xi)$ coincides with $\gamma(\xi; d_*)$: by the reversibility symmetry (\ref{eq:symmc}) -- with $c=0$ -- it is symmetric around the `midpoint' at which it passes through $\{p=q=0\}$ , so that indeed $\gamma_{\rm s-2f}(\xi) \subset W^u(P^0) \cap W^s(P^0)$, a homoclinic 2-front orbit to $P^0$ (asymptotically close to the skeleton $\Gamma_{\rm s-2f}$). \hfill $\Box$
\\ \\
The orbit $\gamma_{\rm s-2f}(\xi)$ is only the first of a countable family of homoclinic 2-front orbits if the slow piece of the skeleton $\Gamma_{\rm s-2f}$ is part of a closed orbit, i.e. if the connected part of the level set $\{\mathcal{H}_{0}^+(w,q) = \HH_{\rm s-d} \}$ that contains the ($\eps \to 0$ limits of the) touch down and take off points $J_{\rm s-d}$ and $J_{\rm s-o}$ is a closed orbit. In this case, the intermediate heteroclinic 1-front $\gamma_{\rm i-1f}(\xi) \subset W^u(P^0) \cap W^s(\M^{+}_\eps)$ introduced in the proof of Theorem \ref{T:standing2fspots} coincides with one of the 2 heteroclinic orbits $\gamma_{\rm s-p}^j(\xi; \mathcal{H}_{\rm s-d})$ between the critical point $P^0 \in \M^0_\eps$ and the periodic orbit $\gamma_{p,\eps}(X) \subset \M^+_\eps$ established in Theorem \ref{T:statfrontstoperiodic}. Thus, in this case $\gamma_{\rm i-1f}(\xi)$ passes countably many times through the plane $\{q=0\}$. By steadily decreasing $d$, we can now determine a sequence of critical values $d^i_*$, $i=0, 1, 2, ...$ such that the associated orbits $\gamma(\xi; d^i_*)$ pass through $\{p=q=0\}$ after $i$ preceding passages through $\{q=0\}$. Thus, the {\it primary} orbits $\gamma_{\rm s-2f}(\xi)$ of Theorem \ref{T:standing2fspots} correspond to $\gamma(\xi; d^0_*)$ and the `higher order' (stationary) homoclinic 2-front orbits $\gamma_{\rm s-2f}^i(\xi) \subset W^u(P^0) \cap W^s(P^0)$ coincide with the orbits $\gamma(\xi; d^i_*)$ for $i \geq 1$. By the symmetry (\ref{eq:symmc}) -- with $c=0$ -- these orbits are also symmetric around the `midpoint' at which they pass through $\{p=q=0\}$: $\gamma_{\rm s-2f}^i(\xi)$ traces $i$ full circuits over the closed orbit determined by the level set $\{\mathcal{H}_{0}^+(w,q) = \HH_{\rm s-d} \}$ during its passage along $\M^+_\eps$ -- see the sketches in Fig. \ref{fig:sketchessection1}(c,d).
\begin{theorem}
\label{T:higherorder2fronts}
Let $(a,\Psi,\Phi,\Omega,\Theta) \in \SSS_{s-2f}$ such that (\ref{asshomM+eps}) holds, and let $\Gamma_{\rm s-2f}^{\rm ext} \subset \mathbb{R}^4$ be the extension of the singular skeleton $\Gamma_{\rm s-2f}$ of Theorem \ref{T:standing2fspots} that includes the entire closed orbit on $\M^+_\eps$ determined by the level set $\{\mathcal{H}_{0}^+(w,q) = \HH_{\rm s-d} \}$. Then, for $\eps > 0$ sufficiently small, there is a countable family of symmetric 2-front orbits $\gamma_{\rm s-2f}^i(\xi) = (b_{\rm s-2f}^i(\xi),p_{\rm s-2f}^i(\xi),w_{\rm s-2f}^i(\xi),q_{\rm s-2f}^i(\xi)) \subset W^u(P^0) \cap W^s(P^0)$ of (\ref{eq:fast}) with $c=0$ ($i=0,1,2,...$) that merges with $\Gamma_{\rm s-2f}^{\rm ext}$ as $\eps \downarrow 0$ for all $i \geq 1$ ($\gamma_{\rm s-2f}^0(\xi) = \gamma_{\rm s-2f}(\xi)$ of Theorem \ref{T:standing2fspots} merges with $\Gamma_{\rm s-2f}$ as $\eps \downarrow 0$). The associated stationary patterns $(B(x,t), W(x,t)) = (b_{\rm s-2f}^i(x),w_{\rm s-2f}^i(x))$ in (\ref{eq:scaled}) represent stationary localized vegetation spots embedded in bare soil with an increasing number of spatial oscillations in the vegetated area.
\end{theorem}
\noindent
The construction of stationary homoclinic 2-front gap patterns -- localized bare soil areas surrounded by vegetation -- goes along exactly the same lines as the above construction of localized spot patterns. The main difference is that the homoclinic orbits-to-be-constructed are $\subset W^u(P^{+,s}) \cap W^s(P^{+,s})$ so that the structure of orbits taking off and touching down now has to start out from the saddle $P^{+,s} \in \M^+_\eps$. Nevertheless, the construction of the skeleton structure $\Gamma_{\rm g-2f}$ is completely similar to that of $\Gamma_{\rm s-2f}$. Therefore, we only provide the essence of the construction of $\Gamma_{\rm g-2f}$.
\\ \\
First, we need to assume that there is a critical point $P^{+,s} \in \M^+_0$ of saddle type. The skeleton structure $\Gamma_{\rm g-2f}$ consists of a piece of $W^u(P^{+,s}) \subset \M^+_0$ from $P^{+,s}$ up to the ($\eps \to 0$ limit of the) take off point $J_{g-o}^{+, 0}$ from $\M^+_0$, followed by (a piece of) $W^u(\M^+_0) \cap W^s(\M^0_0)$ (\ref{WsWu0}) up to the ($\eps \to 0$ limit of the) touch down point $J_{g-d}^{+, 0} \in \M^0_0$ (that has the same $(w,q)$-coordinates as $J_{g-o}^{+, 0}$ in the limit $\eps \to 0$). Note that the take off/touch down points $J_{g-0}^{+, 0}$/$J_{g-d}^{+, 0}$ differ essentially from their counterparts as $J_{\rm s-d}$/$J_{\rm s-0}$ (\ref{defJs-d}) considered so far: while $J_{\rm s-d}$/$J_{\rm s-0}$ concerned the evolution of $W^u(P^0) \cap W^s(\M^+_\eps)$/$W^s(P^0) \cap W^u(\M^+_\eps)$ along $\M^+_\eps$ in forwards/backwards `time', $J_{g-0}^{+, 0}$/$J_{g-d}^{+, 0}$ govern the orbits of $W^u(P^{+,s}) \cap W^s(\M^0_\eps)$/$W^s(P^{+,s}) \cap W^u(\M^0_\eps)$ along $\M^0_\eps$. Nevertheless, the coordinates of all take off/touch down points are at leading order determined by their $\eps \to 0$ limits (\ref{WsWu0}) with $w^\pm_h(0) = 9/(2 + 9a)$ (\ref{wafoc}). The next piece of $\Gamma_{\rm g-2f}$ consists of a symmetric part of a ($\cosh$-type) orbit along $\M^0_0$ of the (linear) slow reduced flow (\ref{eq:SRSzero}) up to the ($\eps \to 0$ limit of the) take off point $J_{g-o}^{0,+}$, which is again followed by a fast jump backwards along (a piece of) $W^s(\M^+_0) \cap W^u(\M^0_0)$ to the ($\eps \to 0$ limit of the) touch down point $J_{g-d}^{0,+} \in \M^+_0$. The final piece is the symmetrical counterpart of the first piece: the flow of (\ref{eq:SRSplus}) along $\M^+_0$ from the final touch down point back towards $P^{+,s}$ -- see Fig. \ref{fig:gapsandperiodics}a for a sketch of a projection of $\Gamma_{\rm s-2f}$ (without its fast jumps).
\\ \\
The (open) region $\SSS_{g-2f}$ is defined by those $(a,\Psi,\Phi,\Omega,\Theta)$-combinations for which $\Gamma_{\rm g-2f}$ can be constructed. We note that there cannot be points in the intersection of $\SSS_{s-2f}$ as defined in Theorem \ref{T:standing2fspots} and $\SSS_{g-2f}$: the (projections of the) take off/touch down points $J_{\rm s-d}$/$J_{\rm s-0}$ lie on the level set $\{\mathcal{H}_{0}^+(w,q) = \HH_{\rm s-d} \}$ for which $\HH_{\rm s-d} < \HH_0^{+,s}$, the value of $\mathcal{H}_{0}^+(w,q)$ for $P^{+,s}$ and its (un)stable manifolds (\ref{manifoldsP+s}). By construction, $J_{g-0}^{+, 0}, J_{g-d}^{+, 0} \subset \{\mathcal{H}_{0}^+(w,q) = \HH_0^{+,s} \}$ -- compare Fig. \ref{fig:3skeletons}a to Fig. \ref{fig:gapsandperiodics}a. This also implies that $\overline{\SSS_{s-2f}} \cap \overline{\SSS_{g-2f}} \neq \emptyset$, in fact, $\partial \SSS_{s-2f} \cap \partial \SSS_{g-2f} \supset \R_{\rm s-1f}$ as defined in Theorem \ref{T:stand1fronts}: both the homoclinic spots $\gamma_{\rm s-2f}(\xi)$ of Theorem \ref{T:standing2fspots} and the gaps $\gamma_{\rm g-2f}(\xi)$ of (upcoming) Theorem \ref{T:standing2fgaps} merge with the heteroclinic cycle spanned by the standing 1-front $\gamma_{\rm s-1f}(\xi)$ of Theorem \ref{T:stand1fronts} and its symmetrical counterpart as $\SSS_{s-2f} \cup \SSS_{g-2f}$ approaches $\R_{\rm s-1f}$.

\begin{theorem}
\label{T:standing2fgaps}
Let $(a,\Psi,\Phi,\Omega,\Theta) \in \SSS_{g-2f}$ and $\Gamma_{\rm g-2f} \subset \mathbb{R}^4$ be the singular (gap) skeleton constructed above. Then, there is for $\eps > 0$ sufficiently small a symmetric 2-front orbit $\gamma_{\rm g-2f}(\xi) = (b_{\rm g-2f}(\xi),p_{\rm g-2f}(\xi),w_{\rm g-2f}(\xi),q_{\rm g-2f}(\xi)) \subset W^u(P^{+,s}) \cap W^s(P^{+,s})$ of (\ref{eq:fast}) with $c=0$ that merges with $\Gamma_{\rm g-2f}$ as $\eps \downarrow 0$. The associated stationary pattern $(B(x,t), W(x,t)) = (b_{\rm g-2f}(x),w_{\rm g-2f}(x))$ in (\ref{eq:scaled}) represents a stationary localized bare soil gap embedded in vegetation.
\end{theorem}
\noindent
Of course the proof of this Theorem goes exactly along the lines of the proof of Theorem \ref{T:standing2fspots}. The main difference between the cases of (stationary, symmetric, homoclinic) 2-front spots and (stationary, symmetric, homoclinic) 2-front gaps is that there cannot be periodic orbits on $\M^0_\eps$ -- the slow reduced flow (\ref{eq:SRSzero}) on $\M^0_0$ is linear -- so that there cannot be any higher order localized gap patterns (as in Theorem \ref{T:higherorder2fronts} for localized spots). We refer to Figs. \ref{fig:basicpatterns}c and \ref{fig:numerics_fairy_circle} for numerical observations of  -- (most likely) stable -- localized 2-front spot and gap patterns in PDE (\ref{eq:scaled}).
\begin{remark}
\label{rem:2ftofasthom}
We refer to \cite{Kok,KWSW2006} for studies of the process of a homoclinic 2-front orbit between a slow manifold $\M^1_\eps$ and a second slow manifold $\M^2_\eps$ detaching from $\M^2_\eps$ to become a slow-fast homoclinic to $\M^1_\eps$ that only makes 1 homoclinic excursion through the fast field (instead of 2 fast heteroclinic jumps between $\M^1_\eps$ and $\M^2_\eps$). The focus of \cite{Kok} is on the (exchange of) stability between these 2 types of homoclinic patterns and the associated bifurcations -- especially as localized stripes in 2 space dimensions. In the present work the situation is somewhat more involved than in \cite{Kok,KWSW2006}, since the skeleton structures as sketched in Figs. \ref{fig:3skeletons} and \ref{fig:gapsandperiodics}a must become asymptotically small in this transition -- which is not necessary in the setting of \cite{Kok,KWSW2006}. Notice that this implies that here the $W$-component of the homoclinic pulse becomes `small' -- a certain well-defined magnitude in $\eps$ -- during this transition, but that this is not the case for the $B$-component since the orbit still has to make (almost) a full jump between $\M^0_\eps$ and $\M^+_\eps$. See Figs. \ref{fig:numerics_2_front_standing}, \ref{fig:numerics_1_pulse} and \ref{fig:projection_pulse_2front_manifolds} in section \ref{sec:sims}.
\end{remark}
\subsection{Spatially periodic multi-front patterns}
\label{sec:perpatt}
A stationary, non-degenerate, symmetric homoclinic pulse solution of a (reversible) reaction-diffusion system (defined for $x \in \mathbb{R}$) -- such as the spots and gaps of Theorems \ref{T:standing2fspots} and \ref{T:standing2fgaps} -- must be the `endpoint' (within phase space) of a continuous family -- a `band' -- of spatially periodic patterns (as the period/wave length $\to \infty$) -- see for instance \cite{Dreview}. Systems (\ref{eq:scaled}) and (\ref{eq:fast}) indeed have large families of spatially periodic solutions.

\begin{theorem}
\label{T:periodics}
Let $(a,\Psi,\Phi,\Omega,\Theta) \in \SSS_{s-2f} \cup \SSS_{g-2f} \cup \R_{\rm s-1f}$ (Theorems \ref{T:standing2fspots}, \ref{T:standing2fgaps} and \ref{T:stand1fronts}) and let $c =0$. Let there be a $\rho \in \mathbb{R}$, $\rho \neq 0$, such that there is a solution of the reduced slow flow (\ref{eq:SRSzero}) on $\M^0_0$ that connects $(\rho,9/(2 + 9a))$ to $(-\rho,9/(2 + 9a))$ over the branch $\Gamma^0_\rho$ (by definition) and a solution of the (projected) reduced flow (\ref{eq:SRSplus}) on $\M^+_0$ that connects $(-\rho,9/(2 + 9a))$ back to $(\rho,9/(2 + 9a))$ over $\Gamma^+_\rho$ -- see Fig. \ref{fig:gapsandperiodics}b. Then, for $\eps > 0$ sufficiently small, there is a periodic solution $\gamma_{{\rm mf},\rho}(\xi) = (b_{{\rm mf},\rho}(\xi),\rho_{{\rm mf},\rho}(\xi),w_{{\rm mf},\rho}(\xi),q_{{\rm mf},\rho}(\xi))$ of (\ref{eq:fast}) that merges in the limit $\eps \to 0$ with the skeleton structure spanned by $\Gamma^0_\rho$, the fast jump over $W^u(\mathcal{M}_0^{0}) \cap W^{s}(\mathcal{M}_0^{+})$ with $w^+_h = 9/(2 + 9a)$, $q=-\rho$ (\ref{WsWu0}), $\Gamma^0_\rho$ and the fast jump back over $W^u(\mathcal{M}_0^{+}) \cap W^{s}(\mathcal{M}_0^{0})$ with $w^-_h = 9/(2 + 9a)$, $q=\rho$. The associated stationary pattern $(B(x,t), W(x,t)) = (b_{{\rm mf},\rho}(x),w_{{\rm mf},\rho}(x))$ in (\ref{eq:scaled}) represents a stationary spatially periodic multi-pulse spot/gap pattern.
\end{theorem}
\noindent
Note that if $\gamma_{{\rm mf},\rho}(\xi)$ exists for a certain $\rho^*$, there clearly must be a neighborhood of $\rho^*$ for which $\gamma_{{\rm mf},\rho}(\xi)$ also exists: the periodic solutions of (\ref{eq:scaled})/(\ref{eq:fast}) indeed come in continuous families/bands \cite{Dreview}. Typically, there is a `subband' of stable periodic patterns -- see Figs. \ref{fig:basicpatterns}d and \ref{fig:numerics_periodic_multifront} for examples of numerically stable patterns $(b_{{\rm mf},\rho}(x),w_{{\rm mf},\rho}(x))$ in (\ref{eq:scaled}).
\\ \\
{\bf Proof of Theorem \ref{T:periodics}.} This proof can be set up very much along the lines of the proofs of similar results -- the existence of spatially periodic patterns in the (generalized) Gierer-Meinhardt equation -- in \cite{DKvdPloeg01}, therefore we restrict ourselves to the essential ingredients of the proof here.
\\ \\
The approach is similar to that of the proof of Theorem \ref{T:standing2fspots}: we construct an orbit that intersects the plane $\{p=q=0\}$ and apply the reversibility symmetry (\ref{eq:symmc}) -- with $c=0$. However, unlike for the homoclinic orbits in Theorem \ref{T:standing2fspots}, we do not `start out' -- as $\xi \to -\infty$ -- at the critical point $P^0 \in \M^0_\eps$, but choose the initial condition of the orbit-to-be-constructed at the -- exponentially short -- interval $\I_{\overline{b}} = \{p=q=0, w = w^0_\rho, b = \overline{b} \in (0,\overline{b}_M)\}$, where $w^0_\rho$ is the midpoint of $\Gamma^0_\rho$, i.e. $(w^0_\rho,0) = \Gamma^0_\rho \cap \{q=0\}$ and $\overline{b}_M$ is exponentially small in $\eps$, so that orbits $\gamma_{\rho,\overline{b}}(\xi)$ with $\gamma_{\rho,\overline{b}}(0) \in  \I_{\overline{b}}$ remain exponentially close to $\M^0_\eps$ over $\OO(1)$ distances in $b, p$ and $\OO(1/\eps)$ `time' $\xi$ (more precise, and as in the proof of Theorem \ref{T:standing2fspots}: $\overline{b} = \OO({\rm exp}(-\lambda_{f,0}(w^0_\rho)/\eps)$ with $\lambda_{f,0}(w_0) = \sqrt{1-aw_0}$, the unstable eigenvalue of the reduced fast flow (\ref{eq:fastblue}) -- with $c=0$ -- associated to $(0,0)$). Consider the 2-dimensional `strip' $\T_{\rho,\overline{b}}$ spanned by all $\gamma_{\rho,\overline{b}}(\xi)$, $\overline{b} \in  (0,\overline{b}_M)$: as it takes off from $\M^0_\eps$, it is exponentially close to $W^u(\M^0_\eps)$ and thus intersects $W^s(\M^+_\eps)$ transversely along the orbit $\gamma_{\rho,\overline{b}^*}(\xi)$ (the intersection $\T_{\rho,\overline{b}} \cap W^s(\M^+_\eps)$ is 1-dimensional and thus an orbit of (\ref{eq:fast})).
\\ \\
Clearly, $W^s(\M^+_\eps)$ cuts $\T_{\rho,\overline{b}}$ into 2 parts -- distinguished by $\overline{b} \lessgtr \overline{b}^*$ -- due to the `fast saddle' structure of $\M^+_\eps$ (as in the proof of Theorem \ref{T:standing2fspots}). Orbits $\gamma_{\rho,\overline{b}}(\xi) \subset \T_{\rho,\overline{b}}^{\rm r} \subset \T_{\rho,\overline{b}}$ cross through $\{p=0\}$, turn around (in their $b$-components) and return back towards $\M^0_\eps$; the $b$-components of orbits $\gamma_{\rho,\overline{b}}(\xi) \subset \T_{\rho,\overline{b}} \backslash \T_{\rho,\overline{b}}^{\rm r}$ -- the complement of $\T_{\rho,\overline{b}}^{\rm r}$ -- increase beyond $\M^+_\eps$. It depends on the relative magnitudes of $w^0_\rho$ and $9/(2+9a)$ whether $\T_{\rho,\overline{b}}^{\rm r}$ is determined by $\overline{b} \in (0,\overline{b}^*)$ or by $\overline{b} \in (\overline{b}^*,\overline{b}_M)$. If $w^0_\rho > 9/(2+9a)$ -- as in Fig. \ref{fig:gapsandperiodics} -- orbits $\gamma_{\rho,\overline{b}}(\xi)$ that take off `too soon' -- i.e. with $\overline{b} > \overline{b}^*$ -- have $w > 9/(2+9a)$ at take off (Fig. \ref{fig:gapsandperiodics}). Since the unstable manifold $W^u((0,0))$ of the (planar) fast reduced system (\ref{eq:fastblue}) with $w_0 > 9/(2+9a)$ contains a closed homoclinic orbit, it follows that orbits $\gamma_{\rho,\overline{b}}(\xi)$ with $\overline{b} > \overline{b}^*$ follow such a homoclinic orbit through the fast field (at leading order in $\eps$). Hence, they pass through $\{p=0\}$ and turn back towards $\M^+_\eps$: $\T_{\rho,\overline{b}}^{\rm r}$ is spanned by $\gamma_{\rho,\overline{b}}(\xi)$ with $\overline{b} \in  (\overline{b}^*,\overline{b}_M)$. For simplicity, we only consider this case (the arguments run along exactly the same lines in the case that $w^0_\rho < 9/(2+9a)$ and $\T_{\rho,\overline{b}}^{\rm r}$ is determined by $\overline{b} \in (0,\overline{b}^*)$).
\\ \\
We can now copy the main (geometrical) argument of the proof of Theorem \ref{T:standing2fspots}: if $\gamma_{\rho,\overline{b}}(\xi) \subset \T_{\rho,\overline{b}}^{\rm r}$ is too far removed from -- but still exponentially close to --  $\gamma_{\rho,\overline{b}^*}(\xi) \subset W^s(\M^+_\eps)$ -- i.e. if $\overline{b} - \overline{b}^*$ is too large -- it will follow $\M^+_\eps$ -- and thus $\Gamma^+_\rho$ -- over a relatively short distance (Fig. \ref{fig:gapsandperiodics}), take off again from $\M^+_\eps$ and thus cross through $\{p=0\}$ before reaching $\{q=0\}$. By decreasing $\overline{b} - \overline{b}^*$, one can keep $\gamma_{\rho,\overline{b}^*}(\xi)$ sufficiently long close to $\M^+_\eps$ that it first passes through $\{q=0\}$ before crossing $\{p=0\}$: there is a value $\overline{b} = b_\rho$ such that $\gamma_{\rho,b_\rho}(L/2) \in \{p=q=0\}$ for certain $L > 0$ (in fact, $L = \OO(1/\eps)$). It follows by (\ref{eq:symmc}) -- with $c=0$ -- that $\gamma_{{\rm mf},\rho}(\xi) = \gamma_{\rho,b_\rho}(\xi)$, a periodic orbit with period $L$. \hfill $\Box$
\\ \\
As in the case of the homoclinic spots patterns -- Theorem \ref{T:higherorder2fronts} -- we may immediately conclude that there are countably many (families of) higher order periodic patterns if $\Gamma^+_\rho$ (as defined in Theorem \ref{T:periodics}) is part of a periodic orbit on $\M^+_0$ (determined by the level set $\{\mathcal{H}_{0}^+(w,q) = \HH_\rho \}$) -- see Fig. \ref{fig:3skeletons}c. By steadily decreasing $|\overline{b}-\overline{b}^*|$, the orbit $\gamma_{\rho,\overline{b}}(\xi) \subset \T_{\rho,\overline{b}}^{\rm r}$ can be made to pass arbitrarily many times through $\{q=0\}$ before taking off from $\M^+_\eps$.
\begin{corollary}
\label{cor:higherorderperiodics}
Let $(a,\Psi,\Phi,\Omega,\Theta) \in \SSS_{s-2f} \cup \SSS_{g-2f} \cup \R_{\rm s-1f}$ such that (\ref{asshomM+eps}) holds. Let $\rho$, $\Gamma^0_\rho \subset \M^0_0$ and $\Gamma^+_\rho \subset \M^+_0$ be as defined in Theorem \ref{T:periodics}, with $\Gamma^+_\rho$ such that it is part of a closed orbit on $\M^+_\eps$ determined by the level set $\{\mathcal{H}_{0}^+(w,q) = \HH_\rho \}$. Then, for $\eps > 0$ sufficiently small, there is a countable family of symmetric multi-front periodic orbits $\gamma_{{\rm mf},\rho}^i(\xi) = (b_{{\rm mf},\rho}^i(\xi),p_{{\rm mf},\rho}^i(\xi),w_{{\rm mf},\rho}^i(\xi),q_{{\rm mf},\rho}^i(\xi))$ of (\ref{eq:fast}) with $c=0$ ($i=1,2,...$) that merges in the limit $\eps \to 0$ with the extended skeleton structure spanned by $\Gamma^0_\rho$, the fast jump over $W^u(\mathcal{M}_0^{0}) \cap W^{s}(\mathcal{M}_0^{+})$ with $w^+_h = 9/(2 + 9a)$, $q=-\rho$ (\ref{WsWu0}), the full closed orbit of $\{\mathcal{H}_{0}^+(w,q) = \HH_\rho \}$ that contains $\Gamma^+_\rho$ and the fast jump back over $W^u(\mathcal{M}_0^{+}) \cap W^{s}(\mathcal{M}_0^{0})$ with $w^-_h = 9/(2 + 9a)$, $q=\rho$. The associated stationary patterns $(B(x,t), W(x,t)) = (b_{{\rm mf},\rho}^i(x),w_{\rm s-2f}^i(x))$ in (\ref{eq:scaled}) are symmetric periodic spot/gap patterns with an increasing number of oscillations in the vegetated areas.
\end{corollary}
Finally, we note that the families of `higher order' periodic patterns $\gamma_{{\rm mf},\rho}^i(\xi)$ are only the first of further -- more complex -- families containing periodic (and aperiodic) patterns of increasing complexity. We refer to \cite{DKvdPloeg01} for the precise settings and proofs, here we only give a sketch of one specific example. However, this sketch provides the main ideas by which all further orbits may be constructed.
\\ \\
Let $b_\rho^i \in (0,\overline{b}_M)$ be such that the $i$-th periodic pattern $\gamma_{{\rm mf},\rho}^i(\xi)$ of Corollary \ref{cor:higherorderperiodics} is given by $\gamma_{\rho,\overline{b}}(\xi)$ with $\overline{b} = b_\rho^i$ (see the proof of Theorem \ref{T:periodics}). We can now choose $\overline{b}$ so close to $b_\rho^i$ that $\gamma_{\rho,\overline{b}}(\xi) \subset \T_{\rho,\overline{b}}^{\rm r}$ follows $\gamma_{{\rm mf},\rho}^i(\xi)$ along its $i$ circuits over $\M^+_\eps$ -- with $i \geq 1$ -- and its jump back to $\M^0_\eps$. Since $\overline{b} \neq b_\rho^i$, $\gamma_{\rho,\overline{b}}(\xi)$ does not close as it passes along $\I_{\overline{b}}$, instead it keeps on following $\gamma_{{\rm mf},\rho}^i(\xi)$ as it makes it second jump towards $\M^+_\eps$. By the approach of the proofs of Theorems \ref{T:standing2fspots} and \ref{T:periodics}, we can now tune $\overline{b}$ so that it has its second take off from $\M^+_\eps$ precisely and that it passes through $\{p=q=0\}$ while following $\Gamma^+_\rho$ (without making any further circuits over the periodic orbit on $\M^+_\eps$ that contains $\Gamma^+_\rho$). It follows by the application of the reversibility symmetry (\ref{eq:symmc}) that for this value of $\overline{b}$, $\gamma_{\rho,\overline{b}}(\xi)$ is a symmetric periodic orbit that `starts' at $\Gamma^0_\rho \subset \M^0_\eps$, jumps to $\M^+_\eps$ to make $i$ circuits along $\M^+_\eps$, jumps back again to $\M^0_\eps$, follows $\Gamma^0_\rho \subset \M^0_\eps$ to return again to $\M^+_\eps$ where it follows $\Gamma^+_\rho$ and subsequently immediately jumps back again to $\M^0_\eps$ -- from which it repeats the same path, etc.. Note that the associated periodic pattern in (\ref{eq:scaled}) consists of an alternating array of 2 different types of localized vegetation spots. Clearly, this procedure can be further refined to establish the existence of patterns containing arbitrary arrays of arbitrarily many different types of vegetation spots -- under the conditions of Corollary \ref{cor:higherorderperiodics}.
\begin{remark}
\label{rem:travpatts}
We decided to focus in this paper on stationary 2- and multi-front patterns. Of course, (\ref{eq:scaled}) also exhibits traveling multi-front patterns -- see for instance Fig. ~\ref{fig:numerics_2_front_traveling} in which a vegetation spot travels towards a stationary, stable (and attracting) spot of the type established by Theorem \ref{T:standing2fspots}. System (\ref{eq:fast}) can also have homoclinic orbits to $P^0$ for (certain specific values of) $c \neq 0$, i.e. vegetation spots may be traveling with constant speed (without changing shape). An approach along the lines of \cite{DvHeijsterK09} indicates that bifurcations to traveling spots appear when the touch down manifold $\I_{\rm down}$ (Lemma \ref{lem:Idown}) is tangent to a level set $\{\mathcal{H}_{0}^+(w,q) = \HH \}$ of the slow reduced flow on $\M^+_0$ with $c=0$ at the (non-transversal) intersection $\I_{\rm down} \cap \{\mathcal{H}_{0}^+(w,q) = \HH \}$ (recall that $\I_{\rm down}$ is parameterized by $c$). A similar property holds for bifurcations of stationary spatially periodic multi-front patterns into traveling spatially periodic (wave train) patterns. A simple investigation of the relative orientations of $\I_{\rm down}$ and the various possible phase plane configurations of the slow reduced flow on $\M^+_0$ shows that there indeed are parameter combinations $ (a, \Psi, \Phi, \Omega, \Theta)$ at which these bifurcations into traveling 2-/multi-fronts must occur. This bifurcation may have relevant ecological implications, nevertheless, we refrain from going into the details here (and leave this to future work) -- see also section \ref{sec:disc}.
\end{remark}

\section{Simulations and discussion}
\label{sec:SimDis}

\subsection{Simulations}
\label{sec:sims}
The motivation for the numerical simulations presented in this section is threefold: 1) to illustrate some of the analytic results of the previous sections (without doing a systematic search for all constructed patterns) 2) to give a brief outlook beyond the worked out analysis to solution types that may be constructed by the geometric set-up developed here and, finally, 3) to give a flavor of the rich dynamics that PDE (\ref{eq:meron})/(\ref{eq:scaled}) exhibits. All numerical simulations have been carried using MATLAB's `pdepe' routine. The corresponding parameter settings are specified in the captions of the figures. Almost all figures show a snapshot in time of the spatial profile of the PDE solution/pattern after it converged to a stationary or uniformly traveling solution.
\\

\begin{figure}[h!]
	\centering
	\begin{subfigure}[b]{0.48\textwidth}
		\centering
		\includegraphics[width= \textwidth]{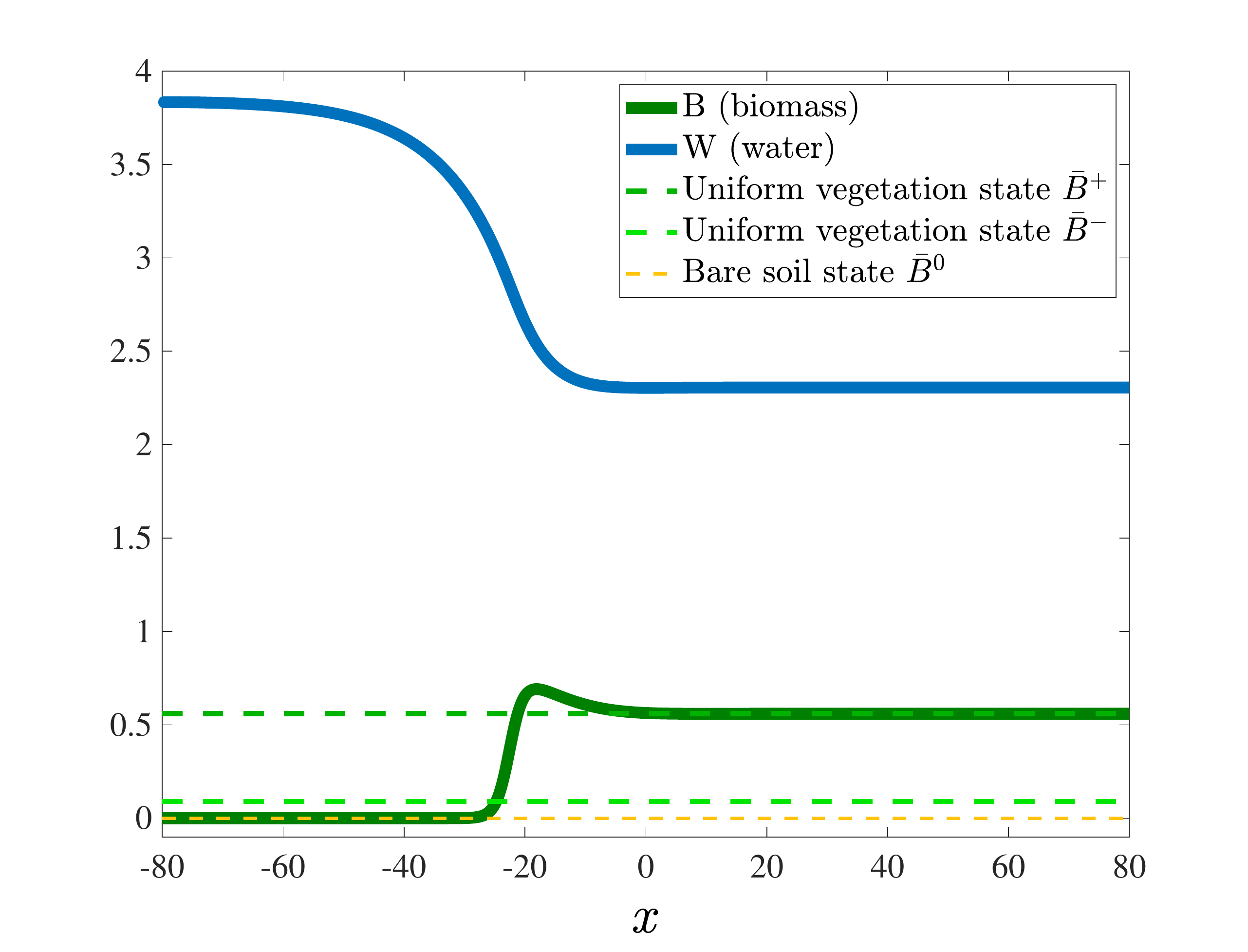}
		\caption{\label{fig:numerics_1_front_ehud}}
	\end{subfigure}%
	\hspace{.4cm}
	\begin{subfigure}[b]{0.48\textwidth}
		\centering
		\includegraphics[scale=0.4]{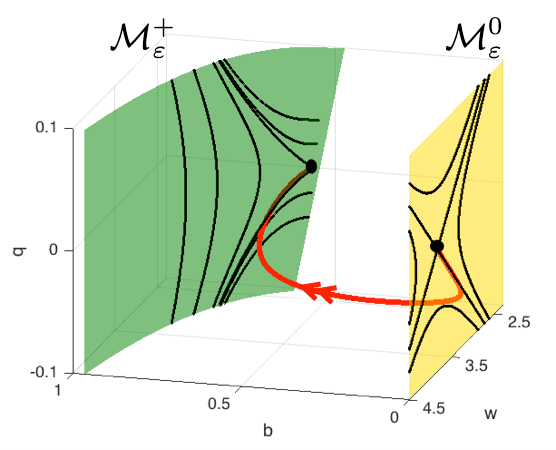}
		\caption{\label{fig:numerics_1_front_ehud_phase_space}}
	\end{subfigure}%
	\caption{ {\bf (a)} Spatial profiles of the $B$- and $W$-components of a traveling front solution of the original, unscaled model (\ref{eq:meron}) -- corresponding to a heteroclinic orbit of (\ref{eq:fast}) as established by Theorem \ref{T:primaryfronts} -- together with the 3 background states. {\bf (b)} The profile from {\bf (a)} as a projection in $(b,w,q)$-space where $ q $ was computed by numerical differentiation. {\it Parameter settings} in (\ref{eq:meron}): $ P = 180, \Lambda = 0.9, K = 0.4, E =  18, M = 15, N = 15, R = 0.7, \Gamma = 12 , D_W = 150, D_B = 1.2 $, corresponding to $\eps^2 = 0.008, a \approx 0.187, \Psi \approx 3.84, \Phi = 1, \Omega \approx 0.235, \Theta \approx 1.71$ in (\ref{eq:scaled}).
	}
	\label{fig:numerics_front_ehud_both}
\end{figure}

The opening figure of this section, Fig. ~\ref{fig:numerics_1_front_ehud}, can be seen as a binding element between the papers that have a more ecological emphasis and motivated the present work -- see \cite{Zelnik2018eco_ind, Zelnik2015} and the references therein -- and the analysis here. It displays a traveling front solution as established by Theorem \ref{T:primaryfronts} for a parameter regime comparable to the one from \cite{Zelnik2015} (with slight adjustment in the parameters to compensate for the choice of a 1-D model -- \ref{Appendix:derivation}). This profile is then shown in Fig. \ref{fig:numerics_1_front_ehud_phase_space} as a projection in $(b,w,q)$-space to illustrate that it indeed starts on the slow manifold $\mathcal{M}^0_\eps$ and then jumps to the $\mathcal{M}^{+}_\eps$ slow manifold. As established by the analysis, the solution first follows the unstable manifold associated to the bare soil state $(\bar{B}^0,\bar{W}^0)$ (as a solution of (\ref{eq:fast})), makes a fast excursion through the fast field to then touch down on  $ \mathcal{M}^{+}_\eps $ following the stable manifold associated to the uniform vegetation state $(\bar{B}^+,\bar{W}^+)$. Note, of course, that this figure contains two approximations: first, the manifold $ \mathcal{M}^{+}_\eps$ is only accurate up to second order in $ \varepsilon $ -- see section \ref{sec:slowfull} -- while the flows on $ \mathcal{M}^{0}_\eps$ and $ \mathcal{M}^{+}_\eps$ are computed numerically (using MATLAB routines).
\\ \\

\begin{figure}[h!]
	\centering
	\begin{subfigure}[b]{0.48\textwidth}
		\centering
		\includegraphics[width= \textwidth]{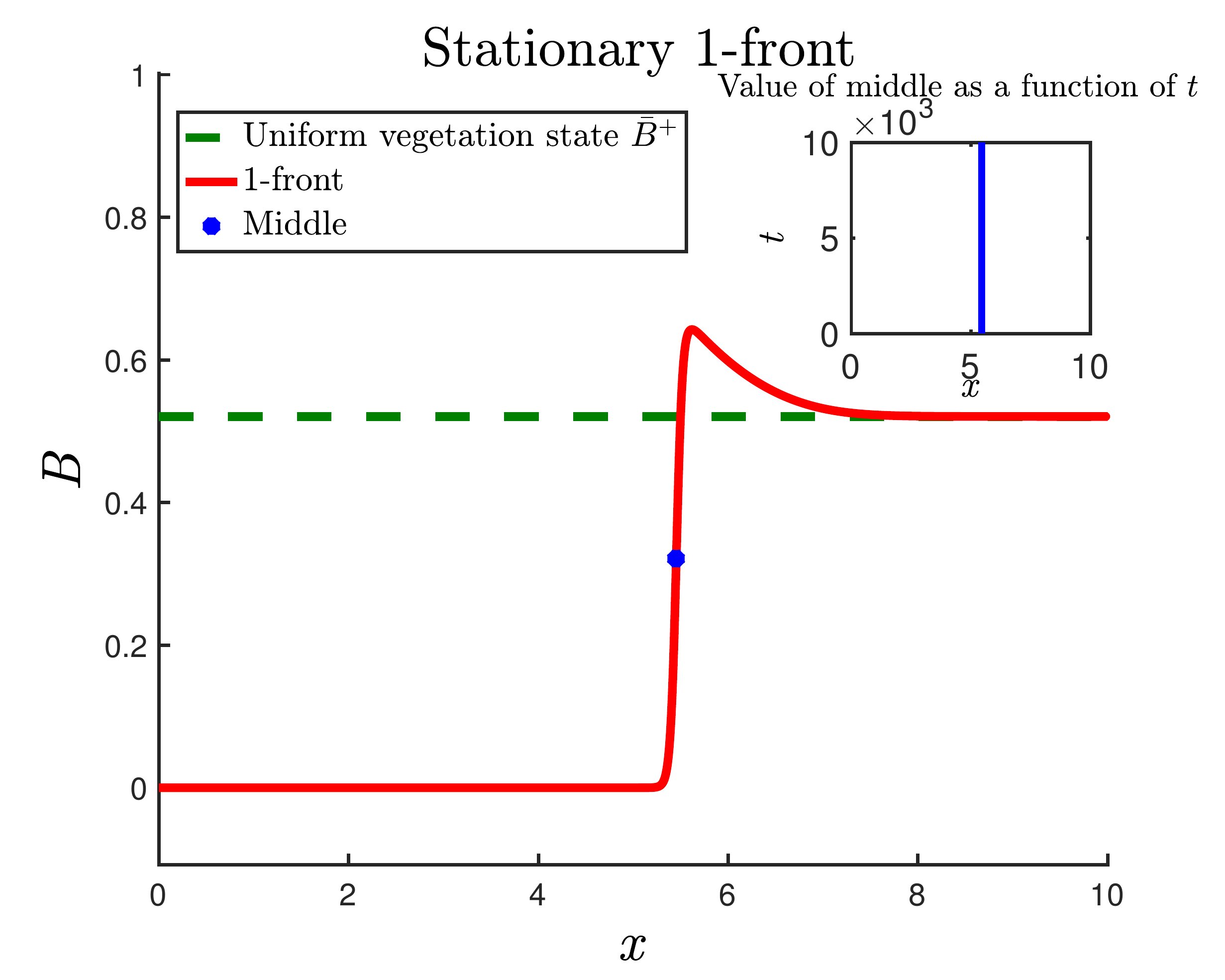}
		\caption{\label{fig:numerics_1_front_standing} $ a =0.0008,
			\Phi = 0.3, \Omega=0.1, \Theta = 0.2, \varepsilon = \sqrt{0.005}$ and $\Psi = \Psi_{\rm s-1f} =  1.6226$ (so that  $(a,\Phi,\Psi_{\rm s-1f},\Omega,\Theta) \in \R_{\rm s-1f}$, Theorem \ref{T:stand1fronts}).}
	\end{subfigure}%
	\hspace{0.4cm}
	\begin{subfigure}[b]{0.48\textwidth}
		\centering
		\includegraphics[width= \textwidth]{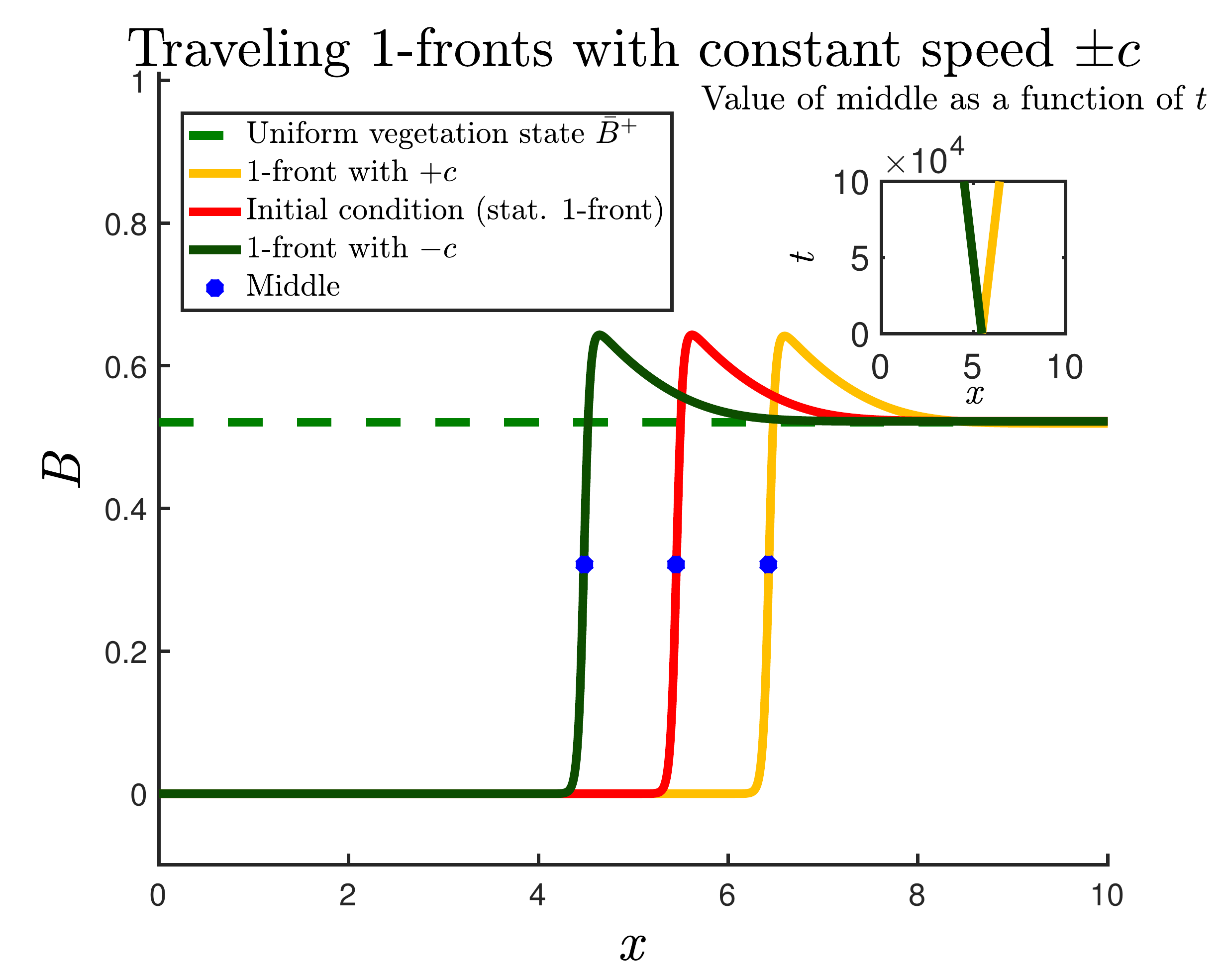}
		\caption{\label{fig:numerics_1_front_traveling} $ a =0.0008,
			\Phi = 0.3, \Omega=0.1, \Theta = 0.2, \varepsilon = \sqrt{0.005}$ and $\Psi_{+c} = 1.6205 < \Psi_{\rm s-1f}$, $\Psi_{-c} =    1.6248 > \Psi_{\rm s-1f}$.}
	\end{subfigure}
	\vspace{0.5cm}
	\begin{subfigure}[b]{0.48\textwidth}
		\centering
		\includegraphics[width= \textwidth]{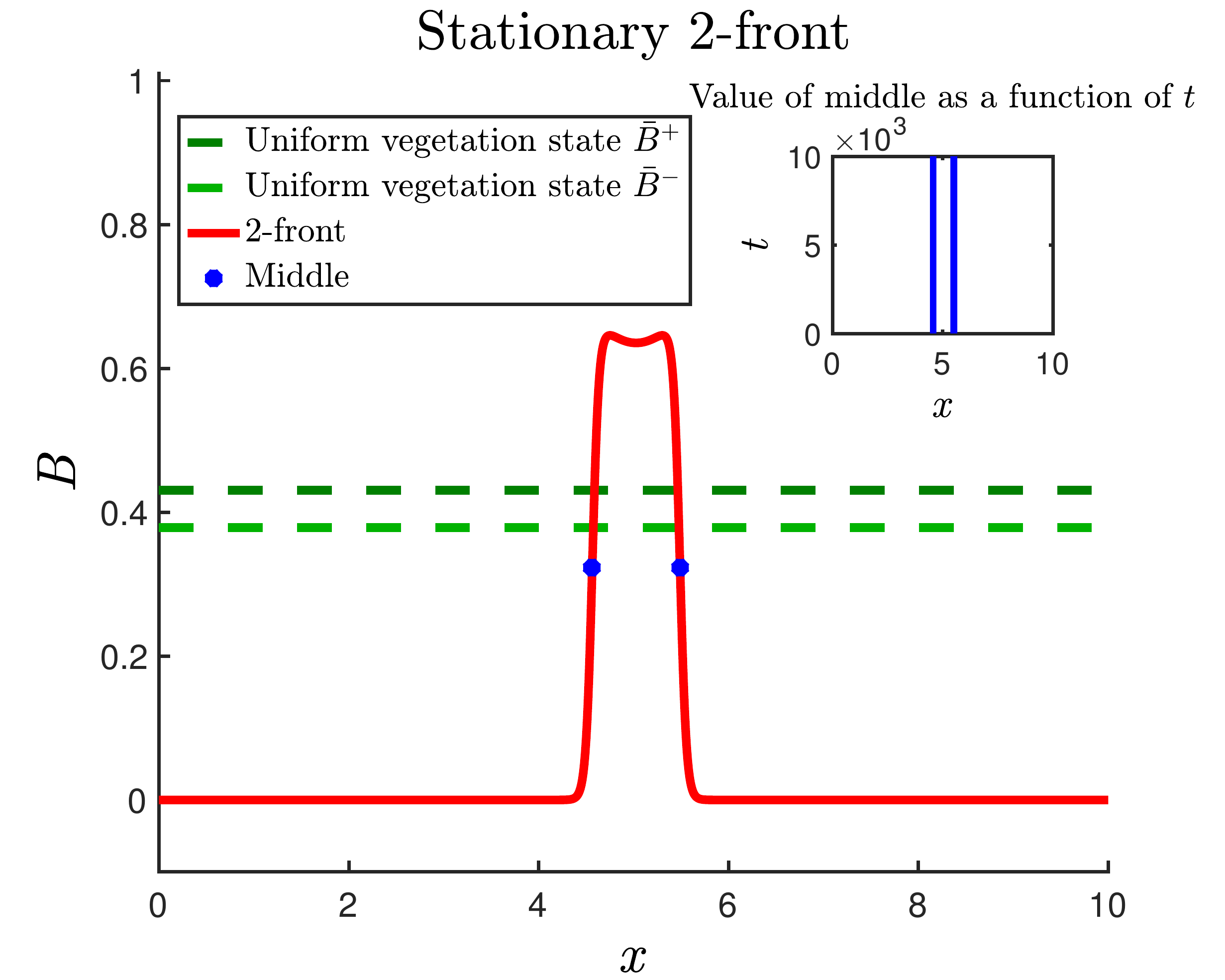}
		\caption{\label{fig:numerics_2_front_standing} $ a =0.032, \Psi =  1.3714 ,
			\Phi = 0.3, \Omega=0.1, \Theta = 0.2, \varepsilon = \sqrt{0.005}$.}
	\end{subfigure}%
	\hspace{0.4cm}
	\begin{subfigure}[b]{0.48\textwidth}
		\includegraphics[width= \textwidth]{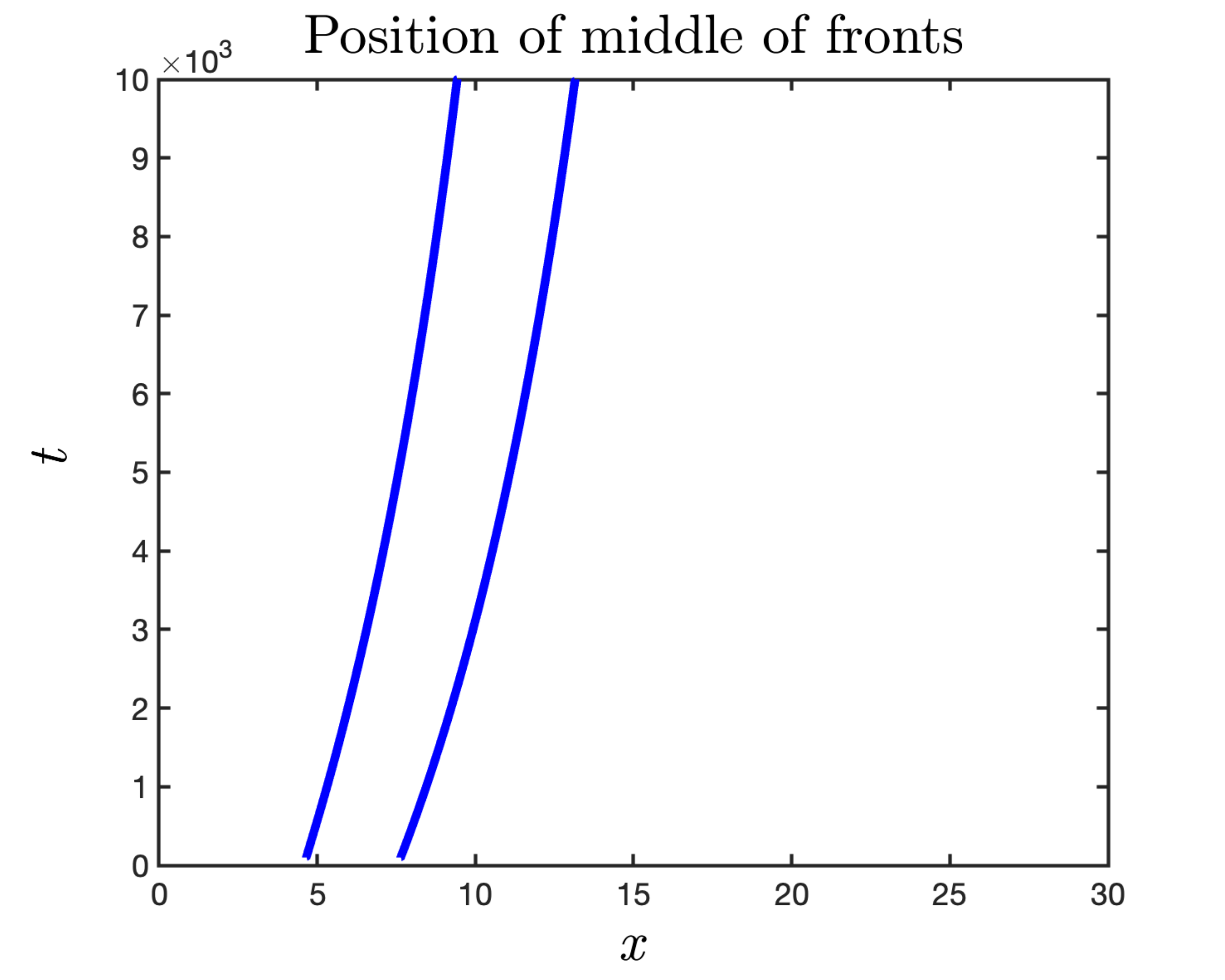}
		\caption{\label{fig:numerics_2_front_traveling} $ a =0.00175, \Psi =  1.6762 ,
			\Phi = 0.3, \Omega=0.1, \Theta = 0.7, \varepsilon = \sqrt{0.1}$.}
	\end{subfigure}%
	\caption{ {\bf (a)} A heteroclinic stationary 1-front pattern of (\ref{eq:scaled}); {\bf (b)} Two traveling 1-front patterns connecting the bare soil state to a homogeneous vegetation state, one invading the bare soil ($c < 0$), the other the vegetation state ($c > 0$); {\bf (c)} A homoclinic stationary 2-front spot pattern; {\bf (d)} Evolution of the middles of the 2 interacting fronts of an evolving 2-front pattern.}
	\label{fig:numerics_standing_and_traveling_waves}
\end{figure}
\begin{figure}[h!]
	\centering
	\begin{subfigure}[b]{0.48\textwidth}
		\centering
		\includegraphics[width= \textwidth]{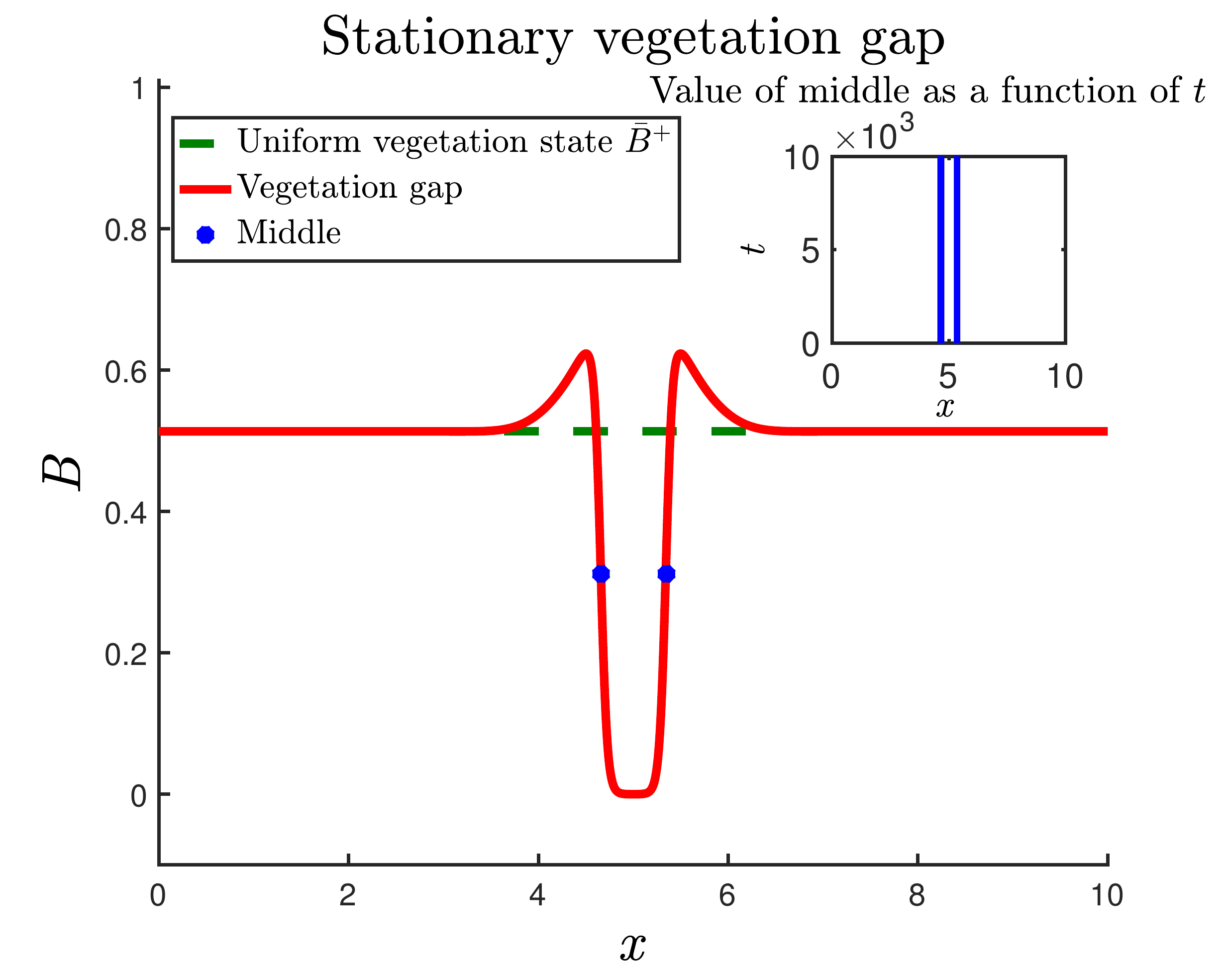}
		\caption{\label{fig:numerics_fairy_circle} $ a =0.032, \Psi =  1.2762 ,
			\Phi = 0.3, \Omega=0.1, \Theta = 0.2, \varepsilon = \sqrt{0.005}$.}
	\end{subfigure}%
	\hspace{0.4cm}
	\begin{subfigure}[b]{0.48\textwidth}
		\centering
		\includegraphics[width= \textwidth]{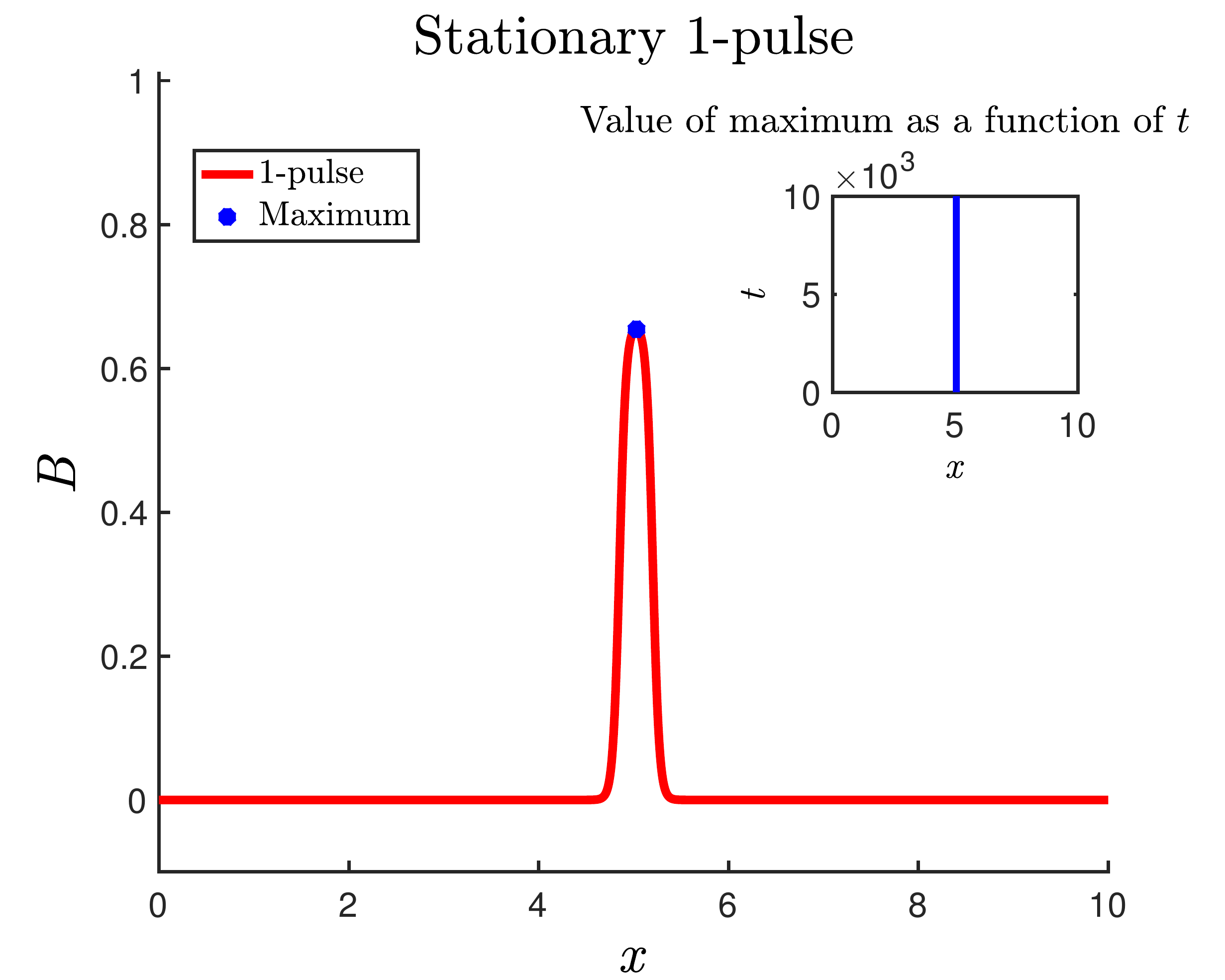}
		\caption{\label{fig:numerics_1_pulse} $ a =0.032, \Psi =  1.2762 ,
			\Phi = 0.3, \Omega=0.1, \Theta = 0.2, \varepsilon = \sqrt{0.005}$.}
	\end{subfigure}%
	\caption{{\bf (a)} A stationary homoclinic vegetation gap pattern (of fairy circle type) that is asymptotic to the stable homogeneous vegatation state $(\bar{B}^+,\bar{W}^+)$. {\bf (b)} A stationary homoclinic spot solution of (classical) 1-pulse Gierer-Meinhardt/Gray-Scott type.}
	\label{fig:numerics_fairy_circle_and_1_pulse}
\end{figure}
\begin{figure}[h!]
	\centering
	\begin{subfigure}[b]{0.48\textwidth}
		\centering
		\includegraphics[width= \textwidth]{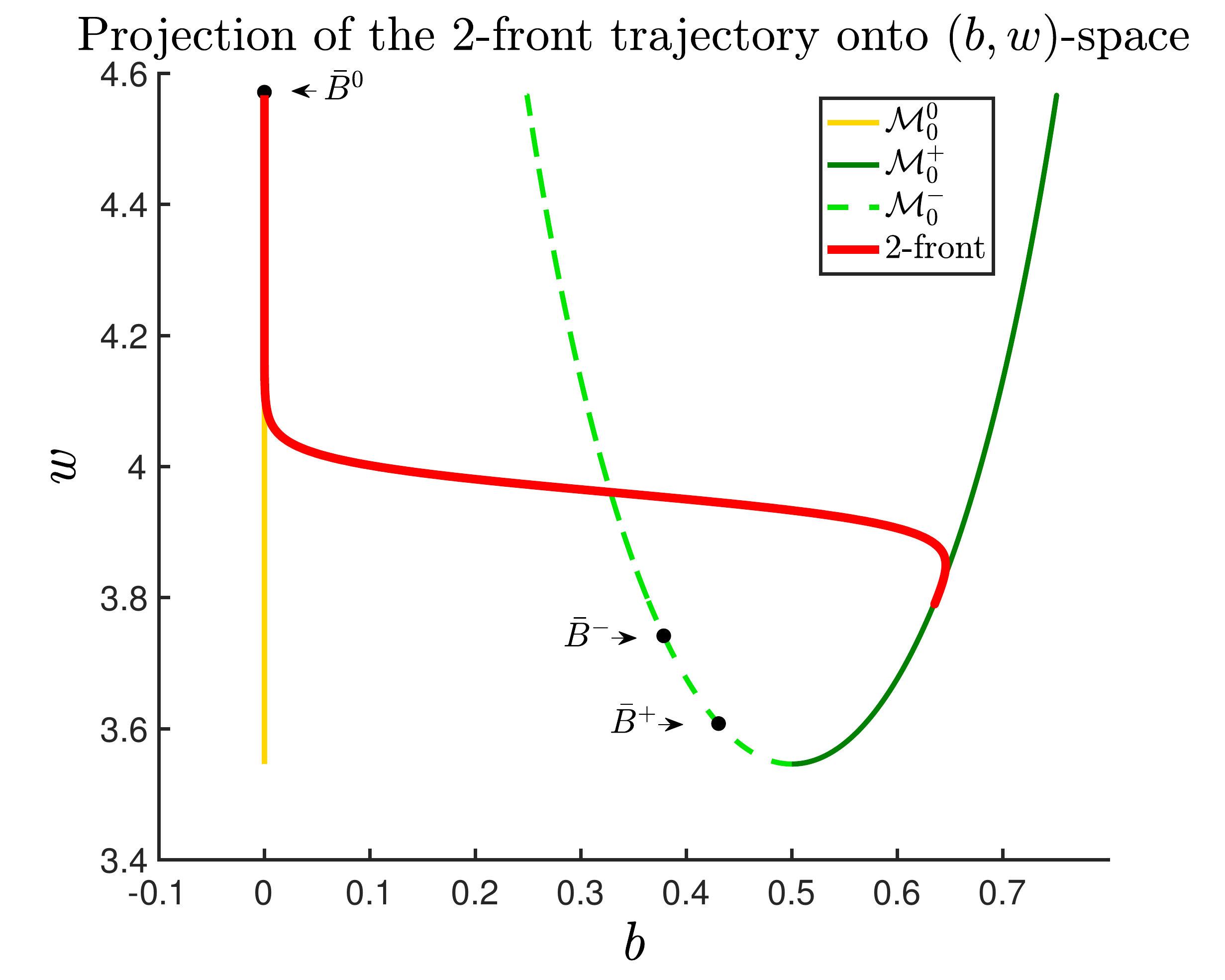}
	\end{subfigure}%
	\hspace{0.4cm}
	\begin{subfigure}[b]{0.48\textwidth}
		\centering
		\includegraphics[width= \textwidth]{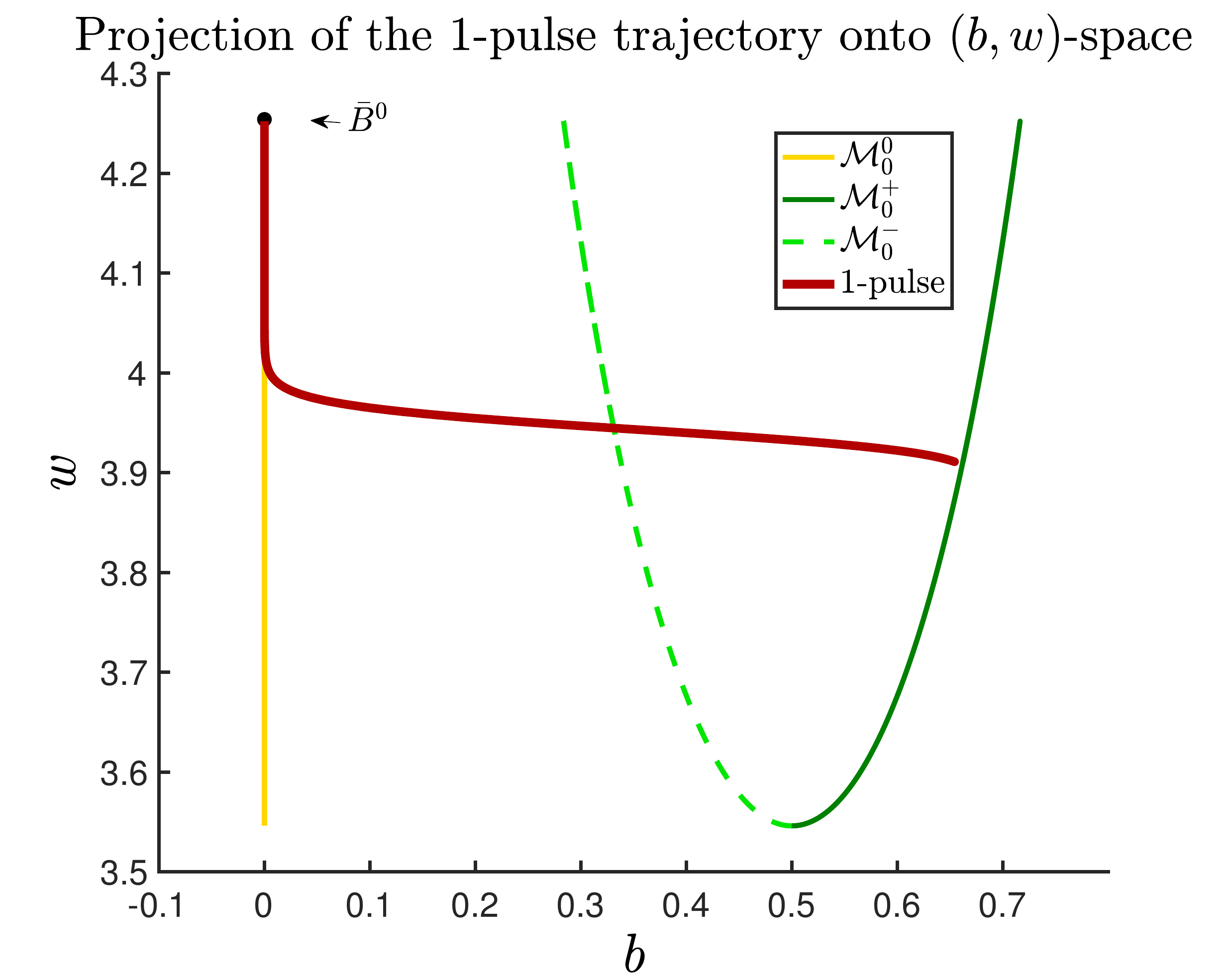}
	\end{subfigure}%
	\caption{Projections into the $(w,b)$-plane of the homoclinc 2-front and 1-pulse solutions of Figs.  ~\ref{fig:numerics_2_front_standing} and ~\ref{fig:numerics_1_pulse} together with the slow manifolds $\M^0_\eps$ and $\M^\pm_\eps$.  Note that the trajectories are symmetric around the middle of the 2-front/1-pulse -- due to the reversibility symmetry (\ref{eq:symmc}) of (\ref{eq:scaled})/(\ref{eq:fast}) -- which results in the red branches shown. }
	\label{fig:projection_pulse_2front_manifolds}
\end{figure}
\begin{figure}[h!]
	\centering
	\begin{subfigure}[b]{0.48\textwidth}
		\centering
		\includegraphics[width=\textwidth]{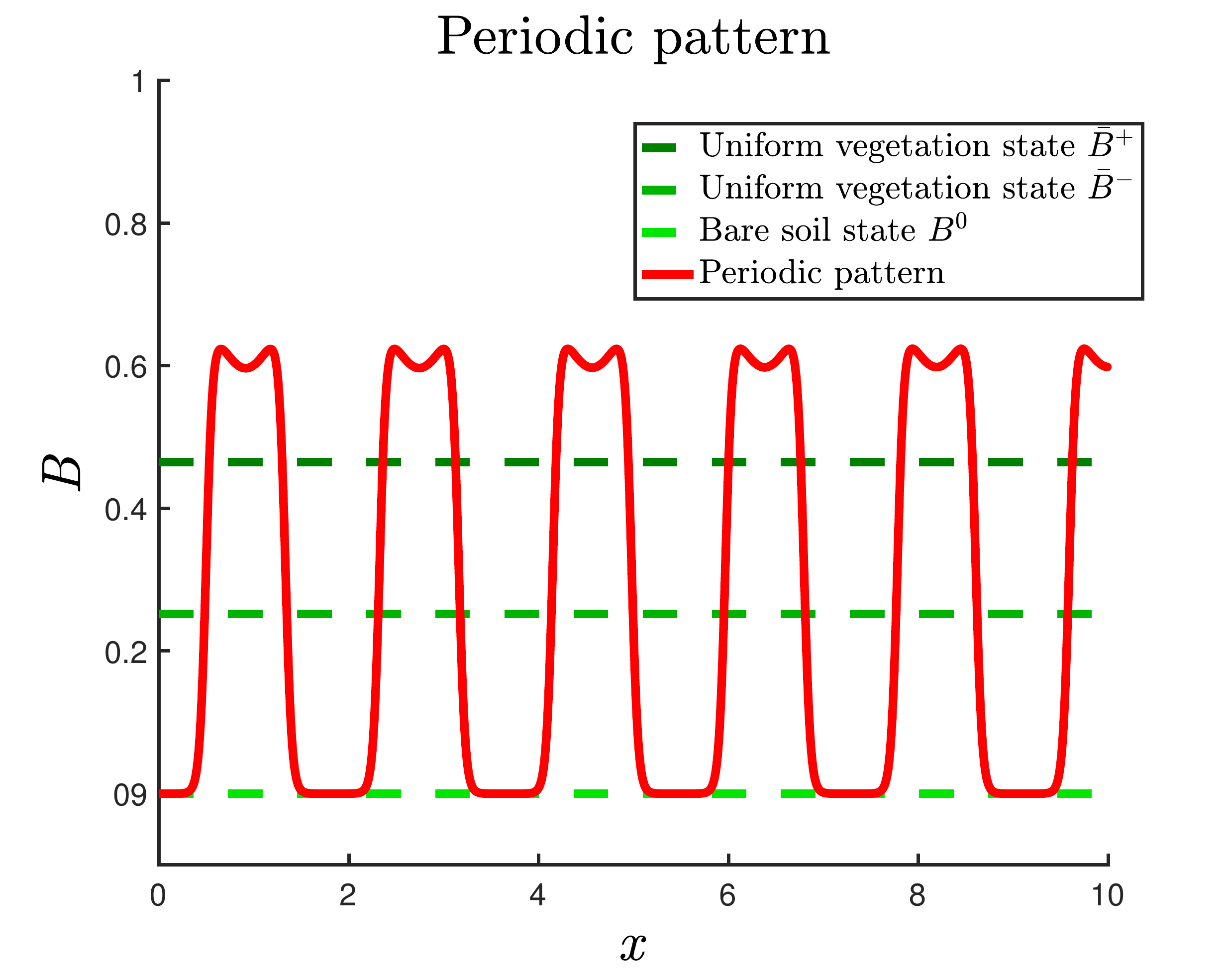}
		\caption{\label{fig:numerics_periodic_front}}
	\end{subfigure}%
	\hspace{0.4cm}
	\begin{subfigure}[b]{0.48\textwidth}
		\centering
		\includegraphics[width=\textwidth]{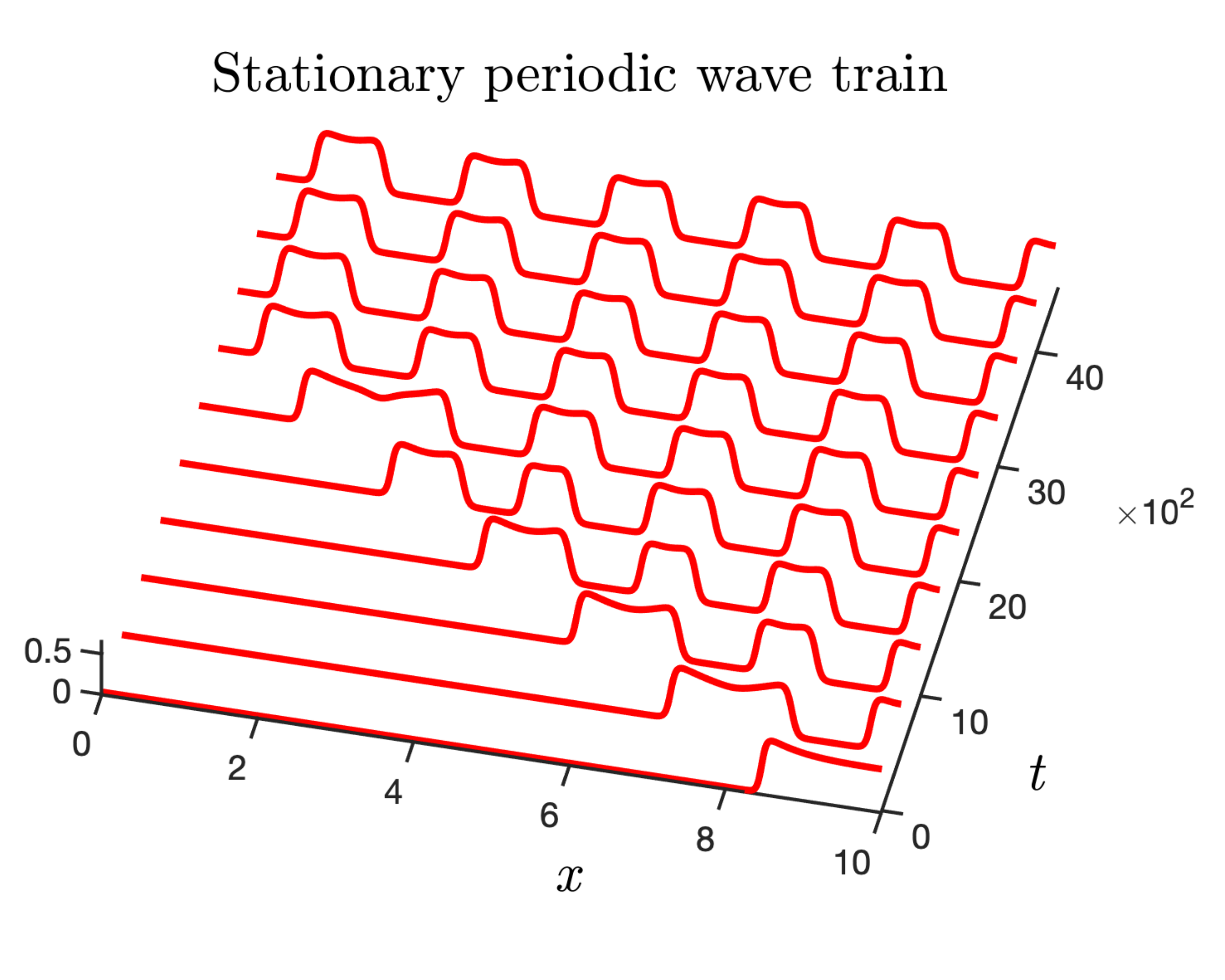}
		\caption{\label{fig:numerics_periodic_front_3D}}
	\end{subfigure}%
	\caption{(a) The standing asymptotic (for $t \to \infty$)  spatially periodic spot/gap pattern generated by the invasion dynamics of Fig. \ref{fig:numerics_periodic_front_3D}. (b) A time/space plot of a vegetation front that invades the bare soil state and leaves the spatially periodic pattern of Fig. \ref{fig:numerics_periodic_front} behind. In both plots: $ a =0.032, \Psi =  1.619 , \Phi = 0.3, \Omega=0.1, \Theta = 0.5, \varepsilon = \sqrt{0.01}$.}
	\label{fig:numerics_periodic_multifront}
\end{figure}

As demonstrated in section \ref{sec:locfront}, heteroclinic 1-front orbits can occur both as traveling -- Theorem \ref{T:primaryfronts} -- or stationary patterns -- Theorem \ref{T:stand1fronts}. We confirm this numerically in Figs. \ref{fig:numerics_1_front_standing} and \ref{fig:numerics_1_front_traveling}. Note that these fronts may either represent the retreat of vegetation by the invasion of the bare soil state into the homogeneous vegetation state -- $c > 0$ -- or the expansion of a (homogeneously) vegetated area into the bare soil state -- $c < 0$. In fact, in order to find the stationary 1-front, we need to tune a single parameter -- $\Phi$ in the statement of Theorem \ref{T:stand1fronts} and $\Psi$ in Fig. \ref{fig:numerics_standing_and_traveling_waves} -- to the border point between the ranges of left-traveling and right-traveling 1-fronts (the `Maxwell point' \cite{BelEtAl12} described by the co-dimension 1 set $\R_{\rm s-1f}$ in Theorem \ref{T:stand1fronts}).
\\ \\
The existence of the homoclinic stationary 2-front pattern depicted in Fig. \ref{fig:numerics_2_front_standing} was established in Theorem \ref{T:standing2fspots}. Note that the level of vegetation on the plateau that determines the spot remains
relatively far away from the value $\bar{B}^+$ of the uniform vegetation state $(\bar{B}^+,\bar{W}^+)$. This is caused by the fact that the homoclinic orbit associated to the spot pattern follows an orbit on the slow manifold $\mathcal{M}^+_\eps$ of   (\ref{eq:fast}) that does not approach the critical point associated to $(\bar{B}^+,\bar{W}^+)$ on $\mathcal{M}^+_\eps$ -- see the sketch of the skeleton structure in Fig. \ref{fig:3skeletons}. From the ecological point of view, a vegetation spot benefits from soil water diffusion from the adjacent water-rich bare soil areas -- see the $W$-profile in Fig. ~\ref{fig:numerics_1_front_ehud} -- besides direct rainfall, and therefore has higher biomass density as compared with uniform vegetation. This also explains (in ecological terms) why the biomass density at the edge of a front, a spot or a gap is higher -- a property also exhibited by `fairy circles' \cite{Wiegand2014,Zelnik2015}. See also the upcoming discussion below of 2-front vegetation gaps (Fig. \ref{fig:numerics_fairy_circle}). In Fig. \ref{fig:numerics_2_front_traveling}, we show the dynamics of the 2 interacting fronts of an evolving 2-front pattern: the distance between the fronts slowly increases while it settles into a stationary standing spot -- see Remark \ref{rem:travpatts} and the discussion in section \ref{sec:disc}.
\\ \\
One of the original motivations to analyze far-from-equilibrium patterns in the Gilad et al. model, was to gain a fundamental understanding of `fairy circles' -- a somewhat subtle phenomenon (for instance) observed in western Namibia \cite{Wiegand2014,Zelnik2015}. The homoclinic stationary 2-front gap patterns of (\ref{eq:scaled}) established by Theorem \ref{T:standing2fgaps} and shown in Fig. \ref{fig:numerics_fairy_circle} indeed have the strongly localized nature of observed fairy circles. Moreover, the spot/gap patterns of Theorem \ref{T:periodics} represent the observed (nearly) periodic arrays of fairy circles (see Fig. \ref{fig:numerics_periodic_front} and notice that the ratio between the lengths of the vegetated state and the bare soil patches typically varies from 0 to $\infty$ in this family (section \ref{sec:perpatt})).  As noted, fairy circle gap patterns have an excess of vegetation at the edge of the gap as distinctive feature -- see for instance the images in \cite{Wiegand2014,Zelnik2015}. In mathematical terms, this means that the connecting fronts are non-monotonous. In the context of the present model, this non-monotonicity is caused by the orientation and curvature of the slow manifold $\mathcal{M}^+_\eps$ relative to the path traced by the gap pattern over $\mathcal{M}^+_\eps$ -- see the projection in Fig. \ref{fig:projection_pulse_2front_manifolds}a for a representation of this `geometrical mechanism' for spot patterns. We refer to   \cite{Fernandez-Oto2014phil_tran_A} for a 1-component model in which the non-monotonicity of the fronts originates from nonlocal effects.
\\

\begin{figure}[h!]
	\centering
	\begin{subfigure}[b]{0.48\textwidth}
		\centering
		\includegraphics[width= \textwidth]{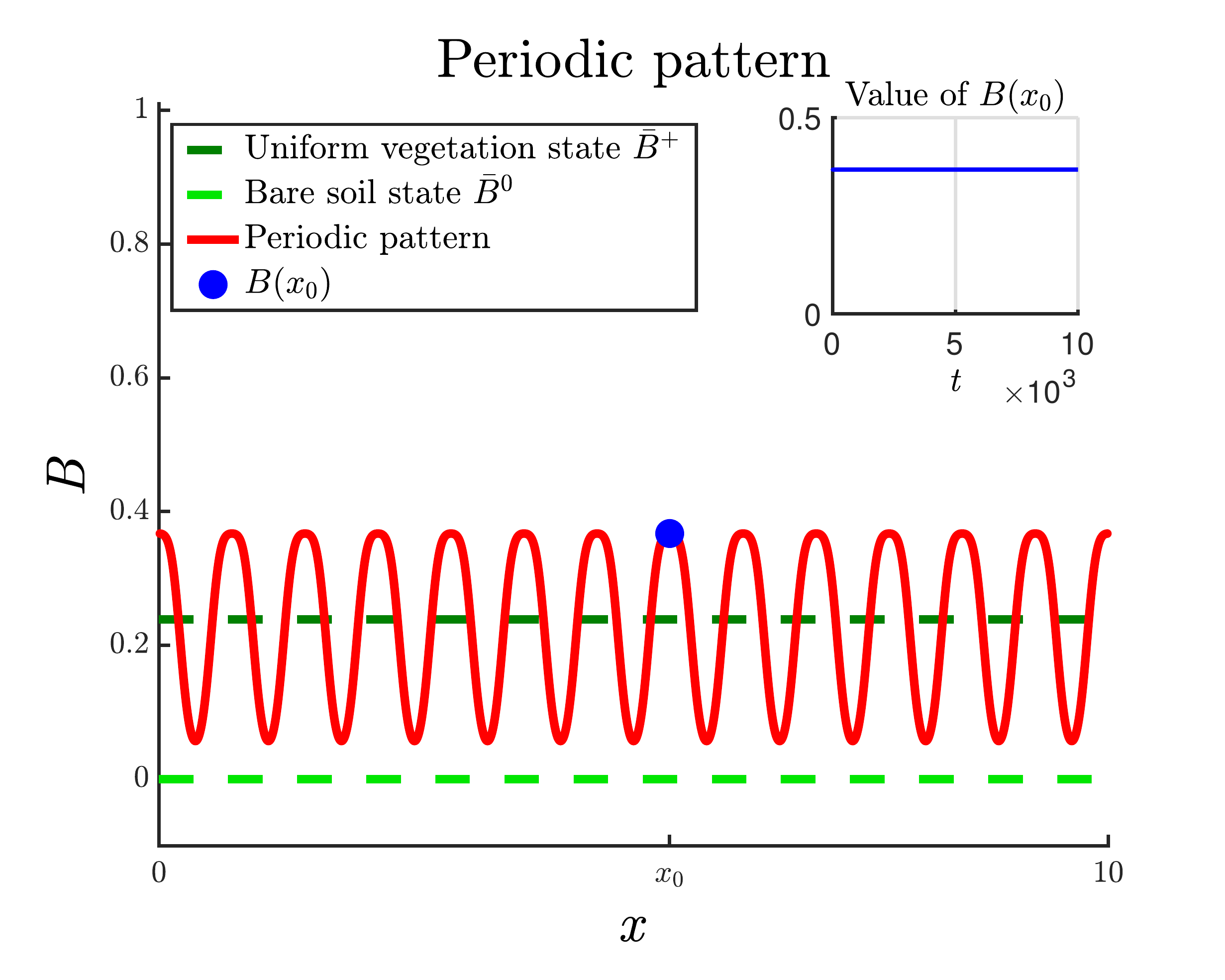}
		\caption{\label{fig:numerics_turing} $ a =0.25, \Psi = 0.42392; ,
			\Phi = 0.059, \Omega=0.4, \Theta = 0.5, \varepsilon = \sqrt{0.2}$.}
	\end{subfigure}%
	\hspace{0.4cm}
	\begin{subfigure}[b]{0.48\textwidth}
		\centering
		\includegraphics[width= \textwidth]{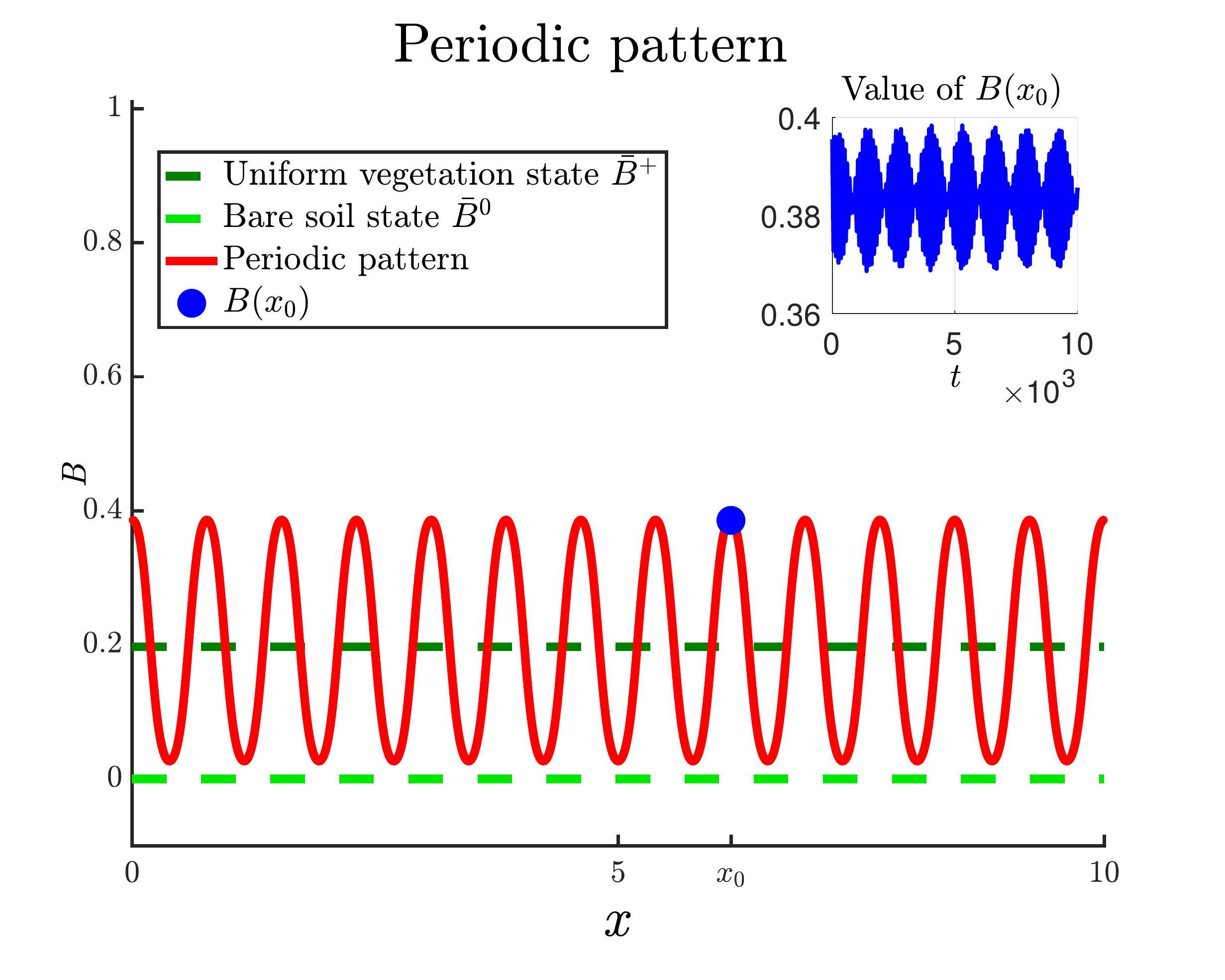}
		\caption{\label{fig:numerics_turing_hopf} $ a =0.25, \Psi =  0.38539 ,
 			\Phi = 0.059, \Omega=0.4, \Theta = 0.5, \varepsilon = \sqrt{0.1}$.}
	\end{subfigure}%
	\caption{{\bf (a)} A small amplitude stationary spatially periodic solution generated by a Turing bifurcation (see Remark \ref{rem:Turing}. {\bf (b)} A pattern that is periodic in space and time that appeared by decreasing $\Psi$ from the Turing pattern of Fig. \ref{fig:numerics_turing} through a Hopf bifurcation (in time).}
	\label{fig:numerics_turing_and_hopf}

\end{figure}

In section \ref{sec:stand2fronts} -- and especially in Remark \ref{rem:2ftofasthom} -- we discussed the bifurcation of the homoclinic slow-fast-slow-fast-slow 2-front spot pattern of Theorem \ref{T:standing2fspots} into a homoclinic slow-fast-slow pulse pattern as it `detaches' from $\M^+_\eps$. Such a (numerically stable) `detached' spot pattern of pulse type is shown Fig. \ref{fig:numerics_1_pulse}. In Fig. \ref{fig:projection_pulse_2front_manifolds}, the detachment process is shown by projections of both the 2-front spot pattern of Fig. \ref{fig:numerics_2_front_standing} and the pulse pattern of Fig. \ref{fig:numerics_1_pulse} on the $(w,b)$-plane: as the parameter $\Psi$ -- which is linearly related to the rainfall $P$ in the original model (\ref{def:pars}) -- is decreased below a critical value $\Psi_* = 1.2952$ the vegetated plateau disappears and the 2-front spot solution transforms into a 1-pulse solution. Note that this 1-pulse solution is of the `classical' Klausmeier-Gray-Scott (and/or Gierer-Meinhardt \cite{DGK01}) type already discussed in the introduction of section \ref{sec:stand2fronts}: its existence may be established by the methods of \cite{DGK01,DVeerman2015} and the references therein.
We also found that the spatially periodic spot/gap patterns of Theorem \ref{T:periodics} may have quite a large domain of attraction: Fig. \ref{fig:numerics_periodic_front_3D} shows the evolution of a traveling vegetation front into the bare soil state that leaves behind a spatially periodic spot/gap pattern -- Fig. \ref{fig:numerics_periodic_front}. This behavior may possibly be related to the existence of a Turing bifurcation -- see Remark \ref{rem:Turing} -- of the uniform vegetation state and calls for further studies. Finally, we show in Fig. \ref{fig:numerics_turing} such a numerically obtained, almost sinusoidal, small amplitude Turing pattern that bifurcated from a (destabilized) uniform vegetation state and note that there are paths through parameter space on which this Turing pattern evolves into a (periodic) multi-pulse pattern -- built from homoclinic 1-pulses of Gierer-Meinhardt/Gray-Scott type as in Fig. \ref{fig:numerics_1_pulse} and typically observed in Klausmeier-type models \cite{SewaltD17,Siteur2014eco_comp,vdSteltetal13} -- that subsequently touches down on $\M^+_\eps$ like the solitary pulses of Figs. \ref{fig:numerics_1_pulse} and \ref{fig:numerics_2_front_standing}, to indeed take the shape of the periodic fairy circle-type spot/gap pattern of Theorem \ref{T:periodics} and Fig. \ref{fig:numerics_periodic_front}. By further tuning parameters it may also happen that the stationary, spatially periodic, Turing pattern undergoes a Hopf bifurcation (in time), resulting in an oscillating pattern that is periodic both in space and in time -- see Fig. \ref{fig:numerics_turing_hopf}. (Note, however, that it is not clear whether this may occur for ecologically feasible parameters -- see \cite{Tzuk2019}.)


\newpage
\begin{remark}
\label{rem:parsFig1}
The basic configurations shown in Fig. \ref{fig:basicpatterns} have been obtained by the same procedures as used in this section, with the following parameter settings. The traveling 1-front: $ a =0.0008, \Psi =   1.6248, \Phi = 0.3, \Omega=0.1, \Theta = 0.2, \varepsilon = \sqrt{0.005}$; the stationary homoclinic 2-front vegetation spot: $ a =0.032, \Psi =  1.3714, \Phi = 0.3, \Omega=0.1, \Theta = 0.2, \varepsilon = \sqrt{0.005}$; the stationary homoclinic 2-front vegetation gap: $ a =0.032, \Psi =  1.2762, \Phi = 0.3, \Omega=0.1, \Theta = 0.2, \varepsilon = \sqrt{0.005}$; the stationary spatially periodic multi-front:  $a =0.032, \Psi =  1.619, \Phi = 0.3, \Omega=0.1, \Theta = 0.5, \varepsilon = \sqrt{0.01}$.
\end{remark}

\subsection{Discussion}
\label{sec:disc}

Of course, the potential relevance of the various singular slow/fast patterns constructed in this paper is ultimately determined by their stability as solutions of PDE (\ref{eq:meron})/(\ref{eq:scaled}). In general, this is a seriously hard problem to study analytically. However, the singularly perturbed nature of the patterns considered here enables us to explicitly and rigorously analyze the (spectral) stability of the constructed (multi-)front patterns. In fact, the explicit `control' we established on the slow-fast structure of the (multi-)fronts provides the perfect (and necessary) starting point for a spectral stability analysis along the lines of (for instance) \cite{DIN2004,CdRS2016} and \cite{dRDR2016} (for the spatially periodic patterns). This is especially the case for the basic front/spot/gap/periodic patterns of Theorems \ref{T:primaryfronts}, \ref{T:stand1fronts}, \ref{T:standing2fspots}, \ref{T:standing2fgaps}, \ref{T:periodics} shown in Fig. \ref{fig:basicpatterns}.
\\ \\
However, the question whether the non-basic, `higher order' patterns (for instance) sketched in Fig. \ref{fig:sketchessection1} can be stable also requires novel mathematical insights and methods. All constructed higher order patterns involve the existence of persisting periodic solutions on the slow manifold $\M^+_\eps$ -- see Theorems \ref{T:frontstoperiodic}, \ref{T:count1fronts}, \ref{T:statfrontstoperiodic}, \ref{T:higherorder2fronts} and Corollary \ref{cor:higherorderperiodics}. Therefore, the structure of the spectrum associated to the stability of the higher order patterns essentially depends on the preliminary question about the spectrum and stability of the persistent periodic solutions on the slow manifold of Theorem \ref{T:slowflow} -- and especially of their homoclinic (or heteroclinic) limits also considered in Theorem \ref{T:slowflow}. In fact, this issue is not (at all) specific for the explicit model here. We claim that higher order patterns of the type sketched in Fig. \ref{fig:sketchessection1} will generically appear as potentially stable solutions in a fully general class of singularly perturbed reaction-diffusion models that includes (\ref{eq:scaled}),
\begin{equation}
\label{eq:gen}
\left\{	
\begin{array}{lcrcl}
U_t &=& U_{xx} & + & F(U,V),\\
V_t &=& \frac{1}{\varepsilon^2}V_{xx} & + & G(U,V).
\end{array}
\right.
\end{equation}
By going into a traveling framework -- and thus introducing $\xi = x - ct$, $U(x,t) = u(\xi)$, $V(x,t) = v(\xi)$, $p(\xi) = u_\xi(\xi)$, $q = v_\xi(\xi)/\eps$ as in section \ref{sec:Intro} -- (\ref{eq:gen}) is reduced into the 4-dimensional form of (\ref{eq:fast}). By taking the $\eps \to 0$ limit, we find that the 2-dimensional (reduced) slow manifolds are determined by  $F(v_0, u) = 0$ (and $p =0$, $(v_0,q_0) \in \mathbb{R}^2$) -- see section \ref{sec:setup} -- which generically determines $J \geq 1$ branches, locally given by graphs,
\[
\M^j_0 = \{(u,p,v,q) \in \mathbb{R}^4: u=f_j(v), p=0\}, \; j = 1,2, ..., J,
\]
with $f_j(v)$ such that $F(f_j(v),v) \equiv 0$ (cf. (\ref{def:M00Mpm0}) and note that $J = 3$ for (\ref{eq:scaled})). For those (parts of) $\M^j_0$ that are normally hyperbolic, $\M^j_0$ persists as $\M^j_\eps$, that is approximately given by,
\[
\M^j_\eps = \{(u,p,v,q) \in \mathbb{R}^4: u=f_j(v) + \eps cq u_{1}^j(v)  + \OO(\eps^2), p= \eps q p_{1}^j(v) + \OO(\eps)^2\}
\]
with
\[
u_1^j(v) = - f_j'(v)/\frac{\partial F}{\partial u}(f_j(v),v), \; p_1^j(v) = f_j'(v).
\]
(cf. (\ref{M+eps}), (\ref{p1b1})). Thus, completely analogous to the analysis in section \ref{sec:slowfull}, we find that the slow flow on a persisting, normally hyperbolic 2-dimensional slow manifold $\M^j_\eps$ is given by a planar Hamiltonian system perturbed by a nonlinear friction term,
\[
v_{XX} + G(f_j(v),v) + \eps c q \left[1- \frac{\partial G}{\partial u}(f_j(v),v) u_1^j(v)\right] = \OO(\eps^2),
\]
with $X = \eps \xi$ (cf. (\ref{eq:slow+eps})). Typically, the unperturbed $\eps \to 0$ limit $v_{XX} + G(f_j(v),v)$ -- i.e. the reduced slow flow on $\M^j_0$ -- is nonlinear and has families of periodic solutions and homoclinic or heteroclinc orbits to critical points on $\M^j_\eps$ that correspond to (potentially stable \cite{Dreview}) homogeneous background states $(U(x,t),V(x,t) \equiv (\bar{U},\bar{V})$ of PDE (\ref{eq:gen}) -- as is the case for (\ref{eq:SRSplus}) on $\M^+_0$. Thus, indeed, the situation is completely similar to that of section \ref{sec:slowfull}: using Melnikov-type arguments persistence results equivalent to Theorem \ref{T:slowflow} may be deduced, also in the present general setting. The geometric framework of orbits `jumping up and down' between two (normally hyperbolic) slow manifolds $\M^j_\eps$ and $\M^k_\eps$ presented in section \ref{sec:touchdown} is based on the persistence of both the stable and unstable manifolds $W^{s,u}(\M^{j,k}_\eps)$ of $\M^{j,k}_\eps$ and thus of the intersections $W^{u}(\M^{j}_\eps) \cap W^{s}(\M^{k}_\eps)$ and $W^{s}(\M^{j}_\eps) \cap W^{k}(\M^{k}_\eps)$. Therefore, we may use the arguments, methods and insights of section \ref{sec:locfront} to deduce the equivalents of the `higher order' existence Theorems \ref{T:frontstoperiodic}, \ref{T:count1fronts}, \ref{T:statfrontstoperiodic}, \ref{T:higherorder2fronts} and Corollary \ref{cor:higherorderperiodics} in the setting of general system (\ref{eq:gen}). Moreover, this also implies that bifurcation scenarios as sketched in Fig. \ref{fig:cPhi2} appear generically (where we notice that the sketch in Fig. \ref{fig:cPhi2} was just a first example -- many other scenarios may occur). In fact, the geometrical setting allows us to (for instance) explicitly establish the existence of heteroclinic networks of orbits jumping between various slow manifolds $\M^j_\eps$ and (slowly) following periodic orbits on $\M^j_\eps$ in between its fast jumps -- like the networks considered in \cite{Radem2005,Radem2010} and the references therein. Thus, the above noted preliminary (and essential) issue of the spectrum associated to the stability of the persisting periodic and homoclinic solutions on $\M^+_\eps$ of Theorem \ref{T:slowflow} also has a fully general -- and thus fundamental -- counterpart, with a similar relevance for the higher order patterns (almost) heteroclinic to these orbits. In other words, insight in the spectrum associated to the stability of the orbits on $\M^j_\eps$ established by a generalization of Theorem \ref{T:slowflow} for (\ref{eq:gen}) is expected to yield explicit insight in the stability and bifurcations of the higher order patterns in PDE (\ref{eq:gen}) established by the generalizations of Theorems \ref{T:frontstoperiodic}, \ref{T:count1fronts}, \ref{T:statfrontstoperiodic}, \ref{T:higherorder2fronts}, etc. and the subsequent more complex `networks'.
\\ \\
This will be the subject of future work, both in the setting of explicit system (\ref{eq:scaled}) -- which will also include a systematic numerical search for the higher order patterns sketched in Fig. \ref{fig:sketchessection1} -- and in the general setting of (\ref{eq:gen}).
\\ \\
To optimally embed the present analysis in the ecological context, we first need to obtain insight in the ranges of the (scaled) parameters $(a, \Psi, \Phi, \Omega, \Theta)$ of (\ref{eq:scaled}) that may correspond to ecologically realistic settings of the unscaled model . Finding such parameter ranges is possible (more than in the Klausmeier model) since (\ref{eq:meron}) is directly linked to the more elaborate 3-component model of Gilad et al. \cite{GiladEtAl2004} -- \cite{Zelnik2015}, \ref{Appendix:derivation} -- and thus to concrete underlying ecological mechanisms \cite{GiladEtAl2007,Meron2015,ShefferEtAl2013}. A crucial question for the potential ecological relevance of the above discussed higher order patterns is whether there are realistic values of $(\Lambda, K, E, M, P, N, R, \Gamma)$ for which there are 2 critical points on $\M^+_\eps$, i.e. for which $\C^2 - 4 \A \D > 0, \C < 0, \D > 0$ (section \ref{sec:crithom}) -- where $(\A,\C,\D)$ is related to $(\Lambda, K, E, M, P, N, R, \Gamma)$ in a rather nonlinear fashion by (\ref{scalingpars}), (\ref{eq:scaled}), (\ref{def:pars}) (\ref{defABC}), (\ref{defWD}). Naturally, this will be part of our upcoming work on (\ref{eq:scaled}).
\\ \\
Each of the higher order invasion fronts established by Theorem \ref{T:count1fronts} and sketched in Figs. \ref{fig:sketchessection1}a and \ref{fig:sketchessection1}b travels with a different speed -- in fact, the (discrete) family may even limit on a stationary front pattern (Remark \ref{rem:counthetsc=0}). Thus, when stable, these invasion fronts may introduce the possibility of slowing down gradual desertification. Moreover, stationary multi-front patterns may bifurcate into traveling patterns with the same structure -- see Remark \ref{rem:travpatts} for a brief sketch of the underlying geometrical mechanism. When stable, the appearance of such traveling multi-front patterns -- either localized spots or spatially periodic wave trains -- may have a similar ecological interpretation and relevance: localized vegetated states may even reverse desertification by invading bare soil -- see \cite{ZelnikEtAl2018,Zelnik2018eco_ind,ZelnikEtAl2017}. Together, the various traveling 1-front patterns and traveling multi-spots form an interacting group of invasion patterns within the transition zone between the bare soil state and a homogeneous vegetation state; in principle all entities in this group travel with different speed. Understanding pattern formation in this zone -- and especially also understanding the translation and/or expansion of this zone in terms of the parameters in the model -- may have a direct ecological significance. In mathematical terms, such a study may also be performed by a front interaction analysis along the lines of \cite{vHeijsterDKP10,CBDvHR15,CBvHRI19} -- although the dynamics generated by (\ref{eq:scaled}) in this ecological transition zone is expected to be richer than that of the generalized FitzHugh-Nagumo model considered there.
\\ \\
To truly obtain ecological relevance, we must consider the model in 2 space dimensions.  Clearly, the extension to more than 1 space dimension does pose fundamental challenges, moreover 2-dimensional systems show much richer dynamical behaviors associated with propagating fronts -- see for instance \cite{HM1994,HM1997}. However, the results obtained here form a foundation upon which aspects of the step from 1 to 2 space dimensions can be taken -- see for instance \cite{vHeijsterSandstede2014,SieroEtAl2015,Ward2018,WeiWinter2014} and the references therein. By extending the patterns constructed here trivially in the second spatial direction, the above mentioned stability analysis can be directly extended to include the stability (and bifurcations) of planar (multi-)fronts/interfaces (where we for simplicity neglect the (technical) fact that (\ref{eq:scaled}) takes a somewhat different form in $\mathbb{R}^2$, \cite{Meron2015,Zelnik2015}, \ref{Appendix:derivation}). Unlike in the extended Klausmeier model \cite{SewaltD17,SieroEtAl2015}, localized stripes are of 2-front type and may thus be expected to possibly be stable -- see \cite{BCD2019} for a rigorous treatment in a generalized Klausmeier type model (posed on a sloped terrain without a diffusion term for the water component -- like the original Klausmeier model \cite{Klausmeier1999}). Naturally, the interfaces will evolve and their curvature driven dynamics may be studied analytically along the lines of \cite{MW2017}. Especially in the above discussed multi-front transition region between bare soil and homogeneous vegetation, the ecosystem dynamics generated by the model may be very rich and complex -- see for instance \cite{HM1994,HM1997,HM1998}.
\\ \\
As a final direction of possible future research, we note that our results may be used to establish the existence -- and later stability -- of localized patterns in the original 3-component model of Gilad et al. \cite{GiladEtAl2004}. Since (\ref{eq:meron}) and thus (\ref{eq:scaled}) -- is obtained from the nonlocal, 3-component model of Gilad et al. (see (\ref{eq:gilad})) in a systematic way -- i.e. by taking several limits (\cite{Meron2015,Zelnik2015}, \ref{Appendix:derivation}) -- it may be expected that it is possible to establish the persistence of patterns constructed here into the nonlocal, 3-component setting, especially since these patterns are constructed geometrically through transversal intersections of invariant manifolds. Once again, this is interesting and relevant both from mathematical and ecological point of view: (asymptotically) small nonlocal and topographical terms may have a significant effects, even on the most simple -- `basic' -- (vegetation) patterns exhibited by a model \cite{BCBD2020,Escaff2015pre}.
\\ \\
{\bf Acknowledgement.} OJ acknowledges the (partial) support of NWO through its {\it Complexity} program and EM acknowledges the support of the Israel Science Foundation under grant number 1053/17.

\bibliographystyle{abbrv}

\bibliography{Newbibnourl}

\addcontentsline{toc}{section}{References}

\appendix

\section{Derivation of the model equations in one spatial dimension}
\label{Appendix:derivation}
We follow \cite{Zelnik2015} to briefly show how model (\ref{eq:meron}) is derived from the original model introduced in \cite{GiladEtAl2007} and given by
\begin{eqnarray}
\label{eq:gilad}
\left\{	
\begin{aligned}
\partial_T B &= G_BB(1-B/K) - MB + D_B \nabla^2B, \\
\partial_T W &= IH - N(1-RB/K)W - G_WW + D_W \nabla^2 W,\\
\partial_T H &= P - IH 	+ D_H\nabla^2(H^2) + 2D_H \nabla H\cdot \nabla Z + 2D_HH\nabla^2 Z,
\end{aligned}
\right.
\end{eqnarray}
where
\begin{subequations}
  \begin{align}
      G_B(\textbf{X},T) & = \Lambda \int_\Omega G(\textbf{X},\textbf{X}^\prime,T) W(\textbf{X}^\prime,T) \textbf{dX}^\prime \\
      G_W(\textbf{X},T) & = \Gamma \int_\Omega G(\textbf{X}^\prime,\textbf{X},T) B(\textbf{X}^\prime,T) \textbf{dX}^\prime \\
      G(\textbf{X},\textbf{X}^\prime,T) & = \left(\frac{1}{\sqrt{2 \pi S_0^2}}\right)^2 \exp \left[ -\frac{|\textbf{X}-\textbf{X}^\prime|^2}{2S_0^2(1 + EB(\textbf{X},T))^2} \right] \label{eq:G} \\
      I & = A \frac{B(\textbf{X},T)+Qf}{B(\textbf{X},T)+Q}
  \end{align}
\end{subequations}
with $\textbf{X} = (X,Y)$ the spatial coordinates of the $2$-dimensional system. The last equation in (\ref{eq:gilad}) describes overland water flow with $H$ being the height of a thin layer of surface water above ground level given by the topography function $Z$. We consider the case of a flat terrain , for which $Z$ = constant, and of high infiltration rates $I$,  both in bare soil and vegetated areas (no infiltration contrast), for which $I$ can be assumed to be independent of $B$. Both conditions are met in the Namibian fairy-circles ecosystems that consist of sandy soil. Since $H$ varies on time scales much shorter than those of $W$ and $B$, these conditions imply fast equilibration of surface water at $H=P/I$. Insertion of this equilibrium value in the equation for $W$ results in the elimination of the surface water equation.
\\ \\
A further simplification we make is related to the nonlocal forms of the biomass growth rate, $G_B$, and the water uptake rate, $G_W$, in (\ref{eq:gilad}). We assume, consistently with the plant species in the Namibian fairy-circles ecosystems, that the roots, described by the root kernel  $G(\textbf{X},\textbf{X}^\prime,T)$,  are laterally confined. We employ this assumption by taking the lateral root extension of a seedling, $S_0$, to be very small. Using the limit $S_0\to 0$ in the integrals in (\ref{eq:G}) we obtain, for a $1$-dimensional system, the simpler algebraic expressions
\begin{align}
G_B(X,T)
& = \Lambda\int_{\Omega}\lim_{S_0\rightarrow 0}  \frac{1}{\sqrt{2\pi S_0^2}}\exp\left[-\frac{|{X} - {X'}|^2}{2S_0^2(1 + EB(X,T))^2}\right]B(X',T)dX'\nonumber\\
&= \Lambda\left(1+ EB(X,T)\right)B(X,T).
\end{align}{
Similarly,
\begin{align}
G_W(X,T)
&= \Gamma\left(1+ EB(X,T)\right)W(X,T).
\end{align}
Inserting these expressions in (\ref{eq:gilad}) we obtain the 2-component model (\ref{eq:meron}). Finally, we note that in \cite{Zelnik2015} this reduction was performed in 2 space dimensions and that the general $n$-dimensional situation is considered in \cite{Meron2015}.


\section{The derivation of the scaled model}
\label{Appendix:scaling}

Introducing the scalings (\ref{eq:biomasswater}) into (\ref{eq:meron}) yields,
\begin{eqnarray}
\label{eq:equation_first_scaling}
\left\{	
\begin{aligned}
\alpha \delta B_t &= \alpha \beta \Lambda WB(1-\alpha B/K)(1+\alpha E B)- \alpha MB + \alpha \gamma^2 D_B B_{xx}, \\
\beta \delta W_t &= P - N\beta W(1-\alpha R B/K) - \alpha \beta\Gamma WB(1+\alpha EB) + \beta \gamma^2 D_W W_{xx},
\end{aligned}
\right.
\end{eqnarray}
which can be brought into the form,
\begin{eqnarray}
\label{eq:equation_second_scaling}
\left\{	
\begin{aligned}
\frac{\delta K}{\alpha^2 \beta \Lambda E}  B_t &= \left(\frac{K}{\alpha^2 E}W- \frac{MK}{\alpha^2 \beta \Lambda E}\right)B + \frac{K}{\alpha  E} \left( E - \frac{1}{K}\right)WB^2 - WB^3 + \frac{\gamma^2K}{\alpha^2 \beta \Lambda E} D_B B_{xx} \\
\frac{ \delta K}{\alpha^2 \beta \Lambda E} W_t &= \frac{PK}{\alpha^2 \beta^2 \Lambda E} -\frac{K}{\alpha^2 \beta \Lambda E}  \left[N\left(1-\frac{\alpha R}{K}B \right)  + \alpha \Gamma B(1+\alpha EB)\right] W + \frac{\gamma^2K}{\alpha^2 \beta \Lambda E} D_W W_{xx}.
\end{aligned}
\right.
\end{eqnarray}
By choosing $\delta$ and $\gamma$ as in (\ref{scalingpars}), we arrive at,
\begin{eqnarray}
\label{eq:progress}
\left\{	
\begin{aligned}
B_t &= \left(\frac{K}{\alpha^2 E}W- \frac{MK}{\alpha^2 \beta \Lambda E}\right)B + \frac{1}{\alpha} \left(K - \frac{1}{E}\right)WB^2 - WB^3 +  B_{xx}, \\
W_t &= \frac{PK}{\alpha^2 \beta^2 \Lambda E} -\frac{K}{\alpha^2 \beta \Lambda E}  \left[N\left(1-\frac{\alpha R}{K} B\right)  + \alpha \Gamma B(1+\alpha EB)\right] W + \frac{D_W}{D_B}W_{xx}.
\end{aligned}
\right.
\end{eqnarray}
Next, we use our freedom in $\alpha$ and $\beta$ to simplify the $B$-equation and scale the factors of the $B$- and $WB^2$-terms (to $-1$ and to $+1$ respectively) -- which is achieved by the choices in (\ref{scalingpars}),
\begin{eqnarray}
\label{eq:final_form}
\left\{	
\begin{aligned}
B_t &= \left(\frac{KE}{(KE - 1)^2}W- 1\right)B + WB^2 - WB^3 +  B_{xx}, \\
W_t &= \frac{\alpha^2P\Lambda E}{M^2K } - \left[N\left(1-\frac{\alpha R}{K} B\right)  + \alpha \Gamma B(1+\alpha EB)\right] \frac{W}{M} + \frac{D_W}{D_B}W_{xx}.
\end{aligned}
\right.
\end{eqnarray}
This is equivalent to (\ref{eq:scaled}) by definitions (\ref{def:aeps}) and (\ref{def:pars}).
\\ \\
Note that our choice to scale the factor of the term $WB^2$ in the $B$-equation to $+1$ implies -- together with the (implicit, natural) assumption that $\tB$ and $B$ have the same signs (\ref{eq:biomasswater}) -- that we have chosen to consider $EK > 1$. Of course, it may happen that $0 < EK \leq 1$. In these cases, either the term $WB^2$ disappears from the equation -- in the `non-generic' case $EK = 1$ -- or its pre-factor can be scaled to $-1$. All of the analysis in this work can also be performed for $EK \leq 1$, without any conceptual differences. However, we chose to focus of $EK > 1$ -- and thus on a $+WB^2$ term in (\ref{eq:scaled}) -- to not further complicate the necessary `algebra'.

\section{\bf Lemma \ref{lem:posigmasmal} and the Bogdanov-Takens bifurcation scenario}
\label{Appendix:BogdanovTakens}
A planar ODE of the form
\[
\left\{
\begin{aligned}
y_{\tilde{X}} & = z\, ,\\
z_{\tilde{X}} & =  \beta_1 + \beta_2 y +  y^2 + s yz + \mathcal{G}(y,z). \, ,
\end{aligned}
\right.
\]
for $ \beta_1, \beta_2 \in \mathbb{R}, s = + 1 $  is known to possess two fixed points - a saddle and a focus - with an unstable periodic orbit (that emerged from the focus in a Hopf bifurcation in the open parameter region)
\[
\mathcal{S}_{BT} = \left\{ (\beta_1, \beta_2) ~|~ \beta_2 <  0 \, , \beta_1 < -\frac{6}{25} \beta_2^2 + o(\beta_2^2) \right\} \, .
\]
The right border $ \{ \beta_2 < 0, \beta_1 = 0 \} $ marks the Hopf bifurcation, while the left border $ \{ \beta_2 <  0 \, , \beta_1 = -\frac{6}{25} \beta_2^2 + o(\beta_2^2)  \} $ describes the region where a homoclinic orbit emerged from the periodic orbit (whose period tends to infinity towards that border).\\
\\
Here, we denote the slow system \eqref{eq:slow+eps} by
\[
\left\{
\begin{aligned}
w_{X} &= q\, ,\\
q_{X} & = F(w) + \eps c q \rho_1(w) + G(w,q)
\end{aligned}
\right.
\]
where $F(w) = - \mathcal{A} + \left(\mathcal{B}+ a\Theta \right) w + \mathcal{C} w \sqrt{a + \frac14 - \frac{1}{w}}$
and $ G $ accounts for the higher order term, and assume that the parameters are chosen such that both fixed points are on $\mathcal{M}_{\varepsilon}^+$, $ \mathcal{D} $ is close (but  beyond) the saddle-node bifurcation point, that is, $ \mathcal{D} = \frac{\C^2}{4\A} - \sigma^2 \A \, , 0 < \sigma \ll 1 $ and the $w$-coordinate of both emerging fixed points is well within the strip $(4/(1+4a), 1/a)$. For $ \sigma $ sufficiently small, there is a neighborhood of $ w_0^{SN} $ such that the slow ODE has the form
\begin{equation}
\label{eq:slow_expanded}
\left\{
\begin{aligned}
w_{X} &= q\, ,\\
q_{X} & = \mu_1+ \mu_2 w + \mu_3 w^2 + \delta q + \mu_4 w q + \widetilde{G}(w,q)
\end{aligned}
\right.
\end{equation}
where
\[
\mu_1= F(w_0^{SN}) \, , \quad
\mu_2 = F'(w_0^{SN})\, , \quad
\mu_3 = \frac12 F''(w_0^{SN}) \, , \quad
\delta = \eps c \rho_1(w_0^{SN})\, , \quad
\mu_4 = \eps c \rho_1'(w_0^{SN}) \, ,
\]
and  $ \mu_j = \mu_j(\sigma^2),  \delta = \delta(\sigma^2)$.  By assuming the non-degeneracy condition $ \mu_3(\sigma^2) \mu_4(\sigma^2) \neq 0 $, performing a shift and scaling $ q = \tilde{q} - \frac{\delta(\sigma)}{\beta_4(0)} $ (assuming that $ \delta(0) \neq 0 $ for simplicity), $ w = \frac{\mu_3(\sigma^2)}{\mu_4(\sigma^2)^2} y, \tilde{X} = \left| \frac{\mu_3(\sigma^2)}{\mu_4(\sigma^2)} \right| X $, we bring \eqref{eq:slow_expanded} into the form
\[
\left\{
\begin{aligned}
y_{\tilde{X}} &= z\, ,\\
z_{\tilde{X}} & =  \beta_1 + \beta_2 y +  y^2 + s yz + \mathcal{G}(y,z). \, ,
\end{aligned}
\right.
\]
where
\[
\beta_1(\sigma^2) = \frac{\mu_4^4(\sigma^2)}{\mu_3^3(\sigma^2)} \mu_1(\sigma^2) \, , \qquad \beta_2(\sigma^2) = \left( \frac{\mu_4(\sigma^2)}{\mu_3(\sigma^2)} \right)^2 \mu_2(\sigma^2) \, ,
\]
and $s =  \mathrm{sign}(\beta_3(\sigma^2)\beta_4(\sigma^2))$. Hence, in order to conclude the corresponding scenario as described for $ \mathcal{S}_{BT} $ for our original system, it remains to analyze the mapping $ \sigma^2 \mapsto (\beta_1(\sigma^2),\beta_2(\sigma^2) ) $, which we refrain from doing here.

\end{document}